\documentclass[english]{bourbaki}

\date{Mars 2006}
\bbkannee{58{\`e}me ann{\'e}e, 2005-2006} \bbknumero{961}
\title{Rigidity results for Bernoulli actions and their \\ von Neumann algebras}
\subtitle{[after Sorin Popa]}

\author{Stefaan VAES}
\address{CNRS, Institut de Math.\ de Jussieu \\
175 rue du Chevaleret \\
F-75013 Paris (France) \vspace{2mm}\\ Department of Mathematics, K.U.Leuven \\
Celestijnenlaan 200B \\ B-3001 Leuven (Belgium)} \email{stefaan.vaes@wis.kuleuven.be}

\theoremstyle{plain}
\newtheorem{definition}{Definition}[section]
\newtheorem{proposition}[definition]{Proposition}
\newtheorem{theorem}[definition]{Theorem}
\newtheorem{corollary}[definition]{Corollary}
\newtheorem*{corolla}{Corollary}
\newtheorem*{theor}{Theorem}
\newtheorem{lemma}[definition]{Lemma}
\theoremstyle{definition}
\theoremstyle{remark}
\newtheorem{remark}[definition]{Remark}

\newtheorem{notation}[definition]{Notation}
\newtheorem{notations}[definition]{Notations}
\newtheorem{terminology}[definition]{Terminology}
\newtheorem{convention}[definition]{Convention}

\renewcommand{\Im}{\operatorname{Im}}

\newcommand{\ups}{\nu}

\newcommand{\cP}{\mathcal{P}}
\newcommand{\pr}{\operatorname{pr}}
\newcommand{\Ind}{\operatorname{Ind}}
\newcommand{\cI}{\mathcal{I}}
\newcommand{\cG}{\mathcal{G}}
\newcommand{\F}{\mathbb{F}}
\newcommand{\cartesian}[1]{\prod_{#1}}
\newcommand{\N}{\mathbb{N}}
\newcommand{\lspan}{\operatorname{span}}
\newcommand{\dis}{\displaystyle}
\newcommand{\Z}{\mathbb{Z}}
\newcommand{\SL}{\operatorname{SL}}
\newcommand{\om}{\omega}
\newcommand{\cN}{\mathcal{N}}
\newcommand{\cF}{\mathcal{F}}
\newcommand{\Mtil}{\widetilde{M}}
\newcommand{\cL}{\mathcal{L}}
\newcommand{\vphih}{\widehat{\varphi}}
\newcommand{\la}{\langle}
\newcommand{\ra}{\rangle}
\newcommand{\vphi}{\varphi}
\newcommand{\Ntil}{\widetilde{N}}
\newcommand{\Sp}{\operatorname{Sp}}
\newcommand{\Ptil}{\widetilde{P}}
\newcommand{\B}{\operatorname{B}}
\newcommand{\M}{\operatorname{M}}
\newcommand{\Tr}{\operatorname{Tr}}
\newcommand{\qtil}{\widetilde{q}}
\newcommand{\id}{\mathord{\text{\rm id}}}
\newcommand{\eps}{\varepsilon}
\newcommand{\si}{\sigma}
\newcommand{\recht}{\rightarrow}
\newcommand{\al}{\alpha}
\newcommand{\cU}{\mathcal{U}}
\newcommand{\Ad}{\operatorname{Ad}}
\newcommand{\R}{\mathbb{R}}
\newcommand{\ot}{\otimes}
\newcommand{\C}{\mathbb{C}}
\newcommand{\be}{\beta}

\newcommand{\vtil}{\widetilde{v}}
\newcommand{\tetil}{\widetilde{\theta}}
\newcommand{\Ker}{\operatorname{Ker}}
\newcommand{\Out}{\operatorname{Out}}
\newcommand{\cM}{\mathcal{M}}
\newcommand{\cA}{\mathcal{A}}
\newcommand{\tauh}{\widehat{\tau}}
\newcommand{\diag}{\operatorname{diag}}

\newcommand{\fun}{\mathcal{F}}
\newcommand{\ovt}{\overline{\otimes}}
\newcommand{\Char}{\operatorname{Char}}
\newcommand{\cZ}{\mathcal{Z}}
\newcommand{\actson}{\curvearrowright}
\newcommand{\cR}{\mathcal{R}}
\newcommand{\inside}{\preceq}
\newcommand{\Inn}{\operatorname{Inn}}
\newcommand{\Aut}{\operatorname{Aut}}
\newcommand{\oth}{\mathbin{\widehat{\otimes}}}
\newcommand{\Cliff}{\operatorname{Cliff}}
\newcommand{\T}{\mathbb{T}}

\newcommand{\og}{}
\newcommand{\fg}{}

\newcounter{teller}

\newenvironment{somoptel}{\begin{list}{(\arabic{teller})}{\usecounter{teller}\settowidth{\labelwidth}{(1)}\setlength{\leftmargin}{\labelwidth}\addtolength{\leftmargin}{\labelsep}\setlength{\rightmargin}{0cm}\setlength{\labelsep}{0.8ex}\addtolength{\leftmargin}{0.2ex}}}{\end{list}}

\begin{document}
\maketitle


\begin{center}
\parbox[t]{13cm}{\small \tableofcontents}
\end{center}

\mbox{}\vspace{-1cm}

\section{Introduction}

Suppose that a countable group $G$ acts freely and ergodically on the standard probability space $(X,\mu)$
preserving the probability measure $\mu$. We are interested in several
types of \lq isomorphisms\rq\ between such actions. Two actions are
said to be
\begin{enumerate}
\item \emph{conjugate} if there exists a group isomorphism and a
measure space isomorphism satisfying the obvious conjugacy formula;
\item \emph{orbit equivalent} if there exists a measure space
isomorphism sending orbits to orbits, i.e.\ the \emph{equivalence
relations} given by the orbits are isomorphic;
\item \emph{von Neumann equivalent} if the crossed product von
  Neumann algebras are isomorphic.
\end{enumerate}
Note that the \emph{crossed product} construction\footnote{The crossed product von Neumann algebra
  $L^\infty(X,\mu) \rtimes G$ contains a copy of $L^\infty(X,\mu)$ and
  a copy of the group $G$ by unitary elements in the algebra, and the
  commutation relations between both are given by the action of $G$ on
  $(X,\mu)$.} has been introduced by Murray and von Neumann \cite{MvN}, who
  called it the \emph{group measure space construction}.

It is clear that conjugacy of two actions implies orbit equivalence. Since the crossed product von Neumann algebra can be defined directly from the
equivalence relation given by the orbits, orbit equivalence implies von Neumann equivalence. \emph{Rigidity results} provide the converse
implications for certain actions of certain groups. This is a highly non-trivial matter. Dye \cite{Dy1,Dy2} proved that all free ergodic measure
preserving actions of groups with polynomial growth on the standard probability space are orbit equivalent. This result was extended to all
\emph{amenable groups} by Ornstein and Weiss \cite{OW}. Finally,
Connes, Feldman and Weiss \cite{CFW} showed that every ergodic amenable
probability measure preserving countable equivalence relation is
generated by a free $\Z$-action and is hence unique.
Summarizing, for amenable group actions all information on the group, except its
amenability, gets lost in the passage to the equivalence relation.

Concerning the relation between orbit equivalence and von Neumann equivalence, it was noted by Feldman
and Moore \cite{FM} that the \emph{pair} $L^\infty(X,\mu) \subset L^\infty(X,\mu) \rtimes G$ remembers the equivalence relation. The abelian
subalgebra $L^\infty(X,\mu)$ is a so-called \emph{Cartan subalgebra}. So, in order to deduce orbit equivalence from von Neumann equivalence, we need
certain uniqueness results for Cartan subalgebras, which is an extremely hard problem. Connes and Jones \cite{CJ2} gave the first examples of non
orbit equivalent, yet von Neumann equivalent actions.

In this talk, we discuss Popa's recent breakthrough rigidity results
for Bernoulli actions\footnote{Every discrete group $G$ acts on
  $(X,\mu) = \cartesian{g \in G} (X_0,\mu_0)$, by shifting the
  Cartesian product. Here
$(X_0,\mu_0)$ is the standard non-atomic probability space and the
action is called the Bernoulli action of $G$.} of Kazhdan groups. These results open a new era
in von Neumann algebra theory, with striking applications in ergodic
theory. The heart of Popa's work is his \emph{deformation/rigidity
  strategy}: he discovered families of von Neumann algebras with a rigid
subalgebra but yet with just enough deformation properties in order for the rigid part to be uniquely determined inside the ambient algebra (up to
unitary conjugacy). This leads to far reaching classification results for these families of von Neumann algebras. Popa considered the
deformation/rigidity strategy for the first time in \cite{P4}. In
\cite{P5}, he used it to deduce orbit equivalence from mere von
Neumann equivalence between certain group actions and to give the
first examples of II$_1$ factors with trivial fundamental group,
through an application of
Gaboriau's $\ell^2$ Betti numbers of equivalence relations \cite{Gab}.
Deformation/rigidity
arguments are again the crucial
ingredient in the papers \cite{P0,P1,P2,P3} that we discuss in this talk and they are applied in \cite{IPP}, in the study of amalgamated free products.
These ideas may lead to many more applications in von Neumann algebra and
ergodic theory (see e.g.\ the new papers
\cite{Ioana,PV} written since this talk was given). 

In the papers discussed in this talk, the \emph{rigidity} comes from the group side and is given by Kazhdan's property~(T) \cite{DK,Ka} and more
generally, by the relative property (T) of Kazhdan-Margulis (see \cite{dHV} and Valette's Bourbaki seminar \cite{VBour} for details): the groups
dealt with, contain an infinite normal subgroup with the relative
property (T) and are called \emph{$w$-rigid groups}. Popa discovered a strong
\emph{deformation property} shared by the Bernoulli actions, and called it \emph{malleability}. In a sense, a Bernoulli action can be continuously
deformed until it becomes orthogonal to its initial position. In order to exploit the tension between the deformation of the action and the rigidity
of the group, yet another technique comes in. Using \emph{bimodules} (Connes' correspondences), Popa developed a very strong method to prove that two
subalgebras of a von Neumann algebra are unitarily conjugate. Note that he used this bimodule technique in many different settings, see
\cite{IPP,OP,P1,P2,P5,P6}.

The following are the two main results of \cite{P0,P1,P2} and are discussed below.
The \emph{orbit equivalence superrigidity theorem}
states that the equivalence relation given by the orbits of a
Bernoulli action of a $w$-rigid group, entirely remembers the group
and the action. The \emph{von Neumann
strong rigidity theorem} roughly says that whenever a Bernoulli action
is von Neumann equivalent with a free ergodic action of a $w$-rigid
group, the actions are actually conjugate.
It is the first
theorem in the literature deducing conjugacy of actions out of von
Neumann equivalence. The methods and ideas behind
these far reaching results are fundamentally \emph{operator
  algebraic} and yield striking theorems in \emph{ergodic theory}.

\subsection*{Some important conventions}

All probability spaces in this talk are standard Borel spaces. All actions of countable groups $G$ on $(X,\mu)$ are supposed to preserve the
probability measure $\mu$. All statements about elements of $(X,\mu)$
only hold almost everywhere. A \emph{$w$-rigid group} is a countable group that admits an infinite normal subgroup with the relative property (T).

\subsection*{Orbit equivalence superrigidity}

In \cite{P0}, the deformation/rigidity technique leads to the following orbit equivalence superrigidity theorem.
\begin{theor}[Theorem \ref{thm.OEsuperrigidity}]
Let $G \actson (X,\mu)$ be the Bernoulli action of a $w$-rigid group $G$
as
above. Suppose that $G$ does not have finite normal subgroups.
If the restriction to $Y \subset X$ of the equivalence relation given
by $G \actson X$ is given by the orbits of $\Gamma \actson Y$ for some group $\Gamma$ acting freely and
ergodically on $Y$, then, up to measure zero, $Y=X$ and the actions of
$G$ and $\Gamma$ are conjugate through a group isomorphism.
\end{theor}
The theorem implies as well that the restriction to a Borel set of
measure $0 < \mu(Y) < 1$, of the Bernoulli action of a $w$-rigid group $G$ without finite normal subgroups, yields an ergodic
probability measure preserving countable equivalence relation that cannot be generated by a free action of a group. The first examples of this
phenomenon -- answering a question of Feldman and Moore -- were given by Furman in \cite{Fur2}. Dropping the ergodicity, examples were given before
by Adams in \cite{Adams}, who also provides examples in the Borel setting.

Popa proves the orbit equivalence superrigidity for the Bernoulli action of $G$ on $X$ using his even stronger \emph{cocycle superrigidity theorem}:
any $1$-cocycle for the action $G \actson X$ with values in a discrete group $\Gamma$ is cohomologous to a homomorphism of $G$ to $\Gamma$. The
origin of orbit equivalence rigidity and cocycle rigidity theory lies in Zimmer's pioneering work. Zimmer proved in \cite{Zim3} his celebrated
cocycle rigidity theorem and used it to obtain the first orbit equivalence rigidity results (see Section 5.2 in \cite{Zim}). Since Zimmer's theorem
deals with cocycles taking values in linear groups, he obtains orbit equivalence rigidity results where both groups are assumed to be linear (see
\cite{Zim2}). Furman developed in \cite{Fur1,Fur2} a new technique and obtains an orbit equivalence superrigidity theorem with quite general ergodic
actions of higher rank lattices on one side and an arbitrary free ergodic action on the other side. Note however that Furman's theorem nevertheless
depends on Zimmer's cocycle rigidity theorem. We also mention the orbit equivalence superrigidity theorems obtained by Monod and Shalom \cite{MS} for
certain actions of direct products of hyperbolic groups. An excellent overview of orbit equivalence rigidity theory can be found in Shalom's survey
\cite{Shal}.

Zimmer's cocycle rigidity theorem was a deep generalization of Margulis' seminal superrigidity theory \cite{Ma1}. In particular, the mathematics
behind involve the theory of algebraic groups and their lattices. On the other hand, Popa's technique to deal with $1$-cocycles for Bernoulli actions
is intrinsically operator algebraic.

As stated above, Popa uses his powerful \emph{deformation/rigidity strategy}
to prove the cocycle superrigidity theorem. Leaving aside several
delicate passages, the
argument goes as follows. A $1$-cocycle $\gamma$ for
the Bernoulli action $G \actson X$ of a $w$-rigid group $G$, can be
interpreted in two ways as a $1$-cocycle for the diagonal action $G
\actson X \times X$, either as $\gamma_1$, only depending on the first
variable, either as $\gamma_2$, only depending on
the second variable. The malleability of the Bernoulli action (this is
the deformation property) yields a continuous path joining $\gamma_1$
to $\gamma_2$. The relative property (T) implies that, in cohomology,
the $1$-cocycle remains essentially constant along the continuous
path. This yields $\gamma_1=\gamma_2$ in cohomology and the weak
mixing property allows to conclude that $\gamma$ is cohomologous to a homomorphism.

Let $(\si_g)$ be the Bernoulli action of a $w$-rigid group $G$ on $(X,\mu)$. Popa's cocycle superrigidity theorem covers his previous result
\cite{P4,PS} identifying the $1$-cohomology group $H^1(\si)$ with the character group $\Char G$. This result allows to compute as well the
$1$-cohomology for quotients of Bernoulli actions, yielding the
following result of \cite{P3}.

\begin{theor}[Theorem \ref{prop.non-orbit-equiv}]
Let $G$ be a $w$-rigid group. Then, $G$ admits a continuous family of non stably\footnote{See Definition \ref{def.stable-OE}.}  orbit equivalent
actions.
\end{theor}

Note that Popa does not only prove an existence result, but explicitly
exhibits a continuous family of mutually non orbit equivalent actions. The
existence of a continuum of non orbit equivalent actions of an
infinite property (T) group had been established before in a
non-constructive way by Hjorth \cite{Hjo}, who exhibits a continuous
family of actions such that every action in the family is orbit
equivalent to at most countably many other actions of the family.

Finally note that the first concrete computations of $1$-cohomology
for ergodic group actions are due to Moore \cite{Moore} and Gefter \cite{Gefter}.

\subsection*{Von Neumann strong rigidity}

The culmination of Popa's work on Bernoulli actions is the following \emph{von Neumann strong rigidity} theorem of \cite{P2}; it is the first theorem
in the literature that deduces conjugacy of the actions from isomorphism of the crossed product von Neumann algebras.
\begin{theor}[Theorem \ref{thm.strong-rigidity}]
Let $G$ be a group with infinite conjugacy classes and $G \actson
(X,\mu)$ its Bernoulli action as above. Let $\Gamma$ be a $w$-rigid group that
acts freely and ergodically on $(Y,\eta)$. If
$$\theta : L^\infty(Y) \rtimes \Gamma \recht p(L^\infty(X) \rtimes
G)p$$
is a $^*$-isomorphism
for some projection $p \in L^\infty(X) \rtimes G$, then $p=1$, the
groups $\Gamma$ and $G$ are isomorphic and the actions of $\Gamma$ and
$G$ are conjugate through this isomorphism.
\end{theor}
Note that in the conditions of the theorem, there is an assumption on the
action on one side and an assumption on the group on the other
side. As such, it is not a superrigidity theorem: one would like to
obtain the same conclusion for any free ergodic action of any group $\Gamma$
and for the Bernoulli action of a $w$-rigid ICC group $G$.

Another type of von Neumann rigidity has been obtained by Popa in \cite{P5,P6}, deducing orbit equivalence from von Neumann equivalence. We just
state the following particular case. Consider the usual action of $\SL(2,\Z)$ on $\T^2$. Whenever a free and ergodic action of a group $\Gamma$ with
the Haagerup property is von Neumann equivalent with the $\SL(2,\Z)$ action on $\T^2$, it actually is orbit equivalent with the latter. One should
not hope to deduce a strong rigidity result yielding conjugacy of the actions: Monod and Shalom (\cite{MS}, Theorem 2.27) proved that any free
ergodic action of the free group $\F_n$ is orbit equivalent with free ergodic actions of a continuum of non-isomorphic groups. Note that this also
follows from Dye's result \cite{Dy1,Dy2} if we assume that every generator of $\F_n$ acts ergodically.

\subsection*{II$_1$ factors and their fundamental group}

Let $G$ act freely and ergodically on $(X,\mu)$. Freeness and ergodicity imply that
  the crossed product von Neumann algebra
  $M:=L^\infty(X,\mu) \rtimes G$ is a \emph{factor} (the center of the
  algebra $M$ is reduced to the scalars) and the invariant probability
  measure yields a finite
  trace on $M$. Altogether, we get that $M$ is a so-called \emph{type II$_1$ factor}.

Another
  class of II$_1$ factors arises as follows: for any countable group
  $G$, one considers the von Neumann algebra\label{page.groupvnalg} $\cL(G)$ generated by the
  left translation operators on the Hilbert space $\ell^2(G)$. The algebra $\cL(G)$
  always admits a finite trace and it is a factor if and only if $G$
  has infinite conjugacy classes (ICC).

Let $M$ be a II$_1$ factor with normalized trace $\tau$. The \emph{fundamental group} of $M$, introduced by Murray and von Neumann \cite{MvN4}, is
the subgroup of $\R^*_+$ generated by the numbers $\tau(p)$, where $p$ runs through the projections of $M$ satisfying $M \cong pMp$. Murray and von
Neumann showed in \cite{MvN4} that the fundamental group of the
hyperfinite\footnote{The hyperfinite II$_1$ factor is, up to
  isomorphism, the unique II$_1$ factor that contains an increasing
  sequence of matrix algebras with weakly dense union.} II$_1$ factor is $\R^*_+$. They also write that there is no reason to
believe that the fundamental group of every II$_1$ factor is $\R^*_+$. However, only forty years later, this intuition was proved to be correct, in
a breakthrough paper of Connes \cite{C1}. Connes shows that the fundamental group of $\cL(G)$ is at most countable when $G$ is an ICC group with
Kazhdan's property (T). This can be considered as the first rigidity type result in the theory of von Neumann algebras. It was later refined by
Golodets and Nessonov \cite{GN} to obtain II$_1$ factors with countable fundamental group containing a prescribed countable subgroup of $\R^*_+$.
However, until Popa's breakthroughs in \cite{P1,P5,P6}, no precise computation of a fundamental group different from $\R^*_+$ had been obtained.

Note in passing that Voiculescu proved in \cite{Voi2} that the
fundamental group of $\cL(\F_\infty)$ contains the positive rationals
and that it was shown to be the whole of $\R^*_+$ by R{\u a}dulescu in \cite{Rad1}.
On the
other hand, computation of the fundamental group of
$\cL(\F_n)$ is equivalent with deciding on the (non)-isomorphism of the
free group factors (see \cite{Dyk,Rad2}), which is a famous open
problem in the subject.

Specializing the problem of Murray and von Neumann, Kadison \cite{K} posed the following question: does there exist a II$_1$ factor $M$ not
isomorphic to $\M_2(\C) \ot M$? This question was answered affirmatively by Popa in \cite{P5}, who showed that, among other examples, $\cL(G)$ has
trivial fundamental group when $G = \SL(2,\Z) \ltimes \Z^2$. For a
more elementary treatment of this example, see \cite{P6}. Note that
Popa shows in \cite{P5} that the fundamental group of
$\cL(G) = \SL(2,\Z) \ltimes L^\infty(\T^2)$ equals the fundamental
group of the equivalence relation given by the orbits of $\SL(2,\Z)
\actson \T^2$. The latter reduces to $1$ using Gaboriau's $\ell^2$ Betti number
invariants for equivalence relations, see \cite{Gab}.
We also refer to the Bourbaki seminar by Connes \cite{C2} on this part of Popa's oeuvre.

In \cite{P1}, Popa goes much further and constructs II$_1$ factors
with an arbitrary countable fundamental group!
\begin{theor}[Theorem \ref{thm.computation-fundamental}]
Given a countable subgroup $S \subset \R^*_+$ and a $w$-rigid ICC
group $G$ with $\cL(G)$ having trivial fundamental group, there
exists an action of $G$ on the hyperfinite II$_1$-factor $\cR$ such
that the crossed product $\cR \rtimes G$ is a II$_1$ factor with
fundamental group $S$.
\end{theor}
The example par excellence of a group $G$ satisfying the conditions of the theorem, is $G = \SL(2,\Z) \ltimes \Z^2$. Other examples include $\Gamma
\ltimes \Z^2$, where $\Gamma$ is any non-amenable subgroup of
$\SL(2,\Z)$ acting on $\Z^2$ by its given embedding in $\SL(2,\Z)$. Again, Popa does not establish a mere existence result: the actions
considered are the so-called Connes-St{\o}rmer Bernoulli actions (see
\cite{CS} and Section \ref{sec.bernoulli} below).

\subsection*{Some comments on proving von Neumann strong rigidity}

We explain how an isomorphism of crossed products forces, in certain
cases, actions to be conjugate.

In a first step, using the \emph{deformation/rigidity
  strategy}, Popa \cite{P1} shows the following result. Suppose that
  $G \actson (X,\mu)$ is the Bernoulli action of
an infinite
  group $G$ and consider the
  crossed product $L^\infty(X,\mu) \rtimes G$. It is shown (see
  Theorem \ref{thm.relative-rigid} below) that any subalgebra of
  $L^\infty(X,\mu) \rtimes G$ with the relative property (T) can
  essentially be unitarily conjugated into $\cL(G)$. Again leaving
  aside several delicate passages, the argument goes as follows.\label{page.explanation} A
  subalgebra $Q \subset L^\infty(X,\mu) \rtimes G$ with the relative
  property (T) is viewed in two ways as a subalgebra of $L^\infty(X
  \times X,\mu \times \mu) \rtimes G$, where $G$ acts diagonally: $Q_1$ only living
  on the first variable of $X \times X$ and $Q_2$ only living on the
  second one. The malleability of the Bernoulli action implies that the
  subalgebras $Q_1$ and $Q_2$ are joined by a continuous path of
  subalgebras $Q_t$. The relative property (T) then ensures that
  $Q_1$ and $Q_2$ are essentially unitarily conjugate. The mixing
  of the action is used to deduce that $Q$ can essentially be
  conjugated into $\cL(G)$.

Note in passing that above result remains true
when the \lq commutative\rq\ Bernoulli action is replaced by a \lq non-commutative\rq\ Connes-St{\o}rmer Bernoulli action, which is the crucial ingredient
to produce II$_1$ factors with prescribed countable fundamental
groups.

Given an isomorphism $\theta : L^\infty(Y) \rtimes \Gamma \recht
L^\infty(X) \rtimes G$, where $G \actson X$ is the Bernoulli action and
the group $\Gamma$ is $w$-rigid, the previous paragraph implies that $\theta$ sends $\cL(\Gamma)$ into $\cL(G)$, after conjugating by a unitary in
the crossed product. Using very precise analytic arguments, Popa \cite{P2} succeeds in proving next that also the Cartan subalgebras $L^\infty(Y)$
and $L^\infty(X)$ can be conjugated into each other with a unitary in
the crossed product (see Theorem \ref{thm.intertwine-cartan}
below). Having at
hand this orbit equivalence and knowing that the group von Neumann algebras can be conjugated into each other, Popa manages to prove conjugacy of the
actions.

An important remark should be made here. The results on Bernoulli actions discussed up to now, use the deformation property called \emph{strong
  malleability} combined with the mixing property of the action. So, they are
valid for all strongly malleable mixing actions. The result on the conjugation of the Cartan subalgebras however, uses a much stronger mixing
property of Bernoulli actions, called the \emph{clustering property}, which roughly means that the Bernoulli action allows for a natural \emph{tail}.
Note in this respect the following conjecture of Neshveyev and St{\o}rmer \cite{NS}: suppose that the abelian countable groups $G$ and $\Gamma$ act
freely and weakly mixingly on the standard probability space and that they give rise to isomorphic crossed products where the isomorphism sends
$\cL(G)$ onto $\cL(\Gamma)$; then, the Cartan subalgebras are conjugate with a unitary in the crossed product\footnote{It is crucial to have
conjugation of the Cartan subalgebras through a unitary in the crossed product, which is the hyperfinite II$_1$ factor. Indeed, thanks to the work of
Connes, Feldman and Weiss \cite{CFW}, two Cartan subalgebras are always conjugate with an automorphism of the hyperfinite II$_1$ factor. But, there exist
continuously many non inner conjugate Cartan subalgebras.}.

\subsection*{Outer conjugacy of actions on the hyperfinite II$_1$
  factor}

The deformation/rigidity technique first appeared\footnote{The paper
  \cite{P4} circulated since 2001 as a preprint of the MSRI and is the
  precursor of the papers \cite{P0,P1,P2,P3,PS} discussed above.} in Popa's paper
\cite{P4} on the computation of several invariants for (cocycle) actions
of $w$-rigid groups on the hyperfinite II$_1$ factor. In fact, many
  ideas exploited in the papers \cite{P0,P1,P2,P3,PS} are already
  present to some extent in the breakthrough paper \cite{P4}.

Recall that two
  actions $(\si_g)$ and $(\rho_g)$ of a group $G$ on a factor are said to be
  \emph{outer conjugate} if there exists an isomorphism $\Delta$ such
  that the conjugate automorphism $\Delta \si_g \Delta^{-1}$ equals
  $\rho_g$ up to an inner automorphism.

The classification up to outer conjugacy of actions of a group $G$ on, say, the
hyperfinite II$_1$ factor is an important subject. This classification
has been completed, first for cyclic groups by Connes \cite{C3,C4}, for
finite groups by Jones \cite{Jon1} and finally, for amenable groups by
Ocneanu \cite{Oc}: any two outer\footnote{An outer action is an action
  $(\si_g)$ such that for $g \neq e$, $\si_g$ is an outer
  automorphism, i.e.\ not of the form $\Ad u$ for a unitary $u$ in the von
  Neumann algebra.} actions of an amenable group $G$ on
the hyperfinite II$_1$ factor are outer conjugate (even cocycle
conjugate).

Away from amenable groups, Jones proved in \cite{Jon2} that any
non-amenable group admits at least two non outer conjugate actions on
the hyperfinite II$_1$ factor. Apart from actions, one also studies
\emph{cocycle actions} of a group $G$ on a factor $N$: families of
automorphisms $(\si_g)_{g \in G}$ such that $\si_g \si_h = \si_{gh}$
modulo an inner automorphism $\Ad u_{g,h}$, where the unitaries
$u_{g,h}$ satisfy a $2$-cocycle relation.

In the previously cited works on amenable group actions, it is shown as well that any cocycle action of an amenable group on the hyperfinite II$_1$
factor is outer conjugate to a genuine action. Popa generalized this result to arbitrary II$_1$ factors in \cite{P7}. In \cite{CJ}, Connes and Jones
constructed, for any infinite property (T) group $G$, examples of cocycle actions of $G$ on the free group factor $\cL(\F_\infty)$ that are non outer
conjugate to a genuine action.

This brings us to the topic of \cite{P4}. Popa introduces two outer conjugacy invariants for a (cocycle) action on a II$_1$ factor: the fundamental
group and the spectrum. These invariants are computed in \cite{P4} for the Connes-St{\o}rmer Bernoulli actions, yielding the following theorem.

\begin{theor}[Theorems \ref{thm.fundamental-hyperfinite} and \ref{thm.not-conjugate}]
Let $G$ be a $w$-rigid group. Then $G$ admits a continuous family of
non outer conjugate actions on the hyperfinite II$_1$ factor. Also,
$G$ admits a continuous family of cocycle actions on the hyperfinite
II$_1$ factor that are non outer conjugate to a genuine action.
\end{theor}

\subsection*{Further remarks}

We discussed in detail how Popa recovers information on a group action
from the crossed product algebra $L^\infty(X,\mu) \rtimes G$. On the
other hand, to
what extent a group von Neumann algebra
$\cL(G)$ remembers the group
$G$? Very little is known on this problem. Connes' celebrated
theorem \cite{C6} states that all the II$_1$ factors $\cL(G)$ defined by amenable
ICC groups $G$ are isomorphic to the hyperfinite II$_1$
  factor. Indeed, they are all \emph{injective}\label{page.conditional}\footnote{A factor $M
    \subset \B(H)$ is called injective if there exists a conditional
    expectation of $\B(H)$ onto $M$ (which of course need not be weakly
    continuous). A conditional expectation of a von Neumann $M$ onto a
  von Neumann subalgebra $N$ is a unital, positive, $N$-$N$-bimodule map $E : M \recht
  N$.} and Connes shows in
\cite{C6} the uniqueness of the injective II$_1$ factor.  Cowling and Haagerup \cite{CH} have shown that the group von
Neumann algebras $\cL(\Gamma)$ are non-isomorphic if one takes
lattices $\Gamma$ in $\operatorname{Sp}(1,n)$ for different values of $n$.

Some group von Neumann algebras $\cL(G)$ can we written
as well as the crossed product by a free ergodic action (but not all,
since Voiculescu \cite{Voi1} showed that the free group factors cannot be written
in this way). We have for instance
$\cL(\SL(n,\Z) \ltimes \Z^n) = L^\infty(\T^n) \rtimes \SL(n,\Z)$.
Another example consists in writing the Bernoulli action crossed
product $L^\infty(X,\mu) \rtimes G$ as $\cL(\Z \wr G)$, where the
\emph{wreath product} group $\Z \wr G$ is defined as the semidirect product $\Z \wr G := (\bigoplus_{g \in G} \Z) \rtimes G$.
Popa's von Neumann
strong rigidity theorem then implies the following result. It can be
considered as a relative version of
Connes' conjecture \cite{C5}, which states that within the class of
ICC property (T) groups, $\cL(G_1) \cong \cL(G_2)$ if and only if $G_1
\cong G_2$. Popa's result \lq embeds injectively\rq\ the category of
$w$-rigid ICC groups into the category of II$_1$ factors.

\begin{corolla}
When $G$ and $\Gamma$ are $w$-rigid ICC groups, $\cL(\Z \wr G) \cong
\cL(\Z \wr \Gamma)$ if and only if $G \cong \Gamma$.
Moreover,
$\cL(\Z \wr G)$ has trivial fundamental group for any $w$-rigid ICC group $G$.
\end{corolla}

Popa's von Neumann strong rigidity theorem is in fact
more precise than the version stated above. As we shall see in Theorem
\ref{thm.strong-rigidity} below, the strong rigidity theorem allows as
well to compute the group $\Out M$ of outer automorphisms of $M =
L^\infty(X,\mu) \rtimes G$, where $G$ is a $w$-rigid ICC group and $G
\actson (X,\mu)$ its Bernoulli action.
Then, $$\Out M \cong \Char G \rtimes
\frac{\Aut^*(X,G)}{G} \; ,$$ where $\Aut^*(X,G)$ is the group of
measure space isomorphisms $\Delta : X \recht X$ for which there
exists a $\delta \in \Aut G$ such that $\Delta(g \cdot x) = \delta(g)
\cdot \Delta(x)$ almost everywhere. Writing $\Delta_g(x) = g \cdot x$,
one embeds $G \hookrightarrow \Aut^*(X,G)$.
Note moreover that
$\Aut^*(X,G)$ obviously contains another copy of $G$ acting by Bernoulli
shifts \lq on the other side\rq.

In \cite{IPP}, Ioana, Peterson and Popa apply the strategy of deformation/rigidity in the
completely different context of amalgamated free products, yielding
the first examples of II$_1$ factors with trivial outer automorphism
group. Much more is done in \cite{IPP}, where actually a von Neumann
version of Bass-Serre theory is developed.

\subsection*{Acknowledgment}

\noindent {\it I learned about Popa's work during his numerous lectures in Paris
  and his course in the Coll{\`e}ge de France in the fall of 2004. It is
  my pleasure to express my warmest thanks to Sorin Popa for all our
  discussions and for his
  hospitality at UCLA. I thank Georges Skandalis for his
  numerous remarks and the entire days we spent discussing the present
  text. I finally took benefit from
  the remarks of Aur{\'e}lien Alvarez, Claire Anantharaman, \linebreak Etienne
  Blanchard, Alain Connes, David Fisher, Damien Gaboriau, Sergey
  Neshveyev, Mikael Pichot, Alain Valette and Dan Voiculescu. Thanks
  to all of you!}


\section{Preliminaries and conventions}

\subsection*{Von Neumann algebras, traces, almost periodic states and group actions}
Throughout $M,\cM,N,\cN,A,\cA$ denote \emph{von Neumann algebras}. Recall
that a von Neumann algebra is a non-commutative generalization of a
measure space, the algebras $L^\infty(X,\mu)$ being the abelian
examples. By definition, a von Neumann algebra is a weakly closed
unital $^*$-subalgebra of $\B(H)$ for some Hilbert space $H$. Whenever $\cM
\subset \B(H)$ is a von Neumann algebra, the \emph{commutant} of $\cM$
is denoted by $\cM'$ and consists of the operators in $\B(H)$
commuting with all the operators in $\cM$. Von Neumann's
\emph{bicommutant} theorem states that $\cM^{\prime\prime} = \cM$ and
this equality characterizes von Neumann algebras among the unital
$^*$-subalgebras of $\B(H)$. A \emph{factor} is a von Neumann algebra
with trivial center, i.e.\ $\cM \cap \cM' = \C 1$.

A \emph{state} on a von Neumann algebra is a positive linear map $\cM
\recht \C$ satisfying $\om(1) = 1$. All states are assumed to be
\emph{normal}, i.e.\ continuous with respect to the ultraweak topology
on $\cM$ (which is
  equivalent with requiring weak continuity on the unit ball of $\cM$). Hence, normal
states are the counterparts of
probability measures on $(X,\mu)$ absolutely continuous with respect to
$\mu$. A state $\om$ is said to be \emph{tracial} if $\om(xy) =
\om(yx)$ for all $x,y$. A state is said to be \emph{faithful} if the
equality $\om(x)=0$ for $x$ positive implies that $x=0$. States are
always assumed to be faithful.

The algebras denoted $M,N,A$ are supposed to admit a \emph{faithful
  normal trace} and if we specify a state
on $M,N$ or $A$, it is always supposed to be a trace. The terminology
  \emph{finite von Neumann algebra $(N,\tau)$} means a von Neumann
  algebra $N$ with a faithful normal trace $\tau$.

An action of a countable group on $(\cM,\vphi)$ is understood to be an
action by
automorphisms \emph{leaving the state $\vphi$ invariant}. We denote by $(X,\mu)$ the
standard probability space without atoms and an action of a countable
group on $(X,\mu)$ is supposed to preserve the probability measure
$\mu$.

If $G$ acts on $(\cM,\vphi)$ by automorphisms $(\si_g)$, $\cM^\si$ denotes the von Neumann
subalgebra of elements $x \in \cM$ satisfying $\si_g(x) = x$ for all
$g \in G$. The action $(\si_g)$ is said to be \emph{ergodic}
if $\cM^\si=\C 1$.

If $\vphi$ is a faithful normal state on $\cM$, we consider the
\emph{centralizer algebra} $\cM^\vphi$ of $\vphi$ consisting of those
$x \in \cM$ satisfying $\vphi(xy) = \vphi(yx)$ for all $y$. More
generally, for a real number $\lambda > 0$, a \emph{$\lambda$-eigenvector} for $\vphi$ is an element $x \in \cM$
satisfying $\vphi(xy) = \lambda \vphi(yx)$ for all $y \in
\cM$. We say that $\vphi$ is
\emph{almost periodic} (or that $(\cM,\vphi)$ is almost periodic), if
the $\lambda$-eigenvectors span a weakly dense subalgebra of
$\cM$ when $\lambda$ runs through $\R^*_+$. If this is the case, $\Sp(\cM,\vphi)$ denotes the point
spectrum of $\vphi$, i.e. the set of $\lambda > 0$ for which there
exists a non-zero $\lambda$-eigenvector.

A finite von Neumann algebra $(P,\tau)$ is said to be \emph{diffuse}
if $P$ does not contain a minimal projection. A finite $(P,\tau)$ is
diffuse if and only if $P$
contains a sequence of unitaries tending weakly to zero. Equivalently,
$P$ does not have a direct summand that is a matrix algebra. For
instance, the group von Neumann algebra $\cL(G)$ (see page
\pageref{page.groupvnalg} for its definition) is diffuse for any
infinite group $G$.

\subsection*{Crossed products}
Whenever a countable group $G$ acts by $\vphi$-preserving automorphisms $(\si_g)$ on $(\cM,\vphi)$, we denote by $\cM \rtimes G$ the crossed product,
generated by the elements $a \in \cM$ and the unitaries $(u_g)_{g \in G}$ such that $u_g a u_g^* = \si_g(a)$ for all $a \in \cM$ and $g \in G$. We
have a natural conditional expectation (see footnote on page
\pageref{page.conditional}) given by $E : \cM \rtimes G \recht \cM :
E(a u_g) = \delta_{g,e} \, a$ and we extend $\vphi$ to a faithful normal state on
$\cM \rtimes G$ by the formula $\vphi \circ E$. If $\vphi$ is tracial, its extension is tracial.

The crossed product $M$ is a factor (hence, a type II$_1$ factor) in the following (non-exhaustive) list of
examples. If $A \subset M$ is an inclusion of von Neumann algebras, we
denote by $M \cap A'$ the \emph{relative commutant} consisting of
elements in $M$ commuting with all elements of $A$.
\begin{itemize}
\item Suppose that $G$ acts (essentially) freely on
  $(X,\mu)$ and put $M = L^\infty(X) \rtimes G$. Then, $M \cap
  L^\infty(X)' = L^\infty(X)$ and $M$ is a factor if and only if the
  $G$-action is ergodic.
\item Suppose that the ICC group $G$ acts on the finite $(N,\tau)$ and
  put $M = N \rtimes G$. Then, $M \cap \cL(G)' = N^G$ and $M$ is a
  factor if and only if the $G$-action on the center of $N$ is ergodic.
\item Suppose that the group $G$ acts on the II$_1$ factor $(N,\tau)$ such that for
  all $g \neq e$, $\si_g$ is an outer automorphism of $N$, i.e.\ an
  automorphism that cannot be written as $\Ad u$ for some unitary $u
  \in N$. Putting $M = N \rtimes G$, we have $M \cap N' = \C 1$ and in
  particular, $M$ is a factor.
\end{itemize}

\subsection*{$1$-cocycles and $1$-cohomology}
Let the countable group $G$ act on $(X,\mu)$. We denote by $g \cdot x$
the action of an element $g \in G$ on $x \in X$ and we denote by
$(\si_g)$ the corresponding action of $G$ on $A=L^\infty(X)$ given by
$(\si_g (F))(x) = F(g^{-1} \cdot x)$. A
\emph{$1$-cocycle} for $(\si_g)$ with coefficients in a Polish group
$K$ is a measurable map
$$\gamma : G \times X \recht K \quad\text{satisfying}\quad \gamma(gh,x) = \gamma(g,h \cdot x) \; \gamma(h,x)$$
almost everywhere. Two $1$-cocycles $\gamma_1$ and $\gamma_2$ are said
to be \emph{cohomologous} if there exists a measurable map $w : X
\recht K$ such that
$$\gamma_1(g,x) = w(g \cdot x) \gamma_2(g,x) w(x)^{-1}
\quad\text{almost everywhere}.$$
Whenever $K$ is abelian, the $1$-cocycles form a group $Z^1(\si,K)$
and quotienting by the $1$-cocycles cohomologous to the trivial $1$-cocycle, we
obtain $H^1(\si,K)$.
Whenever $K=S^1$, we just write $Z^1(\si)$ and $H^1(\si)$.
Several important remarks should be made. Suppose that the action of $G$ on $(X,\mu)$ is free and ergodic.
\begin{itemize}
\item Write $M = L^\infty(X) \rtimes G$. The group $Z^1(\si)$ embeds in $\Aut(M)$, associating with $\gamma \in Z^1(\si)$, the automorphism $\theta_\gamma$ of $M$ defined by
$\theta_\gamma(a) = a$ for all $a \in L^\infty(X)$ and $\theta_\gamma(u_g) = u_g
\gamma(g, \cdot)$. Passing to quotients, $H^1(\si)$ embeds into $\Out(M)$.
\item $H^1(\si)$ is an invariant for $(\si_g)$ up to stable orbit
  equivalence (see Definition \ref{def.stable-OE}).
\item If $(\si_g)$ is weakly mixing, the group of characters $\Char G$ embeds into $H^1(\si)$ as $1$-cocycles not depending on the space variable $x$.
\end{itemize}

\subsection*{The fundamental group of a II$_1$ factor}
Let $M$ be a II$_1$ factor. If $t > 0$, we define, up to isomorphism,
the \emph{amplification} $M^t$ as follows: choose $n \geq 1$ and a projection
$p \in \M_n(\C) \ot M$ with $(\Tr \ot \tau)(p) = t$. Define $M^t :=
p(\M_n(\C) \ot M)p$. The \emph{fundamental group} of $M$ is defined as
$$\fun(M) = \{ t > 0 \mid M^t \cong M \} \; .$$
It can be checked that $\fun(M)$ is a subgroup of $\R^*_+$.

In Theorem \ref{thm.strong-rigidity}, a large class of non-isomorphic II$_1$ factors with trivial fundamental group are obtained. In Theorem
\ref{thm.computation-fundamental}, II$_1$ factors with a prescribed countable subgroup of $\R^*_+$ as a fundamental group, are constructed.

\subsection*{Quasi-normalizers and almost normal subgroups}

Let $Q \subset M$ be a von Neumann subalgebra of $M$. An element $x
\in M$ is said to \emph{quasi-normalize} $Q$ if there exist
$x_1,\ldots,x_k$ and $y_1,\ldots,y_r$ in $M$ such that
$$x Q \subset \sum_{i=1}^k Q x_i \quad\text{and}\quad Q x \subset
\sum_{i=1}^r y_i Q \; .$$
The elements quasi-normalizing $Q$ form a $^*$-subalgebra of $M$ and
their weak closure is called the \emph{quasi-normalizer} of $Q$ in
$M$. The inclusion $Q \subset M$ is said to be \emph{quasi-normal} if
$M$ is the quasi-normalizer of $Q$ in $M$.

A typical example arises as follows: let $G$ be a countable group and
$H$ an \emph{almost normal subgroup}, which means that $g H g^{-1}
\cap H$ is a finite index subgroup of $H$ for every $g \in
G$. Equivalently, this means that for any $g$ in $G$, $H g H$ is
the union of finitely many left cosets, as well as the union of
finitely many right cosets. So, it is clear that for every almost
normal subgroup $H \subset G$, the inclusion $\cL(H) \subset \cL(G)$
is quasi-normal.

There are some advantages to work with the quasi-normalizer rather
than the normalizer. In Lemma \ref{lem.quasi-normalizer}, the
following is shown: let $Q \subset M$ be an inclusion of finite von Neumann
algebras and let $p$ be a projection in $Q$. If $P$
denotes the quasi-normalizer of $Q$ in $M$, the quasi-normalizer of
$pQp$ in $pMp$ is $pPp$. This is no longer true for the
actual normalizer.

\bigskip

\noindent{\bf More background material} is available in the
appendices. We discuss in Appendix~\ref{sec.basic} the \emph{basic construction} $\la \cN,e_B \ra$ starting from an inclusion $B \subset \cN$ of a
von Neumann algebra $B$ in the centralizer algebra of $(\cN,\vphi)$ (in particular, for an inclusion of finite von Neumann algebras). Appendix~\ref{sec.relativeT} deals with the relative property (T) and its analogue for inclusions of finite von Neumann algebras. In Appendix~\ref{sec.intertwining} is studied the relation between conjugating von Neumann subalgebras with a unitary and the existence of finite-trace
bimodules. Finally, Appendix~\ref{sec.mixing} is devoted to (weakly) mixing actions.

\section{The malleability property of Bernoulli actions} \label{sec.bernoulli}

Popa discovered several remarkable properties of Bernoulli actions. The first one is a deformation property, that he called strong malleability and
that is discussed in this section. This notion of malleability, together with its stunning applications, should be considered as one of the major
innovations of Popa.

As is well known, the Bernoulli actions are mixing (see Appendix \ref{sec.mixing} for definition and results) and this fact is used throughout. But,
Popa exploits as well a very strong mixing property of Bernoulli actions that he called the \emph{clustering property}. This will be used in Section
\ref{sec.Neumann-to-orbit}.

\begin{definition}[Popa, \cite{P1,PS}] \label{def.malleable}
The action $(\si_g)$ of $G$ on $(\cN,\vphi)$ is said to be
\begin{itemize}
\item \emph{malleable} if there exists a continuous action
$(\al_t)$ of $\R$ on $(\cN \ot \cN,\vphi \ot \vphi)$ that commutes with the
diagonal action $(\si_g \ot \si_g)$ and satisfies $\al_1(a \ot 1) = 1
\ot a$ for all $a \in \cN$;
\item \emph{strongly malleable} if there moreover exists an
  automorphism $\be$ of $(\cN \ot \cN,\vphi \ot \vphi)$ commuting with $(\si_g \ot \si_g)$ such that
  $\be \al_t = \al_{-t} \be$ for all $t \in \R$ and $\be(a \ot 1) = a
  \ot 1$ for all $a \in \cN$ and such that $\be$ has period $2$:
  $\be^2 = \id$.
\end{itemize}
\end{definition}

\begin{remark}
In \cite{P1,P2}, Popa uses the term \lq malleability\rq\ for a
 larger class of actions: indeed, instead of extending the action from
 $\cN$ to $\cN \ot \cN$, he allows for a more general extension to
 $\widetilde{\cN}$, which can typically be a graded tensor square $\cN
 \widehat{\otimes} \cN$. This last example occurs when considering
 Bogolyubov actions. See remark \ref{rem.graded-malleable} for details.
\end{remark}

\subsection*{Generalized Bernoulli actions}
The main example of a \emph{strongly malleable action} arises as a (generalized) Bernoulli action. Let $G$ be a countable group that acts on the
countable set $I$. Let $(X_0,\mu_0)$ be a probability space. The action of $G$ on $(X,\mu):=\cartesian{i \in I} (X_0,\mu_0)$ by shifting the infinite
product, is called the (generalized) Bernoulli action. The usual Bernoulli action arises by taking $I = G$ with the action of $G$ by translation.
\begin{convention}
For simplicity, we only deal with Bernoulli actions on the infinite product of \emph{non-atomic} probability spaces and we refer to them as \emph{Bernoulli actions with non-atomic base}. Most of Popa's results also hold for Bernoulli actions on products of atomic spaces. They are no longer malleable
but \emph{sub-malleable}, see Definition 4.2 in \cite{P1} and Remark \ref{rem.atomic-bernoulli}.
\end{convention}
Write $A_0 = L^\infty(\R/\Z)$. To check that the generalized Bernoulli action is strongly malleable, it suffices to produce an action $(\al_t)$ of
$\R$ on $A_0 \ot A_0$ and a period $2$ automorphism $\be$ of $A_0 \ot A_0$ such that $\al_1(a \ot 1) = 1 \ot a$, $\be(a \ot 1) = a \ot 1$ for all $a
\in A_0$ and $\be \al_t = \al_{-t} \be$ for all $t \in \R$. One can then take the infinite product of these $(\al_t)$ and $\be$. Take the uniquely determined map $f :
\R / \Z \recht \,\bigl]-\frac{1}{2},\frac{1}{2}\bigr]$ satisfying $x = f(x) \mod \Z$ for all $x$. Define the measure preserving flow $\al_t$ and the
measure preserving transformation $\be$ on $\R/\Z \times \R/\Z$ by the formulae
$$\al_t (x,y) = (x+tf(y-x),y+tf(y-x)) \quad\text{and}\quad \be(x,y) =
(x,2x-y) \; .$$
For $F \in L^\infty(\R/\Z \times \R/\Z)$, write $\al_t (F) = F \circ
\al_t$ and $\be(F) = F \circ \be$.

Popa gives a more functional analytic argument for the strong malleability of the generalized Bernoulli action. Consider $A_0 \ot A_0$ as being
generated by two independent Haar unitaries $u$ and $v$. We have to construct a one-parameter group $(\al_t)$ and a period $2$ automorphism $\beta$
such that $\al_1(u) = v$, $\be(u) = u$ and $\be \al_t = \al_{-t} \be$. Conjugating $\al_t$ and $\be$ with the automorphism $\si$ determined by
$\si(u) = u, \si(v) = vu$ (note that $u$ and $vu$ are independent generating Haar unitaries), the first requirement changes to $\al_1(u) = vu$ and
the other requirements remain. Taking $\log : \T \recht \; ]\! -\pi,\pi]$, we can now set $\al_t(u) = \exp(t \log v) u$, $\al_t(v) = v$ and $\be(u) =
u$, $\be(v) = v^*$.

\subsection*{Connes-St{\o}rmer Bernoulli actions}
Apart from \lq classical\rq\ Bernoulli actions, also the \lq non-commutative\rq\ Bernoulli actions of Connes and St{\o}rmer \cite{CS} satisfy Popa's
malleability condition. These Connes-St{\o}rmer Bernoulli actions provide the main non-commutative examples of \emph{malleable} actions.

Let $G$ be a countable group acting on a countable set $I$. Let $\vphi_0$ be a faithful normal state on $\B(H)$ for some Hilbert space $H$ (finite or
infinite-dimensional). Define
$$(\cN,\vphi) := \bigotimes_{i \in I} (\B(H),\vphi_0) \; .$$
On $(\cN,\vphi)$, $G$ acts by shifting the tensor factors. To prove the malleability, one has to produce an action $(\al_t)$ of $\R$ on $(\B(H \ot
H),\vphi_0 \ot \vphi_0)$ satisfying $\al_1(a \ot 1) = 1 \ot a$ for all $a \in \B(H)$. Denoting by $P \in \B(H \ot H)$ the orthogonal projection on
the symmetric subspace densely spanned by the vectors $\xi \ot \mu + \mu \ot \xi$ for $\xi,\mu \in H$, we define $U_t = P + e^{i \pi t} (1-P)$ and
$\al_t = \Ad U_t$. Note that Connes-St{\o}rmer Bernoulli actions are not in an obvious way \emph{strongly} malleable. In some cases however, a
generalization of strong malleability holds, see \ref{rem.graded-malleable}.

The state $\vphi_0$ is of the form $\Tr_\Delta$ for some positive
trace-class operator $\Delta$. So, $\vphi$ is almost periodic and
$\Sp(\cN,\vphi)$ is the subgroup of $\R^*_+$ generated by the ratios
$t/s$, where $t,s$ belong to the point spectrum of $\Delta$.

\section[Cocycle and orbit
equivalence superrigidity for Bernoulli actions]{Superrigidity for Bernoulli actions}

In this section, Popa's very strong rigidity results for Bernoulli actions of $w$-rigid groups are proved: according to the philosophy in the
beginning of the introduction, an orbit equivalence rigidity result deduces conjugacy of actions out of their mere orbit equivalence. All these
rigidity result follow from the following cocycle superrigidity theorem.

\begin{theorem}[Popa, \cite{P0}] \label{thm.cocycle}
Let $G$ be a countable group with infinite normal subgroup $H$ such that $(G,H)$ has the relative property (T). Let $G$ act strongly malleably on
$(X,\mu)$ and suppose that its restriction to $H$ is weakly mixing. Then, any $1$-cocycle
$$\gamma : G \times X \recht K$$
with values in a closed subgroup $K$ of the unitary group $\cU(B)$ of a finite von Neumann algebra $(B,\tau)$, is cohomologous to a homomorphism
$\theta : G \recht K$.
\end{theorem}

By regarding $\Gamma \subset \cU(\cL(\Gamma))$, the theorem covers all $1$-cocycles with values in countable groups, which is the crucial ingredient
to prove orbit equivalence rigidity results.

The superrigidity theorem for Bernoulli actions proved below, does not only deal with orbit equivalence, but also with \emph{stable orbit
equivalence}. There are several ways to introduce this concept, one of them being the following (see e.g.\ \cite{Fur2}, where the terminology of weak
orbit equivalence is used).

\begin{definition} \label{def.stable-OE}
Let $G \actson (X,\mu)$ and $\Gamma \actson (Y,\eta)$ be free and ergodic actions. A \emph{stable orbit equivalence} between these actions is given
by a measure space isomorphism $\pi : A \recht B$ between non-negligible subsets $A \subset X$ and $B \subset Y$ preserving the restricted
equivalence relations: $\pi(A \cap (G \cdot x)) = B \cap (\Gamma \cdot \pi(x))$ for almost all $x \in A$.

The \emph{compression constant} of  $\pi$ is
defined as $c(\pi):=\eta(B)/\mu(A)$.

The maps $\pi_i : A_i \recht B_i$ ($i=1,2$) define \emph{the same
  stable orbit equivalence} if $$\pi_2(A_2 \cap (G \cdot x)) \subset
  \Gamma \cdot \pi_1(x) \quad\text{for almost
all}\;\; x \in A_1 \; .$$ Note that this implies that $c(\pi_1)=c(\pi_2)$.
\end{definition}

Suppose that $\pi_i : A_i \recht B_i$ ($i=1,2$) define the same stable orbit equivalence. If, say, $\mu(A_1) \leq \mu(A_2)$, there exist $\phi$ in
the full group\label{note.fullgroup}\footnote{The full group of the equivalence relation defined by $G$-orbits, consists of
  the measure space automorphisms $\Delta : X \recht X$ satisfying $\Delta(x)
\in G \cdot x$ for almost all $x$.} of the equivalence
relation given by the $G$-orbits and $\psi$ in the full group of the
equivalence relation given by the $\Gamma$-orbits such that $\phi(A_1)
\subset A_2$ and $\pi_1$ is the
restriction of $\psi \circ \pi_2 \circ \phi$ to $A_1$.

If $\pi : A \recht B$ defines a stable orbit equivalence between the free and ergodic actions $G \actson (X,\mu)$ and $\Gamma \actson (Y,\eta)$, one
defines as follows a $1$-cocycle $\al : G \times X \recht \Gamma$ for
$G \actson X$ with values in $\Gamma$. By ergodicity, we can
choose a measurable map $\pr_A : X
\recht A$ satisfying $\pr_A(x) \in G \cdot x$ almost everywhere and denote $p = \pi \circ \pr_A$. Freeness of the action $\Gamma \actson Y$, allows
to define
$$\al : G \times X \recht \Gamma : p(g \cdot x) = \al(g,x) \cdot p(x)$$
almost everywhere. Taking another $\pi$ defining the same stable orbit equivalence or choosing another $\pr_A$, yields a cohomologous $1$-cocycle.

Given a free and ergodic action $G \actson (X,\mu)$, there are certain actions that are trivially stably orbit equivalent to $G \actson X$ and we
introduce them in Notation \ref{not.induced}. The superrigidity theorem \ref{thm.OEsuperrigidity} states that for Bernoulli actions of $w$-rigid
groups these are \emph{the only actions} that are stably orbit equivalent to the given Bernoulli action.

\begin{notation} \label{not.induced}
Let $G$ act freely and ergodically on $(X,\mu)$. Suppose that $\theta
: G \recht \Gamma$ is a homomorphism with $\Ker \theta$ finite 
and $\Im \theta$ of finite index in $\Gamma$. Define
$$\Ind_G^\Gamma(X,\theta) := G \backslash (X \times \Gamma)
\quad\text{where $G$ acts on $X \times \Gamma$ by}\quad g \cdot (x,s)
= (g \cdot x,\theta(g)s) \; .$$
The action of $\Gamma$ on $\Ind_G^\Gamma(X,\theta)$ given by $t \cdot
(x,s) = (x,st^{-1})$ is free,
ergodic and finite measure preserving.
We also have a canonical stable orbit equivalence between $G \actson X$ and $\Gamma \actson \Ind_G^\Gamma(X,\theta)$, with compression constant
$[\Gamma : \theta(G)]/|\Ker \theta|$.
\end{notation}

\begin{theorem}[Popa, \cite{P0}] \label{thm.OEsuperrigidity}
Let $G$ be a countable group with infinite normal subgroup $H$ such that $(G,H)$ has the relative property (T). Let $G$ act strongly malleably on
$(X,\mu)$ and suppose that its restriction to $H$ is weakly
mixing.

Whenever $\Gamma$ is a countable group acting freely and ergodically on $(Y,\eta)$ and whenever $\pi$ defines a stable orbit equivalence between $G
\actson X$ and $\Gamma \actson Y$, there exists
\begin{itemize}
\item a homomorphism $\theta : G \recht \Gamma$ with $\Ker \theta$
  finite in $G$ and $\Im \theta$ of finite index in $\Gamma$~;
\item a measure space isomorphism $\Delta : Y \recht
  \Ind_G^\Gamma(X,\theta)$ conjugating the actions $\Gamma \actson Y$
  and $\Gamma \actson \Ind_G^\Gamma(X,\theta)$~;
\end{itemize}
such that $\Delta \circ \pi$ defines the canonical stable orbit equivalence between $G \actson X$ and $\Gamma \actson \Ind_G^\Gamma(X,\theta)$. In
particular, the compression constant $c(\pi)$ equals $[\Gamma:\theta(G)]/|\Ker \theta|$.
\end{theorem}

\begin{remark}
Several special instances of Theorem \ref{thm.OEsuperrigidity} should be mentioned. Suppose that the action $G \actson X$ satisfies the conditions of
Theorem \ref{thm.OEsuperrigidity} and denote by $\cR$ the equivalence relation given by the $G$-orbits.
\begin{itemize}
\item If we suppose moreover that $G$ \emph{does not have finite
    normal subgroups}, we get the following result stated in the introduction.
If the restriction to $Y \subset X$ of the equivalence relation given
    by $G \actson X$ is given by the orbits of $\Gamma \actson Y$ for some group
$\Gamma$ acting freely and ergodically on $Y$, then, up to measure zero, $Y=X$ and the actions of $G$ and $\Gamma$ are conjugate through a group
isomorphism.
\item The amplified equivalence relation\footnote{The
  amplified equivalence relation $\cR^t$ is defined as follows. If $t
  \leq 1$, we restrict $\cR$ to a subset of measure $t$. If $t > 1$,
  we take a restriction of the obvious type II$_1$ equivalence
  relation on $X \times \{1,\ldots,n\}$.} $\cR^t$ can be
  generated by a free action of a group if and only if $t = n/|G_0|$,
  where $n \in \N \setminus \{0\}$ and $G_0$ is a finite normal subgroup
  of $G$. So, we get many examples of type II$_1$ equivalence
  relations that \emph{cannot be generated by a free action of a
  group}. The first such examples were given by Furman \cite{Fur2},
  answering a long standing question of Feldman and Moore.
\item The \emph{fundamental group} of $\cR$ is trivial. Note that this fundamental group is defined as
the group of $t > 0$ such that $t$ is the compression constant for some stable orbit equivalence between $G \actson X$ and itself. If $\pi : A \recht
B$ is a stable orbit equivalence with compression constant $t \geq 1$, Theorem \ref{thm.OEsuperrigidity} implies that $t=n/|\Ker \theta|$, where
$\theta : G \recht G$ has finite kernel, satisfies $n = [G:\theta(G)]$ and where $G \actson X$ is conjugate to $G \actson \Ind_G^G(X,\theta)$. Since
the action $G \actson X$ is weakly mixing, the induction is trivial, i.e.\ $n = 1$. This implies that $t \leq 1$ and hence, $t=1$.
\item The \emph{outer automorphism group} $\Out \cR = \Aut \cR / \Inn \cR$ of $\cR$ can be described as follows. Recall first that $\Aut \cR$ is defined as the group of orbit equivalences $\Delta : X \recht X$ of $G
\actson X$ with itself. The full group (see note on page
\pageref{note.fullgroup}) of $\cR$ is a normal subgroup of $\Aut \cR$
and denoted by $\Inn \cR$.
The subgroup $\Aut^*(X,G) \subset \Aut \cR$ consists of those $\Delta$ satisfying
$$\Delta(g \cdot x) = \delta(g) \cdot
\Delta(x) \;\; \text{almost everywhere,}$$ for some group automorphism $\delta \in \Aut G$. For our given $\cR$, $\Out \cR$ is the image of
$\Aut^*(X,G)$ through the quotient map $\Aut \cR \recht \Out \cR$. Weak mixing then implies that $\Out \cR \cong \Aut^*(X,G)/G$.
\end{itemize}
\end{remark}

\begin{remark}\label{rem.atomic-bernoulli}
Let $G$ be a group with infinite normal subgroup $H$ with the relative property (T). Let $G \actson (X,\mu)$ be a strongly malleable action whose
restriction to $H$ is weakly mixing. Then, the conclusions of Theorems \ref{thm.cocycle} and \ref{thm.OEsuperrigidity} hold as well for all
\emph{quotient actions} $G \actson (Y,\eta)$ provided that the quotient map $X \recht Y$ satisfies a \emph{relative weak mixing} property,
introduced by Popa in \cite{P0} (Definition 2.9). Indeed, if for a measurable map $w : X \recht K$ and a homomorphism $\theta : G \recht K$, the
$1$-cocycle $G \times X \recht K : (g,x) \mapsto w(g \cdot x) \theta(g) w(x)^{-1}$ actually is a map $G \times Y \recht K$, then relative weak
mixing imposes that $w$ is already a map $Y \recht K$.

Hence, the conclusions of Theorems \ref{thm.cocycle} and \ref{thm.OEsuperrigidity} hold for all generalized Bernoulli actions that are free and
weakly mixing restricted to $H$, even starting from an atomic base space.
\end{remark}

In fact, Theorem \ref{thm.OEsuperrigidity} follows from the cocycle superrigidity theorem \ref{thm.cocycle} and the following classical lemma.

\begin{lemma}
Let $G \actson (X,\mu)$ and $\Gamma \actson (Y,\eta)$ be free ergodic actions that are stably orbit equivalent. If the associated
$1$-cocycle is cohomologous to a homomorphism $\theta : G \recht \Gamma$, then the conclusion of Theorem \ref{thm.OEsuperrigidity} holds.
\end{lemma}
\begin{proof}
The proof of the lemma consists of two easy translation statements. In the first paragraph, stable orbit equivalence is translated as \emph{measure
equivalence} (see e.g.\ \cite{Fur2}): we get a natural space with an infinite measure preserving action of $G \times \Gamma$. In a second paragraph,
the conclusion follows using the triviality of the cocycle.

Let $p : X \recht Y$ be the equivalence relation preserving map as in
the construction of the $1$-cocycle $\al$ above. Take symmetrically $q
: Y \recht X$ and the $1$-cocycle $\be : \Gamma \times Y \recht G$.
We denote by $g \cdot x$ the action of $G$ on $X$ and by $s \ast y$ the action of $\Gamma$ on
$Y$. Define commuting actions of $G$ and $\Gamma$ on $X \times \Gamma$
and $Y \times G$ respectively, by the formulae
$$g \cdot (x,s) \cdot t = (g \cdot x,\al(g,x) s t) \quad,\qquad
s \ast (y,g) \ast h = (s \ast y,\be(s,y)gh) \; .$$
Following Theorem 3.3 in \cite{Fur2}, we prove that there is a natural $G \times \Gamma$-equivariant measure space isomorphism $\Theta : X \times \Gamma \recht Y \times G$
satisfying $\Theta(x,s) \in (\Gamma \ast p(x)) \times G$ for almost
all $(x,s)$. Indeed, take measurable maps $X \recht G : x \mapsto g_x$ and $Y \recht \Gamma: y \mapsto s_y$
such that $q(p(x)) = g_x \cdot x$ and $p(q(y)) = s_y \ast y$ almost everywhere. Define
\begin{align*}
\Theta &: X \times \Gamma \recht Y \times G : \Theta(x,s) = (s^{-1} \ast p(x),\be(s^{-1},p(x))g_x) \\
\Theta^{-1} &: Y \times G \recht X \times \Gamma : \Theta^{-1}(y,g) = (g^{-1} \ast q(y),\al(g^{-1},q(y))s_y) \; .
\end{align*}

The assumption of the lemma yields a homomorphism $\theta : G \recht \Gamma$
and a measurable map $w : X \recht \Gamma$ such that $\al(g,x) = w(g \cdot x) \theta(g) w(x)^{-1}$ almost everywhere. So, the map
$\Psi : X \times \Gamma \recht X \times \Gamma : \Psi(x,s) = (x,w(x)s)$
is a measure space isomorphism that is equivariant in the following
sense
$$\Psi(g \cdot x,\theta(g) s t) = g \cdot \Psi(x,s) \cdot t \; .$$
So, $\Theta \circ \Psi$ conjugates the new commuting actions $g(x,s)t = (g \cdot
x,\theta(g)st)$ on $X \times \Gamma$ with the commuting actions on $Y
\times G$ given above. In particular, the new action of $G$ on $X
\times \Gamma$ has a fundamental domain of finite measure. Having a fundamental
domain forces $\Ker \theta$ to be finite, while its being of finite
measure imposes $\theta(G)$ to be of finite index in $G$. Finally, the
new action of $\Gamma$ on the quotient $G \backslash (X \times
\Gamma)$ is exactly $\Gamma \actson \Ind_G^\Gamma(X,\theta)$ and $\Theta \circ \Psi$
induces a conjugacy of the actions $\Gamma
\actson \Ind_G^\Gamma(X,\theta)$ and $\Gamma \actson Y$.
\end{proof}

There is a slightly more general way of writing \lq obviously\rq\ stably orbit equivalent actions, by first restricting $G \actson X$ to $G_0 \actson
X_0$, where $G_0$ is a finite index subgroup of $G$ and $G \actson X$ is induced from $G_0 \actson X_0$. Since the superrigid actions in this talk
are all weakly mixing, they are not induced in this way.

It remains to prove the cocycle superrigidity theorem \ref{thm.cocycle}. This proof occupies the rest of the section and consists of several steps.
\begin{somoptel}\setcounter{teller}{-1}
\item Using the \emph{weak mixing} property and the fact that $\cU(B)$ is a
  \emph{Polish group with a bi-invariant metric}, restrict to the case $K = \cU(B)$.
\end{somoptel}

\noindent The $1$-cocycle $\gamma : G \times X \recht \cU(B)$ is then interpreted
as a family of unitaries $\gamma_g \in \cU(A \ot B)$, where $A =
L^\infty(X,\mu)$. Moreover, strong malleability yields $(\al_t)$ and
$\be$ on $A \ot A$.

\begin{somoptel}\label{page.program}
\item Using the \emph{relative property (T)}, find $t_0 >
  0$ and a non-zero partial isometry $a \in A \ot A \ot B$ satisfying
\begin{equation} \label{eq.joepie}
(\gamma_g)_{13} (\si_g \ot \si_g \ot \id)(a) = a (\al_{t_0} \ot
\id)((\gamma_g)_{13}) \tag{$\ast$}
\end{equation}
for all $g \in H$. We use the notation $(a \ot b)_{13} := a \ot 1 \ot
b$ and extend to $u_{13}$ for all $u \in A \ot B$ by linearity and continuity.
\item Using the \emph{period $2$ automorphism given by the
    strong malleability} and the \emph{weak mixing property of the action
    restricted to $H$}, glue together partial isometries, in order
    to get \eqref{eq.joepie} with $t_0 = 1$, i.e.\ a non-zero partial isometry $a \in A \ot A \ot B$ satisfying
$$(\gamma_g)_{13} (\si_g \ot \si_g \ot \id)(a) = a (\gamma_g)_{23}$$
for all $g \in H$.
\item Deduce from the previous equality, using the \emph{intertwining-by-bimodules
  technique}, a non-zero partial isometry $v \in A \ot B$ and
  partial isometries $\theta(g) \in B$ such that
$$\gamma_g (\si_g \ot \id)(v) = v(1 \ot \theta(g))$$
for all $g \in H$.
\item Using a \emph{maximality argument}, glue together such partial
  isometries $v$ in order to get a unitary $v$ satisfying the same
  formula.
\item Use the \emph{normality} of $H$ in $G$ and the weak
  mixing property of the action restricted to $H$, to extend the formula to
  $g \in G$.
\end{somoptel}

Lemma \ref{lem.reduce-subgroup} covers step (0), Lemma
\ref{lem.crucialstep} covers steps (1), (2) and (3), Lemma
\ref{lem.maximality} covers step (4) and the final step (5) is done in the
proof of the theorem.

To prove step (0) of the program, the
essential property of the Polish group $\cU(B)$ that we retain is the
existence of a \emph{bi-invariant metric} $d(u,v) = \|u-v\|_2$.

\begin{lemma} \label{lem.reduce-subgroup}
Let $G$ act weakly mixingly on $(X,\mu)$. Let $\cG$
be a Polish group with a bi-invariant complete metric $d$ and let $K
\subset \cG$ be a closed subgroup.
Suppose that $\gamma : G \times X \recht K$ is a $1$-cocycle. Let $v : X \recht \cG$ be a measurable map and
$\theta : G \recht \cG$ a homomorphism such that
$$\gamma(g,x) = v(g \cdot x) \theta(g) v(x)^{-1}$$
almost everywhere. Whenever $v_0 \in \cG$ is an essential value of the function $v$, we have $v(x)v_0^{-1} \in K$ almost everywhere and
$v_0\theta(g)v_0^{-1} \in K$ for all $g \in G$.
\end{lemma}
\begin{proof}
Let $v_0$ be an essential value of the function $v$. Changing $v(x)$
into $v(x)v_0^{-1}$ and $\theta$ into $(\Ad v_0) \circ \theta$, we assume that $e$ is an
essential value of $v$ and prove that $\theta(g) \in K$ for all $g \in G$ and $v(x) \in K$ almost everywhere.

Denote by $d$ the bi-invariant metric on the $\cG$. Choose $\eps > 0$ and $g \in G$. Take $W \subset X$ with
$\mu(W) > 0$ such that $d(v(x),1) < \eps/4$ for all $x \in W$. Take $k \in G$ such that $\mu(k \cdot W \cap W) > 0$ and $\mu((gk)^{-1} \cdot
W \cap W) > 0$. If $x \in k \cdot W \cap W$, we have $d(v(x),1),d(v(k^{-1} \cdot x),1) < \eps /4$. It follows that $d(\theta(k^{-1}),K) <
\eps/2$. In the same way, $d(\theta(gk),K) < \eps/2$. Together, $d(\theta(g),K) < \eps$. This holds for all $\eps > 0$ and all $g \in G$ and hence,
$\theta(G)\subset K$.

Let $\eps > 0$. The formula $v(g \cdot x) = \gamma(g,x) v(x) \theta(g)^*$ almost everywhere, yields that $\{ x \in X \mid d(v(x),K) < \eps \}$ is
non-negligible and $G$-invariant, hence, the whole of $X$. It follows that $v(x) \in K$ almost everywhere.
\end{proof}

We fix the
following data and notations.

\begin{somoptel}
\item Let $G$ be a countable group with infinite normal subgroup $H$
  such that $(G,H)$ has the relative property (T). Let $G$ act strongly malleably
on $(X,\mu)$ and suppose that its restriction to $H$ is weakly mixing. Write $A = L^\infty(X)$ and write the action of $G$ on $A$ as
$(\si_g (F))(x) = F(g^{-1} \cdot x)$.
\item Let $\gamma : G \times X \recht \cU(B)$ be a $1$-cocycle with values in
  the unitary group of the \emph{II$_1$ factor} $(B,\tau)$. Remark
  that we can indeed suppose that $B$ is a II$_1$ factor\footnote{Any finite $(B,\tau)$ can
  be embedded, in a trace-preserving way, into a II$_1$ factor, e.g.\ into $\bigl(\bigotimes_{n \in \Z}
  (B,\tau)\bigr) \rtimes \Z$ and $\cU(B)$ is then a closed subgroup of
  the unitary group of this II$_1$ factor.}.
We write $\gamma_g \in \cU(A \ot B)$, given by $\gamma_g(x) =
  \gamma(g,g^{-1} \cdot x)$. The $1$-cocycle relation becomes
$$\gamma_g (\si_g \ot \id)(\gamma_h) = \gamma_{gh} \quad\text{for all}\;\; g,h \in G
\; .$$
\item We denote by $(\rho_g)$ the following action of $G$ by
  automorphisms of $A \ot B$:
$$\rho_g(a) = \gamma_g (\si_g \ot \id)(a) \gamma_g^* \quad\text{for all}\;\; a
\in A \ot B \; .$$
\item We denote by $(\eta_g)$ the unitary representation of $G$ on
  $L^2(X) \ot L^2(B)$ given by
$$\eta_g(a) = \gamma_g (\si_g \ot \id)(a) \quad\text{for all}\;\; a \in A
\ot B \subset L^2(X) \ot L^2(B) \; .$$
\item We denote, for every $t \in \R$, by $(\pi^t_g)$ the unitary
  representation on  $L^2(X \times X) \ot L^2(B)$ of $G$ given by
$$\pi^t_g(a) = (\gamma_g)_{13} (\si_g \ot \si_g \ot \id)(a) (\al_t \ot
\id)((\gamma_g)_{13}^*)$$
for all $a \in A \ot A \ot B \subset L^2(X \times X) \ot
L^2(B)$. Recall the notation $u_{13}$ determined by $(a \ot b)_{13} =
a \ot 1 \ot b$.
\end{somoptel}

We cover steps (1), (2) and (3) of the program in the next lemma.

\begin{lemma} \label{lem.crucialstep}
Let $q \in A \ot B$ be a non-zero projection which is $\rho|_H$-invariant. There exists a non-zero partial isometry $v \in A \ot B$, a projection $p
\in B$ and a homomorphism $\theta : H \recht \cU(pBp)$ such that $vv^* \leq q$, $v^*v = 1 \ot p$ and
$$\gamma_h (\si_h \ot \id)(v) = v (1 \ot \theta(h))$$
for all $h \in H$.
\end{lemma}
\begin{proof}
{\it Step (1)} Note that $1$ is a $\pi^0_G$-invariant vector. The relative property $(T)$ yields a $t_0 = 2^{-n}$ and a non-zero element $a \in A \ot A \ot B$ such
that $a$ is $\pi^{t_0}_H$-invariant and such that $\|a - 1\|_2 \leq \|q\|_2 / 2$. It follows that $a (\al_{t_0} \ot \id)(q_{13}) \neq 0$, which
remains $\pi^{t_0}_H$-invariant. Taking the polar decomposition of $a(\al_{t_0} \ot \id)(q_{13})$, we get a non-zero partial isometry $a \in A \ot A \ot B$ which is
$\pi^{t_0}_H$-invariant and satisfies $a^*a \leq (\al_{t_0} \ot \id)(q_{13})$. Moreover, Proposition \ref{prop.mixing} yields
$$aa^* \;\; , \;
(\al_{-t_0} \ot \id) (a^*a) \in (A \ot B)^{\rho|_H}_{13}  \; .$$ So, we have a projection $\qtil \in (A \ot B)^{\rho|_H}$ such that $\qtil \leq q$ and
$$a^* a = (\al_{t_0} \ot \id)(\qtil_{13}) \; .$$

{\it Step (2)} Whenever $a$ and $b$ are $\pi^{t_0}_H$-invariant, we have that $a(\al_{t_0} \ot \id)(b)$ is $\pi^{2t_0}_H$-invariant and that $(\be \ot \id)(a)$ and
$(\al_{-t_0} \ot \id)(a^*)$ are $\pi^{-t_0}_H$-invariant. So, if we define
$$a_1 = (\al_{t_0} \ot \id)\bigl( (\be \ot \id)(a^*) a \bigr)$$
we get that $a_1$ is $\pi^{2t_0}_H$-invariant and
satisfies
$$a_1 a_1^* = \qtil_{13} \quad\text{and}\quad a_1^* a_1 = (\al_{2t_0}
\ot \id)(\qtil_{13}) \; .$$
Iterating the procedure yields at stage $n$ a partial isometry $b \in
A \ot A \ot B$ which is $\pi^1_H$-invariant and satisfies $bb^* =
\qtil_{13}$ and $b^*b = \qtil_{23}$.

{\it Step (3)} Define the (non-zero) operator $T \in \B(L^2(X)) \ot B$ by
$$(T \xi)(x) = \int_X b(x,y)\xi(y) \; d\mu(y) \quad\text{for all}\;\;
\xi \in L^2(X) \ot B \; .$$
We get
$$[T,\eta_h] = 0 \;\;\text{for}\; h \in H \; , \quad \qtil T = T =
T\qtil \; , \quad \| (\Tr \ot \id)(T^*T)\| < \infty \; .$$
Taking a spectral projection $P$ of $T$, we get an non-zero orthogonal
projection $P$ with the same properties as $T$. It follows that the
range of $P$ is a finitely generated right $B$-submodule of $(L^2(X)
\ot L^2(B))_B$ which is stable under $(\eta_h)_{h \in H}$.

As in Proposition \ref{prop.translate}, we get $n \geq 1$, a non-zero projection $p \in \M_n(\C) \ot B$, a non-zero partial isometry $v \in A \ot
\M_{1,n}(\C) \ot B$ and a homomorphism $\theta : H \recht \cU(p (\M_n(\C) \ot B)p)$ such that
$$\gamma_h (\si_h \ot \id)(v) = v(1 \ot \theta(h)) \;\;\text{for}\; h \in H
\; , \quad \qtil v = v \; , \quad v(1 \ot p) = v \; .$$ Since $v^* v$ is $(\sigma_h \ot \Ad \theta(h))$-invariant for all $h \in H$, it follows from
Proposition \ref{prop.mixing} that $v^*v = 1 \ot p_0$ for some non-zero projection $p_0 \in p (\M_n(\C) \ot B)p \cap \theta(H)'$. Since $p_0$
commutes with $\theta(H)$, we can cut down by $p_0$. Since moreover $\tau(p_0) \leq 1$, we can move $p_0$ into the upper corner of $\M_n(\C) \ot B$
and we have found a non-zero partial isometry $v \in A \ot B$, a non-zero projection $p \in B$ and a homomorphism $\theta : H \recht \cU(pBp)$ such
that $vv^* \leq q$, $v^*v = 1 \ot p$ and
$$\gamma_h (\si_h \ot \id)(v) = v(1 \ot \theta(h))$$
for all $h \in H$.
\end{proof}

We cover step (4) of the program in the following lemma.

\begin{lemma} \label{lem.maximality}
There exists a unitary element $v \in A \ot B$ and a homomorphism
$\theta : H \recht \cU(B)$ such that
$$\gamma_h (\si_h \ot \id)(v) = v (1 \ot \theta(h))$$
for all $h \in H$.
\end{lemma}
\begin{proof}
The proof is a straightforward maximality argument. Consider the set $\cI$
of partial isometries $v \in A \ot B$ for which there exist $p \in B$
and $\theta : H \recht \cU(pBp)$ satisfying
$$v^*v= 1 \ot p \quad\text{and}\quad \gamma_h (\si_h \ot \id)(v) = v (1 \ot \theta(h))$$
for all $h \in H$. Partially order $\cI$ by extension of partial
isometries and let $v$ be a maximal element of $\cI$. Write
$v^*v = 1
\ot p$. If $vv^* \neq 1$, put $q
= 1 - vv^*$. Then, $q \in (A \ot B)^{\rho|_H}$ and Lemma
\ref{lem.crucialstep} yields a non-zero partial isometry $w \in A \ot B$, a projection $e \in B$ and a homomorphism $\theta : H \recht \cU(eBe)$ such
that $ww^* \leq q$, $w^*w = 1 \ot e$ and
$$\gamma_h (\si_h \ot \id)(w) = w(1 \ot \theta(h))$$
for all $h \in H$.
Since $e \inside 1 - p$ in the II$_1$ factor $B$, we contradict the maximality $v$.
\end{proof}

\begin{proof}[Proof of Theorem \ref{thm.cocycle}]
Using Lemma \ref{lem.reduce-subgroup}, it is sufficient to prove the existence of a unitary $v \in A \ot B$ and a homomorphism $\theta : G \recht
\cU(B)$ such that
\begin{equation}\label{eq.solution}
\gamma_g (\si_g \ot \id)(v) = v (1 \ot \theta(g))
\end{equation}
for all $g \in G$. Take $v$ and $\theta$ as given by Lemma \ref{lem.maximality}. Fix $g \in G$ and write
$$\vtil = \gamma_g (\si_g \ot \id)(v) \quad\text{and}\quad \tetil(h) =
\theta(g^{-1}hg) \;\;\text{for}\; h \in H \; .$$
Obviously, $\gamma_h (\si_h \ot \id)(\vtil) = \vtil (1 \ot \tetil(h))$ for
all $h \in H$. It follows that
$$(\si_h \ot \id)(\vtil^* v) = (1 \ot \tetil(h)^*)\vtil^* v (1 \ot
\theta(h))$$ for all $h \in H$. Since $\vtil^* v$ is a unitary, the same proof as the one for Proposition \ref{prop.mixing}, yields a unitary $u \in
B$ such that $\tetil = (\Ad u) \theta$ and $\vtil = v(1 \ot u^*)$. So, for any $g \in G$, we find a unique unitary element $\theta(g) \in \cU(B)$
such that \eqref{eq.solution} holds. By uniqueness, $\theta$ is a homomorphism and we are done.
\end{proof}

\section{Non-orbit equivalent actions and 1-cohomology}

The following theorem is an immediate consequence of Theorem \ref{thm.cocycle}.

\begin{theorem}[Popa, Sasyk, \cite{PS}]
Let $G$ be a countable group with infinite normal subgroup $H$ such that $(G,H)$ has the relative property (T). Let $(\si_g)$ be the Bernoulli action (with non-atomic base) of $G$ on $(X,\mu)$. Then, $H^1(\si)= \Char G$.
\end{theorem}

Through the following lemma, one can easily produce non stable orbit equivalent actions

\begin{lemma}
Let $G$ be a countable group and $K$ a compact abelian group. Let $G
\times K$ act on $(X,\mu)$ and denote by $(\si_g \rho_k)$ the
corresponding action on $A =
L^\infty(X)$. Define
$B = A^K$, the algebra of $K$-fixed points. Denote by $(\si^K_g)$ the restriction of $(\si_g)$ to $B$. Assume that
\begin{itemize}
\item $(\si_g)$ is free and weakly mixing,
\item $(\si^K_g)$ is still free,
\item $H^1(\si) = \Char G$.
\end{itemize}
Then, $H^1(\si^K) = \Char G \times \Sp(K,\rho)$, where
$$\Sp(K,\rho) = \{\al \in \Char(K) \mid \; \exists u \in \cU(A), \rho_k(u) = \al(k) u \;\;\text{for all}\; k \in K \} \; .$$
\end{lemma}
\begin{proof}
Whenever $u \in \cU(A)$ and $\rho_k(u) = \alpha(k) u$ for all $k \in K$, we define $\om_g \in B$ by the formula $\om_g = u \si_g(u^*)$. Using the
weak mixing of $(\si_g)$, it is easy to check that we obtain an embedding $\Char G \times \Sp(K,\rho) \hookrightarrow H^1(\si^K)$. Suppose on the
contrary that the $1$-cocycle $\om$ defines an element of
$H^1(\si^K)$. We regard $\om$ as a $1$-cocycle for $\si$ and since $H^1(\si) = \Char G$, we find that $\om$ is cohomologous to a character of $G$. Subtracting this character from $\om$, we may assume that $\om_g = u \si_g(u^*)$ for some unitary $u
\in \cU(A)$. Since for any $k \in K$, $\om_g$ is $K$-invariant and since $(\si_g)$ is weakly mixing, we conclude that there exists $\alpha : K \recht
S^1$ such that $\rho_k(u) = \alpha(k) u$ for all $k \in K$. But this means that $\om$ is given by an element of $\Sp(K,\rho)$.
\end{proof}

The following proposition immediately follows.

\begin{proposition}[Popa, \cite{P3}] \label{prop.non-orbit-equiv}
Let $G$ be a countable group with infinite normal subgroup $H$ such that $(G,H)$ has the relative property (T). Let $\Gamma$ be any countably
infinite abelian group and $K=\widehat{\Gamma}$. Denote by $(\si_g)$ the Bernoulli action of $G$ on $L^\infty(X,\mu) = \otimes_{g \in
  G} L^\infty(K,\text{\rm Haar})$ and define $(\rho_k)_{k \in K}$ as the diagonal action on $L^\infty(X,\mu)$ of the translation action of $K$ on $L^\infty(K)$.
Define
$(\si^K_g)$ as the restriction of $(\si_g)$ to the $K$-fixed points $L^\infty(X)^K$.

Then, $(\si^K_g)$ is a free and ergodic action of $G$ satisfying $H^1(\si^K) = \Char G \times \Gamma$.
\end{proposition}

\begin{remark}
It follows that any countable group $G$ that admits an infinite normal subgroup $H$ such that $(G,H)$ has the relative property (T), admits a
continuous family of non stably orbit equivalent actions. Indeed, $\Char G$ being compact, an isomorphism $\Char G \times \Gamma_1 \cong \Char G
\times \Gamma_2$ entails a virtual isomorphism between $\Gamma_1$ and $\Gamma_2$. It is not hard to exhibit a continuous family of non virtually
isomorphic countable abelian groups.
\end{remark}

\section{\mbox{Intertwining rigid subalgebras of crossed products}} \label{sec.intertwine}

The major aim of the rest of the talk is to prove Popa's \emph{von Neumann
  strong rigidity theorem} for Bernoulli actions of $w$-rigid groups, deducing conjugacy of actions out of their
  mere von Neumann equivalence. This is more difficult, but
  nevertheless related to the orbit equivalence superrigidity Theorem
  \ref{thm.OEsuperrigidity}. In particular, the crucial lemma
  \ref{lem.intertwine-rigid} below, is the von Neumann counterpart to
  Lemma \ref{lem.crucialstep}, covering steps (1), (2) and (3) of the
  program on page \pageref{page.program}. It states that in a crossed
  product $M:=N \rtimes G$ by a malleable mixing action, a subalgebra
  $Q \subset M$ with the relative property (T), can be essentially
  conjugated into $\cL(G)$.

But, the aim of this section is not only preparation to the von
Neumann strong rigidity theorem. The results are applied as well in
the next section in order to construct II$_1$~factors with prescribed countable
fundamental groups. For this reason, we need to deal with actions on
non-tracial (but almost-periodic) algebras.

We refer to page \pageref{page.explanation} for a rough explanation of the idea of
the proof of Lemma \ref{lem.intertwine-rigid}. It is another
application of Popa's \emph{deformation/rigidity strategy}. The
deformation property of \emph{malleability} is played against the relative property
(T). For this, we need the notion of relative
property (T) for an inclusion $Q \subset M$ of finite von Neumann
algebras (see Definition \ref{def.relativeT}).
The mixing property of the action has several von Neumann algebraic consequences that are used
throughout and proved in Appendix~\ref{sec.mixing}. Finally, in order to
actually conjugate (essentially) $Q$ into $\cL(G)$, Popa's
\emph{intertwining-by-bimodules} technique is used (see Appendix~\ref{sec.intertwining}).

\begin{lemma} \label{lem.intertwine-rigid}
Given a \emph{strongly malleable} mixing action of a countable group $G$ on an almost
periodic $(\cN,\vphi)$, write $N=\cN^\vphi$. Let $Q \subset N \rtimes
G$ be a \emph{diffuse} subalgebra with the \emph{relative property
  (T)}. Denote by $P$ the quasi-normalizer of $Q$ in $N \rtimes G$ and
suppose that there is no non-zero homomorphism from $P$ to an
amplification of $N$.

Then, there exists $\gamma > 0$, $n \geq 1$ and a non-zero partial
isometry $v \in \M_{n,1}(\C) \ot (\cN \rtimes G)$ which is a
$\gamma$-eigenvector for $\vphi$ and satisfies
$$v^*v \in P \cap Q' \; , \quad vPv^* \subset \M_n(\C) \ot \cL(G) \; .$$
\end{lemma}

\begin{proof}
In the course of this proof, we use the following terminology: given subalgebras $Q_1,Q_2$ of a von Neumann algebra, an element $a$ is said to be
$Q_1$-$Q_2$-finite, if there exists finite families $(a_i)$ and $(b_i)$ such that
$$a Q_2 \subset \sum_{i=1}^n Q_1 a_i \quad\text{and}\quad Q_1 a \subset \sum_{i=1}^m b_i Q_2 \; .$$
Hence, the $Q$-$Q$-finite elements are nothing else but the elements quasi-normalizing $Q$.

{\it Step (1), using relative property (T).}
Take $(\alpha_t)$ and $\be$ as in Definition
\ref{def.malleable}. Write $\Ntil = (\cN \ot \cN)^{\vphi \ot \vphi}$
and $\Mtil = \Ntil \rtimes G$. Write $M = N \rtimes G$ and consider
$M$ as a subalgebra of $\Mtil$ by considering $\cN \ot 1 \subset \cN
\ot \cN$. Extend $(\al_t)$ and $\be$ to $\Mtil$. The
relative property (T) yields $t_0=2^{-n}$ and a non-zero element $w
\in \Mtil$ such that $xw = w \alpha_{t_0}(x)$ for all $x \in Q$.

{\it Step (2), finding a non-zero element $a \in \Mtil$ that is
  $Q$-$\alpha_1(Q)$-finite, using the period $2$-automorphism $\be$.}
Denote by $\cP$ the $^*$-algebra of $Q$-$Q$-finite elements in $M$. By definition, $P$
  is the weak closure of $\cP$. Whenever $y \in \cP$, the
  element $\alpha_{t_0}(\beta(w^*)yw)$ is
  $Q$-$\alpha_{2t_0}(Q)$-finite. It suffices to find $y$ such that
  $\beta(w^*)yw$ is non-zero, since we can then continue to find a
  non-zero $Q$-$\alpha_1(Q)$-finite element $a$ in $\Mtil$. Denote by
  $p$ the supremum of all range projections of elements $yw$, where
  $y \in \cP$. We have to prove that $p \beta(w)
  \neq 0$. By construction, $p \in \Mtil \cap P'$ and $pw=w$.
From
Proposition \ref{prop.mixing-two} (and here we use that there is no non-zero homomorphism from $P$ to an amplification of $N$), $\Mtil \cap P'
\subset M$ and so, $p \in M$. But, $\beta$ acts trivially on $M$ and
  we obtain $p \beta(w) = \beta(pw) = \beta(w) \neq 0$.

{\it Step (3), using the intertwining-by-bimodules technique to conclude.} Denote by $f \in \la \Mtil,e_{\alpha_1(M)} \ra \cap Q'$ the orthogonal
projection onto the closure of $Q \, a \, \alpha_1(M)$ in $L^2(\Mtil)$ and remark that $\vphih(f) < +\infty$. Denoting by $\cF : \la (\cN \ot
\cN)\rtimes G,e_{(1 \ot \cN) \rtimes G} \ra \recht \la \cN \rtimes G,e_{\cL(G)} \ra$ the $\vphih$-preserving conditional expectation, it follows that
$$\cF(f) \in \la \cN \rtimes G,e_{\cL(G)} \ra \cap Q'
\quad\text{with}\quad \vphih(\cF(f)) < \infty \; .$$ Moreover, $\cF(f) \neq 0$ since $\cF$ is faithful.

From Proposition \ref{prop.translate}, we get $\gamma > 0$, $n \geq 1$, $p \in \M_n(\C) \ot \cL(G)$, a homomorphism $\theta : Q \recht p(M_n(\C) \ot
\cL(G))p$ and a non-zero partial isometry $w \in M_{1,n}(\C) \ot (\cN \rtimes G)$ such that $w$ is a $\gamma$-eigenvector for $\vphi$ and $xw = w
\theta(x)$ for all $x \in Q$. It follows that $w^*w \in p(M_n(\C) \ot (N \rtimes G))p \cap \theta(Q)'$, which is included in $p(M_n(\C) \ot \cL(G))p$
by Theorem \ref{thm.mixing-one}. Also $ww^* \in M \cap Q'$ and hence, $w^*Qw$ is a diffuse subalgebra of $p(M_n(\C) \ot \cL(G))p$. Applying once more
Theorem \ref{thm.mixing-one}, we get $w^*Pw \subset p(M_n(\C) \ot \cL(G))p$. Since obviously $M \cap Q' \subset P$, we can take $v = w^*$ to
conclude.
\end{proof}
\begin{remark} \label{rem.malleable}
If $P$ is a factor, it is sufficient to assume \emph{malleability}
instead of \emph{strong malleability}. Indeed, looking back at the
proof, let $a \in \Mtil$ be a $Q$-$\al_{t_0}(Q)$-finite element. Then,
$a \al_{t_0}(y a)$ is $Q$-$\al_{2 t_0}(Q)$-finite for every $y \in
\Mtil$ that quasi-normalizes $Q$. Denote by $\Ptil$ the
quasi-normalizer of $Q$ in $\Mtil$. It is then sufficient to show that
$\Ptil$ is factorial, to obtain at least one $y$ such that $a
\al_{t_0}(y a) \neq 0$. As in the proof above, $\Mtil \cap P' \subset
M$. Since $\Ptil$ contains $P$, it follows that $\Mtil \cap \Ptil'
\subset M \cap P' = \cZ(P) = \C 1$. So, we are done.
\end{remark}

In two cases, a unitary intertwiner $v$ can be found. The
first case is easy and follows immediately:
assume $G$ to be ICC and the quasi-normalizer $P$ to be a factor.
It is crucial to allow as well for an amplification in order to apply
the result when dealing with
the fundamental group of the crossed product $N \rtimes G$.

\begin{theorem}[Popa, \cite{P1}] \label{thm.relative-rigid}
Given a malleable mixing action of an ICC group $G$ on an almost periodic $(\cN,\vphi)$, write $N=\cN^\vphi$ and $M = N \rtimes G$. Let $t > 0$ and
let $Q \subset M^t$ be a \emph{diffuse} subalgebra with the \emph{relative property
  (T)}. Denote by $P$ the quasi-normalizer of $Q$ in $M^t$.
Suppose that $P$ is a \emph{factor} and that there is no non-zero homomorphism from $P$ to an
amplification of $N$. Realize $M^t = p (\M_n(\C) \ot M)p$.

Then, there exists $\gamma > 0$, $k \geq 1$ and $v \in \M_{n,k}(\C) \ot
(\cN \rtimes G)$ a $\gamma$-eigenvector for $\vphi$, such that
$$v^*v = p \; , \quad q:=vv^* \in \M_k(\C) \ot \cL(G) \; , \quad vPv^*
\subset \cL(G)^{t\gamma} \; ,$$
where we have realized $\cL(G)^{t\gamma}:=q(\M_k(\C) \ot \cL(G))q$.
\end{theorem}
\begin{proof}
Choose a projection $q \in \M_k(\C) \ot Q$ with trace $s$ where $s=1/t$. Write $Q^s := q(\M_k(\C) \ot Q)q$ and $P^s := q(\M_k(\C) \ot P)q$. We
consider $Q^s \subset P^s \subset M$. Clearly, $Q^s$ is diffuse, $Q^s \subset M$ has the relative property (T) by Proposition \ref{prop.Tcutdown} and
$P^s$ is the quasi-normalizer of $Q^s$ by Lemma \ref{lem.quasi-normalizer}. So, Lemma \ref{lem.intertwine-rigid} (with Remark \ref{rem.malleable})
yields a partial isometry $v$ which is a $\gamma$-eigenvector for $\vphi$ and satisfies $v^*v \in P^s$, $v P^s v^* \subset \cL(G)^{\gamma}$. Since
both $P^s$ and $\cL(G)$ are factors, we can move around $v$ using partial isometries in matrix algebras over $P$ and $\cL(G)$ to conclude.
\end{proof}

In the tracial case, assuming $G$ to be ICC is sufficient.

\begin{theorem}[Popa, \cite{P1}] \label{thm.rigidity}
Given a strongly malleable mixing action of an ICC group $G$ on a finite $(N,\tau)$, let $t > 0$ and let $Q \subset (N \rtimes G)^t$ be a
\emph{diffuse} subalgebra with the \emph{relative property
  (T)}. Denote by $P$ the quasi-normalizer of $Q$ in $(N \rtimes G)^t$ and
suppose that there is no non-zero homomorphism from $P$ to an
amplification of $N$.

Then, there exists a unitary element $v \in (N \rtimes G)^t$ such that
$vPv^* \subset \cL(G)^t$.
\end{theorem}
\begin{proof}
Write $M=N \rtimes G$. Below we prove the
existence of a partial isometry $v \in M^t$ satisfying $v^*v \in P
\cap Q'$ and $vPv^* \subset \cL(G)^t$. Since any projection $p \in P \cap Q'$ of
trace $s$ yields an inclusion $pQ \subset pPp \subset M^{st}$
satisfying the assumptions of the theorem, a maximality argument
combined with the factoriality of $\cL(G)$ then allows to conclude.

Choose a projection $q \in \M_k(\C) \ot Q$ with trace $s$ where
$s=1/t$. Write $Q^s := q(\M_k(\C) \ot Q)q$ and $P^s := q(\M_k(\C) \ot
P)q$ as in the proof of the previous theorem. From Lemma
\ref{lem.intertwine-rigid}, we get a partial isometry $w \in M$
satisfying $w^*w \in P^s \cap (Q^s)'$ and $w P^s w^* \subset \cL(G)$.
Let $e$ be the smallest projection in $P \cap Q'$ satisfying $w^*w
\leq 1 \ot e$. Moving around $w$ using partial isometries in matrix
algebras over $Q$ and
$\cL(G)$, we find a partial isometry $v \in M^t$ satisfying $v^*v = e$
and $v P v^* \subset \cL(G)^t$.
\end{proof}

\begin{lemma} \label{lem.quasi-normalizer}
Let $Q \subset M$ be an inclusion of finite von Neumann algebras and
$p$ a non-zero projection in $Q$. If $P$ denotes the quasi-normalizer
of $q$ in $M$, then $pPp$ is the quasi-normalizer of $pQp$ in $pMp$.
\end{lemma}
\begin{proof}
Denote by $\Ptil$ the quasi-normalizer of $pQp$ in $pMp$. We only
prove the inclusion $pPp \subset \Ptil$, the converse inclusion being
analogous. Let $z$ be a central projection in $Q$ such that $z =
\sum_{i=1}^n v_i v_i^*$ with $v_i$ partial isometries in $Q$ and
$v_i^* v_i \leq p$.

If now $x \in M$ quasi-normalizes $Q$, we write $p_0= pz$ and claim that $p_0 x p_0$ quasi-normalizes $pQp$. Indeed, if $xQ \subset \sum_{k=1}^r Q
x_k$, it is readily checked that
$$p_0 x p_0 \; pQp \subset \sum_{k,i} pQp \; v_i^* x_k p \; .$$
Since the central support of $p$ in $Q$ can be
approximated arbitrary well by such special central projections $z$,
$p_0$ approximates arbitrary well $p$ and we have proved that $pPp \subset \Ptil$.
\end{proof}

\section{Fundamental groups of type II$_1$ factors} \label{sec.fundamental}

Recall that we denote the fundamental group of a II$_1$ factor $M$ by
$\fun(M) \subset \R^*_+$ and that $\Sp(\cN,\vphi) \subset \R^*_+$ denotes the point spectrum of the
modular automorphism group of an almost periodic state $\vphi$ on $\cN$.

\begin{theorem}[Popa, \cite{P1}] \label{thm.computation-fundamental}
Let $G$ be an ICC group that admits an infinite almost normal subgroup
$H$ with the relative property (T). Let $(\si_g)$ be a
malleable mixing action of $G$ on the almost periodic injective
$(\cN,\vphi)$. Denote $M:=\cN^\vphi \rtimes G$. One has
$$\Sp(\cN,\vphi) \subset \fun(M) \subset \Sp(\cN,\vphi) \fun(\cL(G))
\; .$$
In particular, if $\cL(G)$ has trivial fundamental group, $\fun(M) = \Sp(\cN,\vphi)$.
\end{theorem}
\begin{proof}
As shown by Golodets and Nessonov \cite{GN}, the inclusion
$\Sp(\cN,\vphi) \subset \fun(M)$ holds. Indeed,
let $s \in \Sp(\cN,\vphi)$ and take an $s$-eigenvector $v \in
\cN$, that we may suppose to be a partial isometry. Write $p=v^*v$ and
$q = vv^*$. Then, $p,q \in \cN^\vphi \subset M$, $\vphi(q) = s \vphi(p)$ and $\Ad
v$ yields an isomorphism of $pMp$ with $qMq$. Hence, $s \in
\fun(M)$.

Suppose $t \in \fun(M)$ and let $\theta : M \recht M^t$ be a $^*$-isomorphism. Since $H$ is almost normal in $G$, $\cL(G)$ is contained in the
quasi-normalizer of $\cL(H)$ in $M$. Moreover, $\cL(H)$ is diffuse since $H$ is infinite. So, it follows from Theorem \ref{thm.mixing-one} that the
quasi-normalizer of $\cL(H)$ in $M$ is exactly $\cL(G)$ and, in
particular, a factor. Since $\cN^\vphi$ is an injective von Neumann
algebra with finite trace $\vphi$, it follows from Remark \ref{rem.amenable} that there is no non-zero homomorphism from $\cL(G)$ to an amplification of $\cN^\vphi$.

Write $\cM = \cN \rtimes G$, $Q = \theta(\cL(H))$ and $P = \theta(\cL(G))$. Realize $M^t := p(\M_n(\C) \ot M)p$, where $p$ is chosen in $\M_n(\C) \ot
\cL(H)$. By Proposition \ref{prop.from-group-to-vnalg}, the inclusion $Q \subset P$ has the relative property (T). Increasing $n$ if necessary, the
previous paragraph and Theorem \ref{thm.relative-rigid} yield $s \in \Sp(\cM,\vphi)$ and $v \in \M_n(\C) \ot \cM$ such that $v$ is an
$s$-eigenvector for $\vphi$, $v^*v = p$, $q:= vv^* \in \M_n(\C) \ot \cL(G)$ and $vPv^* \subset q (\M_n(\C) \ot \cL(G))q$. We claim that this
inclusion is an equality. Then, we have shown that $\cL(G)$ and $\cL(G)^{t s}$ are isomorphic, which yields $t s \in \fun(\cL(G))$ and
hence, $t \in \Sp(\cN,\vphi) \fun(\cL(G))$. So, this ends the proof.

Changing $q$ to an equivalent projection in $\M_n(\C) \ot
\cL(G)$, we may assume that $q \in \M_n(\C) \ot
\cL(H)$. Write $Q_1 \subset P_1 \subset M$ as
$$Q_1 := \theta^{-1}\bigl(v^* (\M_n(\C) \ot \cL(H)) v\bigr) \quad\text{and}\quad
P_1 := \theta^{-1}\bigl(v^* (\M_n(\C) \ot \cL(G)) v\bigr) \; .$$ The inclusion $Q_1 \subset M = N \rtimes G$ has the relative property (T), $P_1$ is
the quasi-normalizer of $Q_1$ and $\cL(G) \subset P_1$. We have to prove that $\cL(G) = P_1$.

By Theorem \ref{thm.relative-rigid}, there exists a
$w \in \M_{k,1}(\C) \ot \cM$, an $r$-eigenvector for $\vphi$
satisfying $w^*w = 1$ and $wP_1 w^* \subset \cL(G)^r$. Since
$\cL(G) \subset P_1$, Theorem \ref{thm.mixing-one} yields $w \in
\M_{k,1}(\C) \ot \cL(G)$. But then, $\cL(G) = P_1$ and we are done.
\end{proof}

\begin{corollary}
Let $G$ be an ICC group that admits an infinite almost normal subgroup
with the relative property (T). Suppose that $\cL(G)$ has trivial
fundamental group. Let $\Tr_\Delta$ be the faithful normal state on
$\B(H)$ given by $\Tr_\Delta(a) = \Tr(\Delta a)$
and define
$(\cN,\vphi) = \bigotimes_{g \in G} (\B(H),\Tr_\Delta)$, with
Connes-St{\o}rmer Bernoulli action $G \actson (\cN,\vphi)$. Write $M := \cN^\vphi \rtimes G$.

Then, $\fun(M)$ is the subgroup of $\R^*_+$ generated by the ratios
$\lambda/\mu$ for $\lambda,\mu$ belonging to the point
spectrum of $\Delta$. In particular, for every countable subgroup $S
\subset \R^*_+$, there exists a type II$_1$ factor with separable predual whose fundamental group is $S$.
\end{corollary}

Popa showed in \cite{P5} that, among other examples, $\cL(G)$ has
trivial fundamental group when $G = \SL(2,\Z) \ltimes \Z^2$. Note that
Popa shows in \cite{P5} that the fundamental group of
$\cL(G) = \SL(2,\Z) \ltimes L^\infty(\T^2)$ equals the fundamental
group of the equivalence relation given by the orbits of $\SL(2,\Z)
\actson \T^2$. The latter reduces to $1$ using Gaboriau's $\ell^2$ Betti number
invariants for equivalence relations, see \cite{Gab}.

It is an open problem whether there exist II$_1$ factors with
separable predual and uncountable fundamental group different from $\R^*_+$.

\section{\mbox{From von Neumann equivalence to orbit equivalence}} \label{sec.Neumann-to-orbit}

The following is an immediate consequence of Theorem
\ref{thm.rigidity}.

\begin{proposition} \label{prop.intertwine-group}
Let $G$ be an ICC group with a strongly malleable mixing action on the probability space $(X,\mu)$. Write $M = L^\infty(X) \rtimes G$. Let
$\Gamma$ be a countable group that admits an almost normal infinite subgroup $\Gamma_0$ such that $(\Gamma,\Gamma_0)$ has the relative property (T).
Suppose that $\Gamma$ acts on the probability space $(Y,\eta)$.

Let $p$ be a projection in $\cL(G)$ and
$$\theta : L^\infty(Y) \rtimes \Gamma \recht p(L^\infty(X) \rtimes G)p$$
a $^*$-isomorphism. Then, there exists a unitary $v \in pMp$ such that $v \theta(\cL(\Gamma)) v^* \subset p\cL(G)p$.
\end{proposition}
\begin{proof}
We apply Theorem \ref{thm.rigidity}, observing that $\cL(\Gamma)$ is included in the quasi-\linebreak normalizer $P$ of $\cL(\Gamma_0)$ in
$L^\infty(Y) \rtimes \Gamma$. Using Remark \ref{rem.amenable}, it
follows that there is no non-zero homomorphism from $P$ to an amplification of $L^\infty(X)$.
\end{proof}

From now on, specify $G \actson (X,\mu)$ to be the \emph{Bernoulli action}. The
following preliminary result is proved: an isomorphism between crossed
products sending one group algebra into the other, makes the Cartan
subalgebras conjugate. The final aim is Theorem
\ref{thm.strong-rigidity} below, which states that the actions are
necessarily conjugate.

\begin{theorem}[Popa, \cite{P2}] \label{thm.intertwine-cartan}
Let $G$ be an infinite group and, for $\mu_0$ non-atomic, $G \actson (X,\mu)= \cartesian{g \in
  G} (X_0,\mu_0)$, its Bernoulli action.
Let $\Gamma$ be an infinite group
that acts freely and weakly mixingly on the probability space $(Y,\eta)$. Write $A=L^\infty(X)$ and $B=L^\infty(Y)$. Let $p$ be a projection in
$\cL(G)$ and
$$\theta : B \rtimes \Gamma \recht p(A \rtimes G)p$$
a $^*$-isomorphism. Suppose that $\theta(\cL(\Gamma)) \subset p \cL(G)
p$. Then,
\begin{itemize}
\item there exists a partial isometry $u \in A \rtimes G$ satisfying $u^*u
= p$, $e:=uu^* \in A$ and $u \theta(B) u^* = eA$;
\item the equality $\theta(\cL(\Gamma)) = p \cL(G) p$ holds.
\end{itemize}
\end{theorem}

Later on, Proposition \ref{prop.intertwine-group} and Theorem
\ref{thm.intertwine-cartan} are combined to prove that the actions of $\Gamma$ and
$G$ are conjugate through a group isomorphism of $\Gamma$ and $G$. The proof of Theorem \ref{thm.intertwine-cartan} certainly is the most technical
and analytically subtle part of this talk.

\begin{notations}
We fix several notations used throughout the lemmas needed to prove Theorem \ref{thm.intertwine-cartan}.
\begin{itemize}
\item We fix an infinite group $G$ and write $A_0= L^\infty(X_0)$, $\dis (A,\tau)=\bigotimes_{g \in G} (A_0,\tau_0)$. For every finite subset $K \subset G$, we write $\dis A_{K^c} := \bigotimes_{g
\notin K} (A_0,\tau_0)$. Write $M = A \rtimes G$ and denote by $\tau$
the tracial state on $M$.
\item We use $\eta : M \recht L^2(M)$ to identify an element of the algebra $M$ with its corresponding vector in the Hilbert space $L^2(M)$.
\item For a finite subset $K \subset G$, we denote by $e_{\check{K}}$ the orthogonal projection onto the closure of $\lspan\{ \eta(A_{K^c} u_g) \mid g \in G \}$ in $L^2(M)$ and
we denote by $p_{\check{K}}$ the orthogonal projection onto the closure of $\lspan\{ \eta(A u_k) \mid k \in G \setminus K\}$ in $L^2(M)$.
\item We do not write the isomorphism $\theta$. We simply suppose that
  $B \rtimes \Gamma = p(A \rtimes G)p$ in such a way that $\cL(\Gamma)
  \subset p\cL(G)p$. Of course, $\tau$ is as well the trace on $B
  \rtimes \Gamma$, but non-normalized.
\item The elements of $\Gamma$ are denoted by $s,t$ and the action of
  $\Gamma$ on $B$ by $(\rho_s)_{s \in \Gamma}$. The elements of $G$
  are denoted by $g,h$ and the action of $G$ on $A$ by $(\si_g)_{g \in
  G}$.
\item Denote by $(\ups_s)_{s
  \in \Gamma}$ the canonical unitaries generating $\cL(\Gamma)$ and by
  $(u_g)_{g \in G}$ the canonical unitaries generating $\cL(G)$.
\end{itemize}
\end{notations}

We first explain the idea of the proof of Theorem
\ref{thm.intertwine-cartan}. Elements in the image of $e_{\check{K}}$ for $K$
large are thought of as \emph{living far away space-wise}, while
elements in the image of $p_{\check{K}}$ for $K$ large are thought of as
\emph{living far away group-wise}. In order to show that $B$ can be
conjugated into $A$, one shows first that sufficiently many elements of $B$
are not living far away group-wise. This suffices to construct a
$B$-$A$-subbimodule of $L^2(M)$ which is finitely generated as an $A$-module. To obtain elements of $B$
that are not living far away group-wise, two lemmas are used:
\begin{itemize}
\item if an element of $B$ lives far away space-wise, it does not live
  far away group wise (Lemma \ref{lem.one});
\item if $b \in B$ and $s_n \recht \infty$ in $\Gamma$, the elements
  $\rho_{s_n}(b)$ are more and more living far away space-wise (Lemma \ref{lem.two}).
\end{itemize}
To pass from the approximate inequalities in Lemmas
\ref{lem.one}, \ref{lem.two} to exact inequalities, the powerful
technique of ultraproducts is applied. This allows to conjugate $B$ into
$A$ at least on the level of the ultrapower algebra. But this is
sufficient to return to earth and conjugate $B$ into $A$.

\begin{lemma} \label{lem.one}
For every $\eps>0$ there exist finite subsets $K,L \subset G$ such
that
$$\|p_{\check{K}} \eta(x) \|^2 \leq 3 \|(1-e_{\check{L}})\eta(x)\| + \eps$$
for all $x \in B$ with $\|x\|\leq 1$.
\end{lemma}
\begin{proof}
We make the following claim. \\ {\it Claim.} For every $a
\in M$ with $\|a\|\leq 1$ and every $\eps > 0$, there exist $K,L
\subset G$ finite such that
$$|\la a \cdot \eta(x) \cdot a^* , p_{\check{K}} \eta(x) \ra| \leq 3
\|(1-e_{\check{L}})\eta(x)\| + \|E_{\cL(G)}(a)\|_2 + \eps$$
for all $x \in M$ with $\|x\|\leq 1$. To deduce the lemma from this
claim it is then sufficient to prove that $B$ contains unitaries $a$
with $\|E_{\cL(G)}(a)\|_2$ arbitrary small and to use the
commutativity of $B$ in order to get $a \cdot \eta(x) \cdot a^* =
\eta(x)$ for $x \in B$.

To prove the claim, choose $a \in M$ with $\|a\|\leq 1$ and $\eps >
0$. By the Kaplansky density theorem, we may assume that $a \in
\lspan\{ A_{F_0} u_g \mid g \in F_1\}$ for some finite subsets
$F_0,F_1 \subset G$. We may assume as well that $e \in F_1$. Put $L=F_1^{-1}F_0$ and $K=LF_0^{-1}$. It is an
excellent Bernoulli exercise to check that
$$e_{\check{L}}(a \cdot \xi) = e_{\check{L}}(E_{\cL(G)}(a) \cdot \xi) \;\; \text{for $\xi
  \in \Im e_{\check{L}}$,}\qquad e_{\check{L}}(\xi \cdot a) = (e_{\check{L}} \xi) \cdot a \;\; \text{for $\xi
  \in \Im p_{\check{K}}$.}$$
Take $x \in M$ with $\|x\|\leq 1$. We obtain that
\begin{equation}
|\la a \cdot \eta(x) \cdot a^* , p_{\check{K}} \eta(x) \ra|
\leq \|e_{\check{L}}(a \cdot \eta(x))\| +
\|(1-e_{\check{L}})\bigl((p_{\check{K}} \eta(x)) \cdot a\bigr)\| \; .
\tag{$\ast$}
\end{equation}
In ($\ast$), the second term equals
$$\|\bigl((1-e_{\check{L}})p_{\check{K}} \eta(x)\bigr)\cdot a\| \leq
\|(1-e_{\check{L}})\eta(x)\| \; .$$
The first term of ($\ast$), is bounded by
\begin{equation}
\|e_{\check{L}}(a \cdot (e_{\check{L}} \eta(x))) \| +
\|(1-e_{\check{L}})\eta(x)\| \; . \tag{$\ast\ast$}
\end{equation}
In ($\ast\ast$), the first term equals
\begin{align*}
\|e_{\check{L}}( E_{\cL(G)}(a) \cdot (e_{\check{L}} \eta(x))) \|
& \leq \|E_{\cL(G)}(a) \cdot (e_{\check{L}} \eta(x))\| \\ &\leq
\|E_{\cL(G)}(a) \cdot \eta(x)\| + \|(1-e_{\check{L}})\eta(x)\| \\ &\leq
\|E_{\cL(G)}(a)\|_2 + \|(1-e_{\check{L}})\eta(x)\| \;.
\end{align*}
We have shown that
$$|\la a \cdot \eta(x) \cdot a^* , p_{\check{K}} \eta(x) \ra| \leq 3
\|(1-e_{\check{L}})\eta(x)\| + \|E_{\cL(G)}(a)\|_2$$
for all $x \in M$ with $\|x\| \leq 1$, which proves the claim.

It remains to show that, for every $\eps >0$, there exists a unitary
$u \in B$ such that $\|E_{\cL(G)}(u)\|_2 < \eps$. If not, it follows
from Proposition \ref{prop.translate} that there exists $n \geq 1$, a
projection $q \in \M_n(\C) \ot \cL(G)$, a homomorphism $\theta : B
\recht q (\M_n(\C) \ot \cL(G))q$ and a non-zero partial isometry $v
\in \M_{1,n}(\C) \ot p M$ satisfying $v^*v \leq q$ and $b v = v
\theta(b)$ for all $b \in B$. Using Theorem \ref{thm.mixing-one},
$v^*v \in \M_n(\C) \ot \cL(G)$ and we may assume that $v^*v=q$. Then,
$v^*Bv$ is a diffuse subalgebra of $q (\M_n(\C) \ot \cL(G))q$. Since
the normalizer of $B$ in $pMp$ is the whole of $pMp$, it follows from
Theorem \ref{thm.mixing-one} that $v^*Mv \subset q (\M_n(\C) \ot
\cL(G))q$. Since $v^*Mv = q (\M_n(\C) \ot M)q$, this is a contradiction.
\end{proof}

\begin{lemma}\label{lem.two}
For every $b \in B$, $\eps > 0$ and $L \subset G$ finite, there exists
$K \subset \Gamma$ finite such that
$$\|(1-e_{\check{L}})\eta(\rho_s(b))\| < \eps$$
for all $s \in \Gamma \setminus K$.
\end{lemma}
\begin{proof}
We again make a claim.\\ {\it Claim.} For every $a \in M$ with
$\|a\|\leq 1$, $L \subset G$ finite and $\eps > 0$, there exists $K_1
\subset G$ finite such that
$$\| (1-e_{\check{L}}) \eta(v a w) \| \leq \eps + \|(1-p_{\check{K}_1})\eta(v)\|$$
for all $v,w \in \cL(G)$ with $\|w\| \leq 1$.

The lemma follows easily from the claim: given $K_1
  \subset G$ finite and $\eps > 0$, we can take $K
  \subset \Gamma$ finite such that $\|(1-p_{\check{K}_1})\eta(\ups_s)\| < \eps$ for all $s \in \Gamma \setminus K$. It remains to
  observe that $\rho_s(b) = \ups_s b \ups_s^*$ and $\ups_s \in
  \cL(\Gamma) \subset \cL(G)$.

To prove the claim, choose $a \in M$ with $\|a\|\leq 1$ and $\eps >
0$. By the Kaplansky density theorem, we may assume that $a \in
\lspan\{ A_{F} u_g \mid g \in G\}$ for some finite subset
$F \subset G$. Given $L \subset G$ finite, we put $K_1=L F^{-1}$ and
leave as an exercise to check that
$$(p_{\check{K}_1}\eta(v)) \cdot (aw) \in \Im e_{\check{L}} \;\;\text{for all $v,w \in \cL(G)$.}$$
The claim follows immediately.
\end{proof}

\begin{lemma} \label{lem.conditional}
For every $b \in B$, $\dis E_{\cL(G)}(b) =
\frac{\tau(b)}{\tau(p)}p$. Hence, $\cL(\Gamma) = p\cL(G)p$.
\end{lemma}
\begin{proof}
We have to prove the following: if $b \in B$ and $\tau(b) = 0$, then
$E_{\cL(G)}(b) = 0$. Take such a $b \in B$ with $\tau(b)=0$. Since
$\Gamma$ acts weakly mixingly on $B$, we take a sequence $s_n \recht
\infty$ in $\Gamma$ such that $\rho_{s_n}(b) \recht 0$ in the weak topology.

Combining Lemmas \ref{lem.one} and \ref{lem.two}, we find a finite subset $K \subset G$ and $n_0$ such that $\|p_{\check{K}} \eta(\rho_{s_n}(b))\|^2 \leq \eps$
for all $n \geq n_0$. Denote by $f$ the orthogonal projection of $L^2(M)$ onto the closure of $\eta(\cL(G))$. Since $f$ and $p_{\check{K}}$ commute, we find
that $\|p_{\check{K}} \eta\bigl(E_{\cL(G)}(\rho_{s_n}(b))\bigr)\|^2 \leq \eps$ for all $n \geq n_0$. On the other hand, $E_{\cL(G)}(\rho_{s_n}(b))$ tends
weakly to $0$ and belongs to $\cL(G)$. Hence,
$$\|(1-p_{\check{K}}) \eta\bigl(E_{\cL(G)}(\rho_{s_n}(b))\bigr)\|^2 \recht 0$$
when $n \recht \infty$. We conclude that for $n$ sufficiently large,
$\|E_{\cL(G)}(\rho_{s_n}(b))\|^2_2 \leq 2\eps$. But, for every $n$,
$$\|E_{\cL(G)}(\rho_{s_n}(b))\|_2=
\|\ups_{s_n}E_{\cL(G)}(b)\ups_{s_n}\|_2 = \|E_{\cL(G)}(b)\|_2 \; .$$
It follows that $\|E_{\cL(G)}(b)\|_2^2 \leq 2\eps$ for all $\eps > 0$,
which proves that $E_{\cL(G)}(b) = 0$.

Since $pMp = B \rtimes \Gamma$ and $\cL(\Gamma) \subset p\cL(G)p$, it
suffices to apply $E_{\cL(G)}$ to obtain that $p\cL(G)p=\cL(\Gamma)$.
\end{proof}

Let us warm up the ultraproduct machinery to finish the proof of
Theorem \ref{thm.intertwine-cartan}.

\begin{notations}
We introduce the following notations.
\begin{itemize}
\item Let $\om$ be a free ultrafilter on $\N$ and define the ultrapower
algebra $M^\om$, containing $A^\om$ as a maximal abelian subalgebra. Denote by $A^\om_\infty \subset A^\om$ the \emph{tail algebra} for the Bernoulli
action, defined as
$$A^\om_\infty := \bigcap_{\begin{smallmatrix}F \subset G \\ F \;
    \text{finite}\end{smallmatrix}}
(A_{F^c})^\om \; .$$ Observe that $A^\om_\infty$, as a subalgebra of $M^\om$ is normalized
    by the unitaries $(u_g)_{g \in G}$.
\item Denote by $A^\om_\infty \rtimes G$ the von Neumann subalgebra of
  $M^\om$ generated by $A^\om_\infty$ and $\cL(G)$.
\item We define
$\chi := B^\om \cap p(A^\om_\infty \rtimes G)p$.
\end{itemize}
\end{notations}

Lemmas \ref{lem.one} and \ref{lem.two} can be reinterpreted to yield elements of $\chi$.

\begin{lemma} \label{lem.hulpresult}
The following results hold.
\begin{enumerate}
\item\label{res.one} A bounded sequence $(b_n)$ in $B$ represents an element of $\chi$ if and only if
$$
\lim_{n \recht \om} \|(1-e_{\check{L}})\eta(b_n)\| = 0 \quad\text{for every finite subset $L \subset G$.}
$$
\item\label{res.two} When $s_n \recht \infty$ in $\Gamma$ and $b \in B$, the sequence $(\rho_{s_n}(b))$ represents an element in $\chi$.
\item\label{res.three} If a bounded sequence $(b_n)$ in $B$ represents an element of $\chi$, then $b_n - \tau_1(b_n)p$ tends to $0$ weakly.
Here $\tau_1:=\tau(\cdot)/\tau(p)$ denotes the normalized trace on $pMp$.
\end{enumerate}
\end{lemma}
\begin{proof}
\eqref{res.one} If $(a_n) \in A^\om_\infty$ and $g \in G$, clearly $\lim_{n \recht \om}  \|(1-e_{\check{L}})\eta(a_n u_g)\| = 0$. Hence, the same holds if we
replace $(a_n u_g)$ by any element of $A^\om_\infty \rtimes G$. Conversely, let $b \in B^\om$ be represented by the bounded sequence $(b_n)$ in $B$
such that \eqref{res.one} holds. For any finite $K \subset G$, define $z_K \in M^\om$ by the sequence $\bigl( \sum_{g \in K} E_A(b_n u_g^*) u_g
\bigr)$. Our assumption yields that $z_K \in A^\om_\infty \rtimes G$ for all $K$. From Lemma \ref{lem.one} it follows that $\|z_K - b\|_2 \recht 0$,
if $K \recht G$. Hence, $b \in A^\om_\infty \rtimes G$.

\eqref{res.two} This follows using Lemma \ref{lem.two} and statement \eqref{res.one}.

\eqref{res.three} Using Lemma \ref{lem.conditional}, it suffices to check that $b_n - E_{\cL(G)}(b_n)$ tends to $0$ weakly. This is true for any
$(b_n)$ in $A^\om_\infty \rtimes G$.
\end{proof}

In the
next lemma, $\chi$ is shown to be sufficiently big.

\begin{lemma} \label{lem.chi-commutant}
One has $pM^\om p \cap \chi' = B^\om$.
\end{lemma}
\begin{proof}
We first claim that the action $(\rho_s)_{s \in \Gamma}$ is $2$-mixing
(see Definition \ref{def.two-mixing}). We have to prove that for all
$a,b,c \in B$,
$$|\tau(a \rho_s(b) \rho_t(c)) - \tau(a)\tau(\rho_s(b)\rho_t(c))|
\recht 0$$
when $s,t \recht \infty$.

Suppose that the bounded sequence $(d_n)$ represents an element $d \in \chi$. By \eqref{res.three} in Lemma \ref{lem.hulpresult}, $d_n - \tau_1(d_n)p
\recht 0$ weakly and hence,
$$|\tau_1(a d_n) - \tau_1(a) \tau_1(d_n)| \recht 0$$
for all $a \in B$. Fix $a,b,c \in B$ and take sequences $s_n,t_n
\recht \infty$ in $\Gamma$. From \eqref{res.two} in Lemma \ref{lem.hulpresult}, we get that the
sequences $(\rho_{s_n}(b))$ and $(\rho_{t_n}(c))$ represent elements of $\chi$. Since $\chi$ is a von Neumann algebra, the sequence
$(\rho_{s_n}(b)\rho_{t_n}(c))$ represents an element of $\chi$ as well. Applying the previous paragraph to this sequence, we have proved the claim.
Combining the $2$-mixing of the action $(\rho_s)_{s \in \Gamma}$ with Lemma \ref{lem.two-mixing}, we are done.
\end{proof}

\begin{proof}[Proof of Theorem \ref{thm.intertwine-cartan}]
We first claim there exists a non-zero $a \in p\la M^\om,e_{A^\om}\ra^+p \cap \chi'$ with $\tauh(a) < \infty$. As usual, $\tauh$ denotes the
semi-finite trace on the basic construction $\la M^\om,e_{A^\om}\ra$,
see Appendix \ref{sec.basic}.

There exists a finite subset $K \subset G$ such that
$$\lim_{n \recht \om} \|p_{\check{K}} \eta(b_n)\| \leq \frac{1}{2}$$
for all $(b_n)$ in the unit ball of $\chi$. Indeed, if not, write $G$ as an increasing union of finite subsets $K_n$ and choose $b_n \in B$ with
$\|b_n\|\leq 1$, $\|(1-e_{\check{K}_n})\eta(b_n)\|\leq 1/n$ and $\|p_{\check{K}_n} \eta(b_n)\| > 1/2$. This yields a contradiction with Lemma \ref{lem.one}.

Define the projection $f_K \in \la M^\om,e_{A^\om} \ra$ as $f_K = \sum_{g \in K} u_g^* e_{A^\om} u_g$. Clearly $\tauh(f_K) < \infty$. Denote by $a$
the (unique) element in the ultraweakly closed convex hull of $\{b f_K b^* \mid b \in \cU(\chi) \}$. By construction $\tauh(a) < \infty$ and $a \in
\chi'$. To obtain the claim, we have to show that $a \neq 0$. Whenever $(b_n)$ represents $b \in \cU(\chi)$, we have
$$\tauh(e_{A^\om} b f_K b^* e_{A^\om}) = \lim_{n \recht \om} \|(1-p_{\check{K}})\eta(b_n)\|^2 \geq 3/4 \; .$$
Hence, $\tauh(e_{A^\om} a e_{A^\om}) \neq 0$ and $a \neq 0$. This proves the claim stated in the beginning of the proof.

It follows from Lemma \ref{lem.chi-commutant} and Theorem \ref{thm.intertwine-abelian} that there exists a non-zero partial isometry $v \in M^\om$
satisfying $v^*v \in B^\om$, $vv^* \in A^\om$ and $v B^\om v^* \subset A^\om$. Take partial isometries $v_n \in M$ such that $e_n:=v_n^* v_n \in B$,
$v_n v_n^* \in A$ and $(v_n)$ represents $v$. It follows that there exists $n$ such that
$$\|v_n b v_n^* - E_A(v_n b v_n^*)\|_2 < \frac{1}{2}\|e_n\|_2$$
for all $b \in B$ with $\|b\|\leq 1$. Indeed, if not, we find a sequence of elements $b_n \in B$ with $\|b_n\|\leq 1$ and $\|v_n b_n v_n^* - E_A(v_n b
v_n^*)\|_2 \geq \frac{1}{2}\|e_n\|_2$. Since $(b_n)$ defines an element in $B^\om$, taking the limit $n \recht \om$ yields a contradiction.

If we write $f=v_n v_n^* \in A$, $A_1 := fA$ and $B_1:= v_nBv_n^*$ as subalgebras of $fMf$, we have, after normalization of the trace, $\|b -
E_{A_1}(b)\|_2 \leq \frac{1}{2}$ for all $b \in B_1$ with $\|b\|\leq 1$. Hence, \eqref{stat.four} in Proposition \ref{prop.translate} is satisfied
and an application of Theorem \ref{thm.intertwine-abelian} concludes the proof of Theorem \ref{thm.intertwine-cartan}.
\end{proof}

\section{Strong rigidity for von Neumann algebras} \label{sec.strong-rigidity}

Suppose that $G$ acts on $(A,\tau)$ by $(\si_g)_{g \in G}$ and $\Gamma$ on $(B,\tau)$ by $(\rho_s)_{s \in \Gamma}$. A \emph{conjugation} of both
actions is a pair $(\Delta,\delta)$ of isomorphisms $\Delta : B \recht A$, $\delta : \Gamma \recht G$ satisfying
$\Delta(\rho_s(b))=\si_{\delta(s)}(\Delta(b))$, for all $b \in B$ and $s \in \Gamma$. Associated with the conjugation $(\Delta,\delta)$ is of course
the obvious isomorphism of crossed products $\theta_{\Delta,\delta} : B \rtimes \Gamma \recht A \rtimes G$.

Whenever $G$ acts on $(A,\tau)$ and $\alpha : G \recht S^1$ is a character, we have an obvious automorphism $\theta_\alpha$ of the crossed product $A
\rtimes G$ defined as fixing pointwise $A$ and $\theta_\alpha(u_g) = \alpha(g) u_g$.

\begin{theorem}[Popa, \cite{P2}] \label{thm.strong-rigidity}
Let $G$ be an ICC group acting and $G \actson (X,\mu)$ its Bernoulli action (with non-atomic
base). Let $\Gamma$ be a countable group that admits an almost normal
infinite subgroup $\Gamma_0$ such that $(\Gamma,\Gamma_0)$ has the relative property (T). Suppose that $\Gamma$ acts freely on the probability space
$(Y,\eta)$. Let $p$ be a projection in $L^\infty(X) \rtimes G$ and
$$\theta : L^\infty(Y) \rtimes \Gamma \recht p(L^\infty(X) \rtimes G)p$$
a $^*$-isomorphism. Then, $p=1$ and there exist a unitary $u \in L^\infty(X) \rtimes G$, a conjugation $(\Delta,\delta)$ of the actions through a
group isomorphism $\delta : \Gamma \recht G$ and a character $\alpha$ on $G$ such that
$$\theta = \Ad u \circ \theta_\alpha \circ \theta_{\Delta,\delta} \; .$$
\end{theorem}

Theorem \ref{thm.strong-rigidity} admits the following corollary stated in the
introduction.

\begin{corollary}
Let $G$ be a $w$-rigid group and denote by $M_G := L^\infty(X) \rtimes
G$ the crossed product of the Bernoulli action $G \actson (X,\mu)$ with non-atomic base. Then, for $w$-rigid ICC groups $G$ and
$\Gamma$, we have $M_G \cong M_\Gamma$ if and only if $G \cong
\Gamma$. Moreover, all $M_G$ for $G$ $w$-rigid ICC, have trivial
fundamental group.
\end{corollary}

The corollary is an immediate consequence of Theorem \ref{thm.intertwine-cartan} and the orbit equivalence superrigidity Theorem
\ref{thm.OEsuperrigidity}. Indeed, let $G$ and $\Gamma$ be $w$-rigid ICC groups with Bernoulli actions on $(X,\mu)$ and $(Y,\eta)$,
respectively. If $p$ is a projection in $L^\infty(X) \rtimes G$ and $\theta : L^\infty(Y) \rtimes \Gamma \recht p (L^\infty(X) \rtimes G) p$ is a
$^*$-isomorphism, we have to prove that $p=1$ and that $\Gamma$ and $G$ are isomorphic. Combining Proposition \ref{prop.intertwine-group} and Theorem
\ref{thm.intertwine-cartan}, we may assume that $p \in L^\infty(X)$ and $\theta(L^\infty(Y)) = L^\infty(X)p$. Hence, $\theta$ defines a stable orbit
equivalence between $\Gamma \actson Y$ and $G \actson X$. So, Theorem \ref{thm.OEsuperrigidity} allows to conclude.

Refining the reasoning above, Theorem \ref{thm.strong-rigidity} is proved. First, taking a further reduction, it is shown that we may assume that the
action $\Gamma \actson Y$ is weakly mixing. So, Proposition \ref{prop.intertwine-group} and Theorem \ref{thm.intertwine-cartan} can be applied and
yield a stable orbit equivalence of $\Gamma \actson Y$ and $G \actson X$. Associated with this stable orbit equivalence is a cocycle. The unitary
that conjugates $\cL(\Gamma)$ into $\cL(G)$ (its existence is the contents of Proposition \ref{prop.intertwine-group}) is reinterpreted as making
cohomologous this cocycle to a homomorphism into $\cU(\cL(G))$. Using
the weak mixing property through an application of Lemma \ref{lem.reduce-subgroup}, the
homomorphism can be assumed to take values in $G$ itself. This yields the conjugacy of the actions.

\begin{proof}[Proof of Theorem \ref{thm.strong-rigidity}]
Write $A=L^\infty(X)$ and $B=L^\infty(Y)$. Write $M=A \rtimes G$ and identify through $\theta$, $B \rtimes \Gamma = p(A \rtimes G)p$. First applying
Proposition \ref{prop.intertwine-group}, we may assume that $p \in
\cL(G)$ and $\cL(\Gamma) \subset p \cL(G) p$. We claim that there
exists a finite index subgroup $\Gamma_1 \subset \Gamma$ and a
$\Gamma_1$-invariant projection $p_1 \in B \cap \cL(G)$ such that the
$\Gamma$-action on $B$ is induced from the $\Gamma_1$-action on $p_1
B$ obtained by restriction, and such that the $\Gamma_1$-action on $p_1 B$
is weakly mixing.

Whenever $V \subset B$ is a finite-dimensional $\Gamma$-invariant subspace, it follows from Theorem \ref{thm.mixing-one} that
$V \subset p\cL(G)p$. Also, $B \cap \cL(G)$ is a $\Gamma$-invariant von Neumann subalgebra of $B$. By the ergodicity of the $\Gamma$-action on $B$,
this invariant subalgebra is either diffuse or atomic. If it is diffuse and since it commutes with $B$, it would follow from Theorem
\ref{thm.mixing-one} that $B \subset p\cL(G)p$ and hence, $pMp \subset p\cL(G)p$, a contradiction. So, $B \cap \cL(G)$ is atomic, hence
finite-dimensional, and it suffices to take a minimal projection $p_1
\in B \cap \cL(G)$. This proves the claim.

It now suffices to prove the theorem under the additional assumption
that the action of $\Gamma$ on $B$ is weakly mixing. We apply Theorem \ref{thm.intertwine-cartan}. Conjugating again, we obtain the following situation: a projection $q \in A$
and a partial isometry $v \in M$ such that $vv^*=p \in \cL(G)$,
$v^*v=q$ and $B \rtimes \Gamma = q(A \rtimes G)q$ in such a way that
$B = qA$ and $v \cL(\Gamma) v^* = p \cL(G) p$. The theorem follows from
Proposition \ref{prop.nes-sto} below.
\end{proof}

In the proof of Theorem \ref{thm.strong-rigidity}, we used the
following proposition. It is a weaker version of Theorem 5.2 in
\cite{P2}, but sufficient for our purposes. It also
provides a generalization
and simpler proof for the main result in \cite{NS} by Neshveyev and St{\o}rmer.

\begin{proposition}[Popa, \cite{P2}] \label{prop.nes-sto}
Let $G$ be an infinite group that acts freely and weakly mixingly on
$(X,\mu)$. Let $\Gamma$ be an infinite group that acts freely on
$(Y,\eta)$. Write $A =L^\infty(X)$ and $B = L^\infty(Y)$. Suppose
that $q$ is a projection in $A$ such that
$$B \rtimes \Gamma = q(A \rtimes G)q \quad\text{with}\quad B = qA \;
.$$
Suppose that there exists a partial isometry $v \in A \rtimes G$
satisfying $v^* v = q$, $vv^* = p \in \cL(G)$ and $v \cL(\Gamma) v^* =
p\cL(G)p$.
\begin{itemize}
\item If $G$ has no finite normal subgroups, $q=1$.
\item If $q=1$, there exists $w \in \cU(\cL(G))$ such that, writing
  $\vtil = wv$, $\vtil$ normalizes $B=A$ and $\vtil \ups_s \vtil^* =
  \alpha(s) u_{\delta(s)}$ for some $\alpha \in \Char(\Gamma)$ and
  some group isomorphism $\delta : \Gamma \recht G$.
\end{itemize}
\end{proposition}
We write this rather pedantic formulation of the proposition, to cover
at the same time its application in the proof of Theorem
\ref{thm.strong-rigidity} (where $G$ is ICC and hence, without finite
normal subgroups) and the result of \cite{NS} (where $G$ is an any abelian
group, but $q=1$ from the beginning).

\begin{proof}
We make use of the canonical embedding $\eta: A \rtimes G \recht A \ovt \ell^2(G)$ of the crossed product into the Hilbert-W$^*$-module $A \ovt
\ell^2(G)$ given by $\eta(u_g a) = a \ot \delta_{g^{-1}}$ for all $g \in G$ and $a \in A$. Here $(\delta_g)_{g \in G}$ is the canonical orthonormal
basis of $\ell^2(G)$. We identify $A \ovt \ell^2(G) = L^\infty(X,\ell^2(G))$ and we make act $\cL(G)$ on $\ell^2(G)$ on the left and the right: $u_g
\delta_h = \delta_{gh}$ and $\delta_h u_g = \delta_{hg}$. At the same time, we regard $\cL(G) \subset \ell^2(G)$.

Denote $S^1 G:= S^1 \times G$ that we identify in the obvious way with a closed subgroup of $\cU(\cL(G))$. We identified $Y \subset X$ such that
$\Gamma$ acts on $Y$, $B = L^\infty(Y)$, $A=L^\infty(X)$ and $q = \chi_Y$.
We have the orbit equivalence $q(A \rtimes G)q = B \rtimes \Gamma$ with $B = qA$. This yields a one-cocycle $\gamma : \Gamma
\times Y \recht S^1 G$ given by
$$\eta(z \ups_s)(x) = \eta(z)(s \ast x) \; \gamma(s,x)$$
for all $z \in A \rtimes G$ and where we use $s \ast x$ to denote the action of an element $s \in \Gamma$ on $x \in Y$. We claim that the partial
isometry $v$ makes $\gamma$ cohomologous to a homomorphism.

Observe that $E_{\cL(G)}(v a v^*) = \tau(p)^{-1} \tau(a) p$ for all $a
\in B$. Indeed,
$$E_{\cL(G)}(v a v^*) = \tau(p)^{-1} E_{v\cL(\Gamma)v^*}(v a v^*) =
\tau(p)^{-1} v E_{\cL(\Gamma)}(a) v^* = \tau(p)^{-1} \tau(a) p \; .$$ We first study the function $w:=\tau(p)^{1/2}\eta(v) \in
L^\infty(Y,\ell^2(G))$. Suppose that $w_0 \in \cL(G)$ is an essential value of this function. We find a decreasing sequence of non-zero projections
$q_n$ in $B$ such that $\|\tau(p)^{1/2}\eta(v)q_n - q_n \ot w_0\|_\infty \recht 0$, where we use the uniform norm for functions in
$L^\infty(X,\ell^2(G))$. So, we have a sequence $\eps_n \recht 0$ such that $\|(\tau(p)^{1/2}v-w_0)q_n\|_2 \leq \eps_n \|q_n\|_2$, where we use the
norm of $L^2(M)$. In $L^1(M)$, we obtain that $\tau(q_n)^{-1} \|\tau(p) v q_n v^* - w_0 q_n w_0^* \|_1 \recht 0$. Applying $E_{\cL(G)}$ it follows
that $\| p - w_0 w_0^*\|_1 \recht 0$ and hence $w_0 w_0^* = p$. We have shown that for almost all $y \in Y$,
$$w(y) \in \cL(G) \quad\text{and}\quad w(y) w(y)^* = p \; .$$
Since we can replace $v$ by $w_0^* v$, we may assume that $p$ is an essential value of the function $w$.

Define the homomorphism $\pi : \Gamma \recht \cU(p\cL(G)p) : \pi(s) = v \ups_s v^*$. For every $s \in \Gamma$, $v \ups_s = \pi(s) v$. Applying
$\eta$, this yields,
\begin{equation}\label{eq.cocycletrivial}
w(s \ast x) \; \gamma(s,x) = \pi(s) \; w(x) \quad\text{for almost all $x \in Y$.}
\end{equation}
If $q=1$, Lemma \ref{lem.reduce-subgroup} yields that $\pi(s) \in S^1 G$ for all $s \in \Gamma$ and $w(x) \in S^1 G$ for almost all $x \in X$. The latter implies
that $v$ normalizes the Cartan subalgebra $A=B$. The former allows to define the group isomorphism $\delta : \Gamma \recht G$ and the character
$\alpha : \Gamma \recht S^1$ such that $\pi(s) = \alpha(s) \delta(s)$
for all $s \in \Gamma$. So, we are done in the case $q=1$.

It remains to show that $p=1$ when $G$ has no finite normal subgroups. The orbit equivalence allows as well for an inverse $1$-cocycle: define $W =
\{(g,x) \in G \times Y \mid x \in Y, g \cdot x \in Y \}$. We use the notation $g \cdot x$ to denote the action of an element $g \in G$ on $x \in X$.
Then, the $1$-cocycle $\mu : W \recht S^1 \Gamma$ is well defined and related to $\gamma$ by the formula
$$g = \gamma(\mu_{\text{\rm group}}(g,x),x)\; \mu_{\text{\rm scal}}(g,x)$$
for almost all $(g,x) \in W$. Here we split up explicitly $\mu = \mu_{\text{\rm scal}} \mu_{\text{\rm group}}$. Plugging the previous equality into
\eqref{eq.cocycletrivial} yields
\begin{equation}\label{eq.cocycletrivial-two}
w(g \cdot x)\; u_g = \pi(\mu(g,x)) \; w(x) \quad\text{for almost all $(g,x) \in W$.}
\end{equation}
Since $p$ is an essential value of the function $w$ and since $\pi$ takes values in the unitaries of $p \cL(G) p$, arguing exactly as in the proof of
Lemma \ref{lem.reduce-subgroup}, yields that for any $g \in G$, $p u_g$ is arbitrary close to a unitary and hence, $u_g$ and $p$ commute for all $g
\in G$. So, $p$ is a central projection in $\cL(G)$ and it follows that $w(x) \in \cU(p\cL(G)p)$ for almost all $x \in Y$. Conjugating equation
\eqref{eq.cocycletrivial-two} with $v^*$, implies that the cocycle $\mu : W \recht S^1 \Gamma$ is cohomologous, as a cocycle with values in
$\cU(\cL(\Gamma))$, to the homomorphism $g \mapsto v^* u_g v$. It follows from Lemma \ref{lem.reduce-subgroup} that $v^* u_g v \in S^1 \Gamma$ for
all $g \in G$. On $S^1 \Gamma$, the trace $\tau$ takes the values $0$
and $\tau(p) S^1$. Hence, for all $g \in G$, we have
$$\tau(u_g p) = \tau(u_g vv^*) = \tau(v^* u_g v) \in \{0\} \cup S^1
\tau(p) \; .$$ We also know that $p$ is a central projection in $\cL(G)$.
It is an excellent
exercice to deduce from all this that $p$ is of the form $\sum_{g \in K} \beta(k) u_k$ for some finite normal subgroup $K \subset G$ and an $\Ad
G$-invariant character $\beta \in \Char K$. Hence, $K = \{e\}$, $p=1$ and we are done.
\end{proof}

\section{Outer conjugacy of $w$-rigid group actions on the
  hyperfinite II$_1$ factor}

We discuss some of the results of Popa \cite{P4} on (cocycle) actions of $w$-rigid groups on the hyperfinite II$_1$ factor. As explained in the
introduction, the paper \cite{P4} is the precursor to all of Popa's papers on rigidity of Bernoulli actions.

\begin{definition}
A \emph{cocycle action} of a countable group $G$ on a von Neumann
algebra $N$ consists of automorphisms $(\si_g)_{g \in G}$ of $N$ and
unitaries $(u_{g,h})_{g,h \in G}$ in $N$ satisfying
$$
\si_g \si_h  = (\Ad u_{g,h}) \si_{gh} \; , \quad u_{g,h}
\, u_{gh,k} = \si_g(u_{h,k}) \, u_{g,hk} \; , \quad
\si_e = \id \quad\text{and}\quad u_{e,e} = 1 \; ,
$$
for all $g,h,k \in G$.

A cocycle action $(\si_g)$ of $G$ on $N$ is said to be \emph{outer
  conjugate} to a cocycle action $(\rho_g)$ of $G$ on $M$ if there
  exists an isomorphism $\Delta : N \recht M$ such that $\Delta \si_g \Delta^{-1}= \rho_g \mod \Inn M$ for all $g \in G$.
\end{definition}

Note that a stronger notion of conjugacy exists, called \emph{cocycle
  conjugacy}, where it is imposed that $\Delta \si_g \Delta^{-1} = (\Ad w_g) \rho_g$,
  with unitaries $(w_g)$ making the $2$-cocycles for $\si$ and $\rho$
  cohomologous. In the case of an outer conjugacy between cocycle
  actions on a factor, the associated
  $2$-cocycles are only made cohomologous up to a scalar-valued $2$-cocycle.

Cocycle actions on a II$_1$ factor can be obtained by reducing an action by a projection. Let $(\si_g)$ be an action of $G$ on the II$_1$ factor $N$.
Whenever $p$ is a non-zero projection in $N$, choose partial isometries $w_g \in N$ such that $w_g w_g^* = p$ and $w_g^* w_g = \si_g(p)$. This is
possible because $(\si_g)$ preserves the trace and hence, $p$ and $\si_g(p)$ are equivalent projections since they have the same trace. Define
\begin{equation} \label{eq.reduction}
\si_g^p \in \Aut(pNp) : \si^p_g(x) = w_g \si_g(x) w_g^*
\quad\text{and}\quad u_{g,h} \in \cU(pNp) : u_{g,h} = w_g \si_g(w_h) w_{gh}^*
\; .
\end{equation}
It is easily checked that $(\si_g^p)$ is a cocycle action of $G$ on
the II$_1$ factor $pNp$ and that its outer conjugacy class does not
depend on the choice of $w_g$.

\begin{definition}
Let $(\si_g)$ be an action of the countable group $G$ on the II$_1$
factor $N$. Whenever $t > 0$, the cocycle action $(\si^t_g)$ of $G$
on $N^t$ is defined by reducing the action $(\id \ot \si_g)$ of $G$
on $\M_n(\C) \ot N$ by a projection $p$ with $(\Tr \ot \tau)(p) = t$,
as in \eqref{eq.reduction}.

The \emph{fundamental group} $\fun(\si)$ of the action $\si$ is
defined as the group of $t > 0$ such that $(\si^t_g)$ and $(\si_g)$
are outer conjugate.
\end{definition}
It is clear that $\fun(\si)$ is an outer conjugacy invariant for
$(\si_g)$. The following theorem computes the fundamental group for
Connes-St{\o}rmer Bernoulli actions of
$w$-rigid groups.

\begin{theorem}[Popa, \cite{P4}] \label{thm.fundamental-hyperfinite}
Let $(\cN,\vphi)$ be an almost periodic von Neumann algebra and
suppose that $N:=\cN^\vphi$ is a II$_1$ factor. Suppose that the countable group
$G$ admits an infinite normal subgroup $H$ with the relative property
(T) and that $(\si_g)$ is a malleable action of $G$ on $(\cN,\vphi)$
whose restriction to $H$ is weakly mixing.

If we still denote by $(\si_g)$ the restricted action of $G$ on the II$_1$ factor
$N$, then $\fun(\si) = \Sp(\cN,\vphi)$.
\end{theorem}
\begin{proof}
If $s \in \Sp(\cN,\vphi)$, we take a non-zero partial isometry $v \in
\cN$ which is an $s$-eigenvector for $\vphi$. Denote $p = v^* v$ and
$q = vv^*$. Then, $\Ad v$ outer conjugates $(\si^p_g)$ and
$(\si^q_g)$. Since $s = \frac{\vphi(q)}{\vphi(p)}$, it follows that $s
\in \fun(\si)$.

Conversely, suppose that $s \in \fun(\si)$. We have to prove that $s
\in \Sp(\cN,\vphi)$. We may clearly assume that $0 < s < 1$ and take a
projection $p \in N$ and elements $w_g \in N$ such that $\vphi(p) =
s$, $w_g w_g^* = p$ and $w_g^* w_g = \si_g(p)$ for all $g \in G$ and
such that $\rho_g(x) = w_g \si_g(x) w_g^*$ defines a genuine action of
$G$ on $pNp$ that is conjugate to $(\si_g)$. We only retain that
$(\rho_g)$ is a genuine action and that its 
restriction $\rho|_H$ is weakly mixing.

Let $(\al_t)$ be the one-parameter group on $\cN \ot \cN$ given by the
malleability of $(\si_g)$. As in the proof of Lemma \ref{lem.crucialstep}, the relative
property (T) yields $t_0=1/n$ and a non-zero partial isometry $a \in (\cN \ot
\cN)^{\vphi \ot \vphi}$ such that $aa^* \leq p \ot 1$, $a^*a \leq
\al_{t_0}(p \ot 1)$ and
$$(w_g \ot 1)(\si_g \ot \si_g)(a) = a \al_{t_0}(w_g \ot 1)
\quad\text{for all}\quad g \in H \; .$$
Weak mixing of $\si|_H$ on $\cN$ and of $\rho|_H$ on $pNp$ implies
that $aa^*=p \ot 1$ and $a^*a = \al_{t_0}(p \ot 1)$. Taking
$b:=a \al_{t_0}(a) \cdots \al_{(n-1)t_0}(a)$, we get a partial
isometry $b \in (\cN \ot \cN)^{\vphi \ot \vphi}$ satisfying $bb^* = p
\ot 1$, $b^*b = 1 \ot p$ and
$$(w_g \ot 1)(\si_g \ot \si_g)(b) = b (1 \ot w_g)
\quad\text{for all}\quad g \in H \; .$$ Continuing as in the proof of
Lemma \ref{lem.crucialstep}, Step~(3), we obtain the following data: a non-zero partial
isometry $v \in p \cN \ot \M_{1,n}(\C)$ which is a $\gamma$-eigenvector for $\vphi$ and satisfies $v^*v = 1$ as well as $w_g (\si_g \ot \id)(v) = v
(1 \ot \theta(g))$ for all $g \in H$, where $\theta : G \recht \cU(n)$ is a projective representation. The ergodicity of $\rho|_H$ yields $vv^* = p$
and hence, $\Ad v$ conjugates the actions $\rho|_H$ and $(\rho_g \ot \Ad \theta(g))_{g \in H}$. Since $1 \ot \M_n(\C)$ is an invariant subspace of
the latter, weak mixing of $\rho|_H$ imposes $n=1$. Since $vv^*=p$, $v^*v=1$ and $v$ is a $\gamma$-eigenvector, we conclude that $s = 1/\gamma
\in \Sp(\cN,\vphi)$.
\end{proof}

In Section \ref{sec.bernoulli}, Connes-St{\o}rmer Bernoulli actions were shown to be malleable and mixing. The following corollary is then clear.

\begin{corollary}
Let $G$ be a countable group
that admits an infinite normal subgroup with the relative property
(T). Let $\Tr_\Delta$ be the faithful normal state on $\B(H)$ given by
$\Tr_\Delta(a) = \Tr(\Delta a)$ and define
$(\cN,\vphi) = \bigotimes_{g \in G} (\B(H),\Tr_\Delta)$, with
Connes-St{\o}rmer Bernoulli action $G \actson (\cN,\vphi)$.
Write $\cR := \cN^\vphi$ and denote by
$(\si_g)$ the restricted action of $G$. Then, $\fun(\si)$ is the
subgroup of $\R^*_+$ generated by the ratios $\lambda/\mu$ between
$\lambda,\mu$ in the
point spectrum of $\Delta$.

In particular, $G$ admits a continuous family of non outer conjugate
actions on the hyperfinite II$_1$ factor $\cR$.
\end{corollary}

In Theorem \ref{thm.fundamental-hyperfinite} the following question was studied: when is the cocycle action $(\si^t_g)$ outer conjugate to $(\si_g)$?
Another natural question is: when is the cocycle action $(\si^t_g)$ outer conjugate to a genuine action. The following remark shows that $(\si^t_g)$
is always outer conjugate to a genuine action when $(\si_g)$ is a Connes-St{\o}rmer Bernoulli action on the centralizer of $\otimes_{g \in G}
(\B(H),\vphi_0)$ for $\vphi_0$ \emph{non-tracial}. On the other hand, for $\vphi_0$ the \emph{trace} on $\M_2(\C)$ and $t$ not an integer,
$(\si^t_g)$ is not outer conjugate to a genuine action, see Theorem \ref{thm.not-conjugate} below.

\begin{remark}
Let $(\cN,\vphi)$ be an almost periodic factor with $N:=\cN^\vphi$ a type II$_1$ factor and $\vphi$ non-tracial (note that this implies that $\cN$ is
a factor of type III$_\lambda$ with $0 < \lambda \leq 1$). Suppose that the group $G$ acts on $(\cN,\vphi)$ and denote by $(\si_g)$ the restriction
of this action to $N$. Then, for any $t
> 0$, $(\si^t_g)$ is outer conjugate to a genuine action.

For simplicity of notation, suppose $t \leq 1$ and let $p \in N$ be a
projection with $\vphi(p)=t$. We can write a series $t = \sum_n
\gamma_n$ with $\gamma_n \in \Sp(\cN,\vphi)$. Write $p = \sum_n p_n$
for some mutually orthogonal projections $p_n$ in $N$ with $\vphi(p_n)
= \gamma_n$. Take partial isometries
$v_n \in \cN$ such that $v_n$ is a $\gamma_n$-eigenvector for $\vphi$
and $v_n^* v_n = 1$, $v_n v_n^* = p_n$. Define for $g \in G$, the
element $w_g \in N$ as
$$w_g := \sum_n v_n \si_g(v_n^*) \; .$$
It is easy to check that $w_g w_g^* = p$, $w_g^* w_g = \si_g(p)$ for
all $g \in G$ and $w_g \si_g(w_h) = w_{gh}$ for all $g,h \in
G$. Writing $\si^p_g(x) = w_g \si_g(x) w_g^*$ for $x \in pNp$, it
follows that $(\si^p_g)$ is a genuine action of $G$ on $pNp$ and a way
to write $(\si^t_g)$.
\end{remark}

\begin{theorem}[Popa, \cite{P4}] \label{thm.not-conjugate}
Suppose that the countable group $G$ admits an infinite normal subgroup $H$ with the relative property (T). Denote by $(\si_g)$ the Bernoulli action
of $G$ on $\cR = \otimes_{g \in G} (M_2(\C),\tau)$. For $t > 0$, the cocycle action $(\si^t_g)$ is outer conjugate to a genuine action if and only if
$t \in \N_0$.
\end{theorem}

Observe moreover that it follows from Theorem \ref{thm.fundamental-hyperfinite} that, for different values of $t > 0$, the cocycle actions
$(\si^t_g)$ are mutually non outer conjugate.

\begin{proof}
Given $(\si^t_g)$ outer conjugate to a genuine action $(\rho_g)$, we can start off in the same way as in the proof of
\ref{thm.fundamental-hyperfinite}, but we do not know anymore that $\rho|_H$ is weakly mixing (or even, that $\rho$ is ergodic). So, in order to make
the passage from \lq an intertwiner for $\al_{t_0}$\rq\ to \lq an intertwiner for $\al_{1}$\rq, we need the extra data of \emph{strong malleability},
as in the proof of Lemma \ref{lem.crucialstep}. But, the Connes-St{\o}rmer Bernoulli action $(\si_g)$ is not strongly malleable in the sense of
Definition \ref{def.malleable} in an obvious way. So, we need a more flexible notion, essentially replacing tensor products by graded tensor
products, see Remark \ref{rem.graded-malleable} below.

Let $t > 0$ and suppose that $(\si^t_g)$ is outer conjugate to a
genuine action. So, we can take $k \in \N$, a projection $p \in
\cR \ot \M_k(\C)$ with $(\tau \ot \Tr)(p) = t$ and partial isometries
$w_g \in \cR \ot \M_k(\C)$ such that $w_g w_g^* = p$, $w_g^*w_g = (\si_g
\ot \id)(p)$ and such that $\rho_g(x) = w_g (\si_g \ot \id)(x)
w_g^*$ defines an action of $G$ on $\cR^t := p(\M_k(\C) \ot
\cR)p$. Let $q \leq p$ be any non-zero projection in $\cR^t$ invariant
under $\rho|_H$. We shall prove that $q$ dominates a non-zero
projection $q_0$, invariant under $\rho|_H$ and with $(\tau \ot \Tr)(q) \in \N$. This of course
proves that $(\tau \ot \Tr)(p) \in \N$.

Combining Remark \ref{rem.graded-malleable} and the proof of
Lemma \ref{lem.crucialstep}, we find a non-zero partial isometry $v \in \cR \ot
\M_{k,n}(\C)$ and a projective representation $\theta : G \recht
U(n)$ such that $v^* v = 1$, $v v^* \leq q$ and such that $w_g (\si_g
\ot \id)(v) = v (1 \ot \theta(g))$ for all $g \in H$. Putting $q_0 =
vv^*$, we are done.
\end{proof}

\begin{remark}\label{rem.graded-malleable}
The Connes-St{\o}rmer Bernoulli action $(\si_g)$ of the group $G$ on $N:=\otimes_{g \in G} \M_2(\C)$ satisfies the following form of strong
malleability: the II$_1$ factor $N$ is $\Z/2\Z$-graded, the action $(\si_g)$ commutes with the grading and the graded tensor square $N \oth N$ is
equipped with a one-parameter group of automorphisms $(\al_t)$ and a period $2$ automorphism $\beta$, all commuting with the grading and satisfying
$$\al_1(x \oth 1) = 1 \oth x \; , \quad \be (x \oth 1) = x \oth 1 \quad\text{and}\quad \be \al_t \be =
\al_{-t} \quad\text{for all $x \in N, t \in \R$.}$$ To check that the Bernoulli action indeed admits such a graded strong malleability, it suffices
to construct the grading and $(\al_t)$, $\be$ on the level of $\M_2(\C)$ and take the infinite product.

More generally however, for any real Hilbert space $H_\R$, one
considers the complexified Clifford $^*$-algebra $\Cliff(H_\R)$, generated
by self-adjoint elements $s(\xi)$, $\xi \in H_\R$ with relations
$$s(\xi)^2 = \|\xi\|^2  \quad\text{for
  all} \quad \xi \in H_\R \quad\text{and}\quad \xi \mapsto s(\xi) \quad\text{$\R$-linear.}$$
The $^*$-algebra $\Cliff(H_\R)$ admits an obvious $\Z/2\Z$-grading
  such that the elements $s(\xi)$ have odd degree. Also,
  $\Cliff(H_\R)$ has a natural tracial state yielding the hyperfinite
  II$_1$ factor after completion if $H_\R$ is of infinite
  dimension. Clearly, any orthogonal representation on $H_\R$ extends
  to an action on $\Cliff(H_\R)$ preserving the grading. Finally, we
  have a canonical isomorphism $\Cliff(H_\R \oplus K_\R) \cong
  \Cliff(H_\R) \oth \Cliff(K_\R)$.

If one notes that $\Cliff(\R^2) \cong \M_2(\C)$, one defines $\al_t$
and $\be$ on $\Cliff(\R^2 \oplus \R^2)$ by the formulas
$$\al_t \; (s  \begin{pmatrix} \xi \\ \eta \end{pmatrix})
= s ( \begin{pmatrix} \cos \frac{\pi t}{2} & - \sin
  \frac{\pi t}{2} \\ \sin \frac{\pi t}{2} & \cos \frac{\pi t}{2} \end{pmatrix}
\begin{pmatrix} \xi \\ \eta \end{pmatrix} )
\quad\text{and}\quad \be \; (s \begin{pmatrix} \xi \\ \eta
\end{pmatrix}) = s \begin{pmatrix} \xi \\ -\eta
\end{pmatrix} \; .$$
The above procedure shows that also the so-called \emph{Bogolyubov actions}
are strongly malleable in a graded way.
\end{remark}

\renewcommand{\sectionname}{Appendix}
\renewcommand{\thesection}{\Alph{section}}
\setcounter{section}{0}

\section{The basic construction and Hilbert modules} \label{sec.basic}

Let $(\cN,\vphi)$ be a von Neumann algebra with almost periodic
faithful normal state $\vphi$ and let $B \subset \cN^\vphi$ be a von
Neumann subalgebra of the centralizer algebra. A particularly
interesting case, is the one where $\vphi$ is a trace and where we
consider an inclusion $B \subset (N,\tau)$. We briefly explain the
so-called \emph{basic construction} von Neumann algebra $\la
\cN,e_B \ra$, introduced in \cite{Skau,Chr} and used extensively by
Jones \cite{Jon3} in his seminal work on subfactors.
We refer to \cite{C7,GDJ,Jon3} for further reading and
briefly explain what is needed in this talk.

The basic construction $\la \cN,e_B \ra$ is
defined as the von Neumann subalgebra of $\B(L^2(\cN))$ generated by
$\cN$ and the orthogonal projection $e_B$ of $L^2(\cN)$ onto $L^2(B) \subset
L^2(\cN)$. It can be checked that $\la \cN,e_B \ra$ consists of those
operators $T \in  \B(L^2(\cN))$ that commute with the right module
action of $B$: $T(\xi b) = T(\xi) b$ for all $\xi \in L^2(\cN)$ and
$b \in B$.

The basic construction $\la \cN,e_B \ra$ comes equipped with a canonical normal semi-finite faithful weight $\vphih$ satisfying
$$\vphih(x e_B y) = \vphi(xy)  \quad\text{for all $x,y \in
  \cN$.}$$
If $\vphi$ is a tracial state, $\vphih$ is a semi-finite trace.

Let $(B,\tau)$ be a finite von Neumann algebra with faithful tracial
state $\tau$. Whenever $K$ is a right $B$-module, the commutant $B'$
of $B$ on $K$ is a semi-finite von Neumann algebra that admits a
canonical semi-finite trace $\tau'$, characterized by the
formula
$$\tau'(T T^*) = \tau(T^*T) \quad\text{whenever}\quad T : L^2(B)
\recht K \quad\text{is bounded and right $B$-linear.}$$
Observe that for every bounded right $B$-linear map $T : L^2(B) \recht
K$, the element $T T^*$ belongs to $B'$ and $T^*T$ belongs to $B$, acting on the
left on $L^2(B)$.

When $B$ is a factor, one defines $\dim_B(K) := \tau'(1)$ and calls
$\dim_B(K)$ the \emph{coupling constant}. It is a complete invariant
for countably generated $B$-modules, which means the following: if $\dim_B(K) =
+\infty$, $K$ is isomorphic to $\ell^2(\N) \ot L^2(B)$ as
a right $B$-module and if $\dim_B(K) = t$ and $p \in \M_n(\C) \ot B$
is a projection with $(\Tr \ot \tau)(p) = t$, then $K$ is isomorphic
with $p L^2(B)^{\oplus n}$ as a right $B$-module.

When $(B,\tau)$ is an arbitrary finite von Neumann algebra with faithful tracial
state $\tau$, the situation is slightly more complicated. If $E_\cZ$
denotes the center valued trace, i.e.\ the unique $\tau$-preserving
conditional expectation $E_\cZ : B \recht \cZ(B)$ of $B$ onto the
center of $B$, we know that $E_\cZ(xy) = E_\cZ(yx)$ for all $x,y \in
B$ and that $p \inside q$ if and only if $E_\cZ(p) \leq E_\cZ(q)$ whenever $p$
and $q$ are projections in $B$. Moreover, whenever the Hilbert space $K$ is a right
$B$-module and $\tau$ a faithful tracial state on $B$, we denote by
$B'$ the commutant of $B$ on $K$ as above and construct a normal,
semi-finite positive linear map $$E'_\cZ : (B')^+ \recht \{\;
\text{positive elements affiliated with $\cZ(B)$} \; \}$$ satisfying
$E'_\cZ(x^*x) = E'_\cZ(xx^*)$ for all $x$ and such that
$$E'_\cZ(T T^*) = E_\cZ(T^*T) \quad\text{whenever}\quad T : L^2(B)
\recht K \quad\text{is bounded and right $B$-linear.}$$
The positive affiliated element $E'_\cZ(1)$ of $\cZ(B)$ provides a
complete invariant for countably generated right $B$-modules. First note that the
$B$-module $K$ is finitely generated, i.e.\ of the form $p
L^2(B)^{\oplus n}$ for some projection $p \in \M_n(\C) \ot B$, if and
only if $E'_\cZ(1)$ is bounded. In that case $E'_\cZ(1) = (\Tr \ot
E_\cZ)(p)$.

Note that $\tau' = \tau \circ E'_\cZ$. So,
if $\tau'(1) < \infty$, it follows that $E'_\cZ(1)$ is not
necessarily bounded, but $\tau$-integrable. This implies that
$E'_\cZ(1)z$ is bounded for projections $z \in \cZ(B)$ with trace
arbitrary close to $1$. So, we have shown the following lemma.

\begin{lemma} \label{lem.cut-down-central}
Let $K$ be a right $B$-module and $\tau$ a normal faithful tracial
state on $B$. Denote by $\tau'$ the canonical semi-finite trace on the
commutant $B'$ of $B$ on $K$. If $\tau'(1) < \infty$, there exists for
any $\eps > 0$ a central projection $z \in \cZ(B)$ with $\tau(z) \geq
1 - \eps$ and such that the $B$-module $Kz$ is finitely generated,
i.e.\ of the form $p
L^2(B)^{\oplus n}$ for some projection $p \in \M_n(\C) \ot B$.
\end{lemma}

Returning to the basic construction for the inclusion $B \subset \cN$, with $B \subset \cN^\vphi$, we observe that the restriction of $\vphi$ defines
a tracial state on $B$ and that $\la \cN,e_B \ra$ is the commutant of $B$ on $L^2(\cN)$. Using the previous paragraph, $\la \cN,e_B \ra$ comes
equipped with a canonical semi-finite trace $\vphi'$. If $\vphi$ is tracial on $\cN$, it is easily checked that $\vphih = \vphi'$. If $\vphi$ is no
longer a trace, but an almost periodic state, we denote by $p_\gamma$ the orthogonal projection of $L^2(\cN)$ on the $\gamma$-eigenvectors for
$\vphi$. Note that $p_\gamma$ belongs to $\la \cN,e_B \ra$ because $B \subset \cN^\vphi$. It is easy to check that
$$\vphih(x) = \sum_{\gamma \in \Sp(\cN,\vphi)} \vphih(p_\gamma x
p_\gamma) \quad\text{and}\quad \vphi'(x) = \sum_{\gamma \in
  \Sp(\cN,\vphi)} \gamma^{-1} \vphih(p_\gamma x
p_\gamma)$$ for all $x \in \la \cN,e_B \ra^+$. In particular, $\vphih$ is tracial and a multiple of $\vphi'$ on $p_\gamma \la \cN,e_B \ra p_\gamma$,
for all $\gamma \in \Sp(\cN,\vphi)$.

\section{Relative property (T) and II$_1$ factors} \label{sec.relativeT}

A countable group $G$ has Kazhdan's \emph{property (T)} if every
unitary representation of $G$ that admits  a
sequence of almost invariant unit vectors, admits a non-zero $G$-invariant vector.
More generally, a pair $(G,H)$ consisting of a countable group $G$ with subgroup $H$ is said to have the \emph{relative property (T)} of Kazhdan-Margulis
\cite{dHV,DK,Ka,Ma}, if every unitary representation of $G$ that
admits a sequence of almost invariant unit vectors, admits a non-zero
$H$-invariant vector. The main
example is the pair $(\SL(2,\Z) \ltimes \Z^2, \Z^2)$.

A countable group $G$ is said to be \emph{amenable} if the regular
representation on $\ell^2(G)$ admits a sequence of almost invariant
unit vectors. Hence, an amenable property (T) group is finite and an
amenable group does not have an infinite subgroup with the relative
property (T).

Below, we need the following alternative characterization of relative
property (T) due to Jolissaint (see Theorem 1.2(a3) in \cite{J}). The pair $(G,H)$ has
the relative property (T) if and only if every unitary representation
of $G$ admitting a sequence of almost invariant unit vectors, admits a
non-zero $H$-invariant finite dimensional subspace.

The notion of property (T) has been defined for II$_1$ factors by Connes and Jones \cite{CJ}. Unitary representations of groups are replaced by
\emph{bimodules} (\emph{Connes' correspondences}, see \cite{C5,P8}). Popa \cite{P5} defined the relative property (T) for an inclusion of finite von
Neumann algebras $Q \subset P$ and we explain it in this appendix.

A $P$-$P$ bimodule is a Hilbert space $H$ with a left and a right (normal, unital)
action of $P$. We write $x \xi$, resp.\ $\xi x$ for the left, resp.\ right action of $P$ on $H$.
\begin{terminology}
Let $(P,\tau)$ be a von Neumann algebra with a faithful normal tracial
state $\tau$. If $K$ is a $P$-$P$-bimodule and $(\xi_n)$ a sequence of unit vectors in $K$, we say that
\begin{itemize}
\item $(\xi_n)$ is \emph{almost central} if $\|x \xi_n - \xi_n x \| \recht 0$ for all $x \in P$;
\item $(\xi_n)$ is \emph{almost tracial} if $\|\la \xi_n, \cdot
  \xi_n\ra - \tau\| \recht 0$ and $\|\la \xi_n, \xi_n \cdot \ra -
  \tau\| \recht 0$.
\end{itemize}
\end{terminology}

A vector $\xi$ is said to be \emph{$Q$-central} for some von Neumann subalgebra $Q \subset P$ if $x \xi = \xi x$ for all $x \in Q$.

\begin{definition}[Popa, \cite{P5}] \label{def.relativeT}
Let $(P,\tau)$ be a von Neumann algebra with a faithful normal tracial
state $\tau$. The inclusion $Q \subset P$ is
said to have the \emph{relative property (T)} if any $P$-$P$ bimodule
that admits a sequence of almost central almost tracial unit vectors,
admits a sequence of almost tracial
$Q$-central unit vectors.
\end{definition}

\begin{remark}
One might wonder why almost traciality is assumed in the definition of
relative property (T). In applications (as the ones Popa's work), it is
crucial that an inclusion $Q \subset P$ with the relative property (T) remains relative (T) when cutting down with a projection of $Q$ (see
Proposition \ref{prop.Tcutdown}). Now look at the following example:
we take a II$_1$ factor $P$, two von Neumann subalgebras $Q_1,Q_2
\subset P$ and we consider the
inclusion of $Q_1 \oplus Q_2 \subset \M_2(\C) \ot P$. If one would define naively the relative property (T) by imposing that any $P$-$P$ bimodule
admitting almost central vectors, admits a non-zero $Q$-central vector, then the inclusion $Q_1 \oplus Q_2 \subset \M_2(\C) \ot P$ would have the
relative property (T) if \emph{one of the inclusions} $Q_1 \subset P$, $Q_2 \subset P$ has the relative property (T). And hence, Proposition
\ref{prop.Tcutdown} would not hold.
\end{remark}

\begin{remark} \label{rem.amenable}
A finite von Neumann algebra $(P,\tau)$ with faithful normal tracial
state $\tau$ is said to be \emph{injective} (or \emph{amenable}) if
the coarse Hilbert $P$-$P$-bimodule $L^2(P) \ot L^2(P)$ defined by $a
\cdot \xi \cdot b = (a \ot 1) \xi (1 \ot b)$ contains a sequence of
almost central almost tracial vectors. It is then clear that an
injective $(P,\tau)$ does not contain a diffuse subalgebra $Q \subset
P$ with the relative property (T). More generally, if $Q \subset P$ is
diffuse with the relative property (T), there is no non-zero normal
homomorphism from $P$ to an injective finite von Neumann algebra.
\end{remark}

A lot can be said about relative property (T) in the setting of von
Neumann algebras, see the papers of Peterson and Popa \cite{PP,P5}. In this talk, only two easy
results are shown, which suffices for the applications in the rest of
the talk.

\begin{proposition} \label{prop.from-group-to-vnalg}
Let $G$ be a countable group with subgroup $H$. Then, $(G,H)$ has the
relative property (T) if and only if the inclusion $\cL(H) \subset \cL(G)$ has the relative property (T) in the sense of Definition \ref{def.relativeT}.
\end{proposition}
\begin{proof}
First suppose that $(G,H)$ has the relative property (T). Let $K$ be an $\cL(G)$-$\cL(G)$-bimodule with an almost central almost $\tau$-tracial
sequence of unit vectors $(\xi_n)$, for some faithful normal tracial state $\tau$ on $\cL(G)$. Define the representation $\pi(g)\xi = u_g \xi u_g^*$
of $G$ on $K$. Choose $\eps > 0$. Using the stronger version of relative property (T), we can take a $\pi(H)$-invariant unit vector $\xi$ and $n \in
\N$ such that
$$\|\xi - \xi_n\| < \frac{\eps}{3} \; , \quad \|\la \xi_n, \cdot \xi_n\ra - \tau\| < \frac{\eps}{3} \; ,\quad \|\la \xi_n, \xi_n \cdot \ra - \tau\| <
\frac{\eps}{3} \; .$$ Since a $\pi(H)$-invariant vector is $\cL(H)$-central, we have found an $\cL(H)$-central unit vector $\xi$ satisfying
$$\|\la \xi, \cdot \xi\ra - \tau\| < \eps \; ,\quad \|\la \xi, \xi \cdot \ra - \tau\| <
\eps \; .$$ It follows that $K$ admits a sequence of almost tracial $\cL(H)$-central vectors.

Conversely, suppose that the inclusion $\cL(H) \subset \cL(G)$ has the
relative property (T) in the sense of Definition
\ref{def.relativeT}. Let $\pi : G \recht \cU(K_0)$ be a unitary
representation of $G$ that admits a sequence $(\xi_n)$ of almost
invariant unit vectors. As stated above, it is sufficient to prove that $K_0$
admits a non-zero finite-dimensional $\pi(H)$-invariant subspace. Define $K=\ell^2(G)
\ot K_0$, which we turn into an $\cL(G)$-$\cL(G)$-bimodule by the
formulas
$$u_g \cdot (\delta_h \ot \xi) = \delta_{gh} \ot \pi(g)\xi
\quad\text{and}\quad (\delta_h \ot \xi)\cdot u_g = \delta_{hg} \ot \xi
$$
for all $g,h \in G$, $\xi \in K_0$. It is clear that $(\delta_e \ot \xi_n)$ is a sequence of almost central almost tracial unit vectors. So, $K$
admits a non-zero $\cL(H)$-central vector $\mu$. Considering $\mu$ as an element in $\ell^2(G,K_0)$, we get that $\mu(hgh^{-1}) = \pi(h) \mu(g)$ for
all $h \in H,g \in G$. Take $g \in G$ such that $\mu(g) \neq 0$. Since $\mu \in \ell^2(G,K_0)$, we conclude that $\{hgh^{-1} \mid h \in H\}$ is
finite. But then, the linear span of $\{\mu(hgh^{-1}) \mid h \in H\}$ is a finite-dimensional $\pi(H)$-invariant subspace of $K_0$.
\end{proof}

\begin{proposition} \label{prop.Tcutdown}
Let $P$ be a II$_1$ factor and $Q \subset P$ an inclusion having the relative property (T). If $p \in Q$ is a non-zero projection, $pQp \subset pPp$
has the relative property (T).
\end{proposition}
\begin{proof}
Write $Q_1 = pQp$ and $P_1 = pPp$. Since $P$ is a II$_1$ factor, we can take partial isometries $v_1,\ldots,v_k \in P$ satisfying $v_1 = p$, $v_i^*
v_i \leq p$ and $\sum_{i=1}^k v_i v_i^* = 1$. Let $K_1$ be a $P_1$-$P_1$-bimodule admitting the almost central almost tracial sequence of unit
vectors $(\xi_n)$. Define $K$ as the induced $P$-$P$-bimodule: put a scalar product on $Pp \, K_1 \, pP$ by the formula
$$\la x \xi y^*, a \mu b^* \ra = \la \xi, (x^*a) \mu (b^* y) \ra \quad\text{for all}\;\; x,y,a,b \in Pp, \; \xi,\mu \in K_1 \; .$$
Up to normalization, the sequence $\sum_{i=1}^k v_i \xi_n v_i^*$ is almost central almost tracial in the $P$-$P$-bimodule $K$. Hence, $K$ admits an
almost tracial sequence $(\mu_n)$ of $Q$-central vectors. Up to normalization, $(p\mu_n) = (\mu_n p)$ defines an almost tracial sequence of
$pQp$-central vectors in $K_1$.
\end{proof}

The above proposition remains valid when $(P,\tau)$ is just von
Neumann algebra with faithul tracial state $\tau$, but the proof becomes slightly more
involved.

\section{Intertwining subalgebras using bimodules} \label{sec.intertwining}

The fundamental problem in the whole of this talk is to decide when two von Neumann subalgebras $P,B \subset M$ can be conjugated one into the
other: $u P u^* \subset B$ for some $u \in \cU(M)$. The usage of the
basic construction in this respect goes back to Christensen
\cite{Chr}, who used it to study conjugacy of uniformly close
subalgebras. A major innovation came with the work of Popa
\cite{P1,P5}, who managed to prove conjugacy results for arbitrary
subalgebras, still using the basic construction.

Roughly, Proposition \ref{prop.translate} below says the following.
Let $P,B \subset M$ be von
Neumann subalgebras of a finite von Neumann algebra $(M,\tau)$. Then, the following are equivalent.
\begin{itemize}
\item A corner of $P$ can be conjugated into a corner of $B$.
\item $L^2(M)$ contains a non-zero $P$-$B$-subbimodule which is finitely
  generated as a $B$-module.
\item The basic construction $\la M,e_B \ra$ contains a
  positive element $a$, commuting with $P$ and satisfying $0 <
  \tauh(a) < +\infty$, where $\tauh$ is the canonical semi-finite
  trace on $\la M,e_B \ra$.
\end{itemize}
The relation between the second and the third condition is clear: the
orthogonal projection $p_K$ onto a $P$-$B$-subbimodule $K$ of $L^2(M)$ belongs
to $\la M,e_B \ra \cap P'$ and $\tauh(p_K) < \infty$ is essentially
equivalent to $K$ being a finitely generated $B$-module.

We reproduce from \cite{P1,P5} two results needed in this talk.

\begin{proposition}[Popa, \cite{P1,P5}] \label{prop.translate}
Let $(\cM,\vphi)$ be a von Neumann algebra with an almost periodic faithful normal state $\vphi$. Let $P,B \subset \cM^\vphi$ be von Neumann
subalgebras. Then, the following statements are equivalent.
\begin{enumerate}
\item\label{stat.one} There exists $n \geq 1$, $\gamma > 0$, $v \in \M_{1,n}(\C) \ot \cM$, a projection $p \in \M_n(\C) \ot B$ and a homomorphism $\theta : P \recht p(\M_n(\C) \ot B)p$
such that $v$ is a non-zero partial isometry which is a $\gamma$-eigenvector for $\vphi$, $v^*v \leq p$ and
$$x v = v \theta(x) \quad\text{for all}\quad x \in P \; .$$
\item\label{stat.two} There exists a non-zero element $w \in \cM$ such that $P w \subset \sum_{k=1}^n w_k B$ for some finite family $w_k$ in $\cM$.
\item\label{stat.three} There exists a non-zero element $a \in \la \cM,e_B \ra^+ \cap P'$ with $\vphih(a) < \infty$. Here $\la \cM,e_B \ra$ denotes the basic
construction for the inclusion $B \subset \cM$, with its canonical almost periodic semi-finite weight $\vphih$.
\item\label{stat.four} There is no sequence of unitaries $(u_n)$ in $P$ such that $\|E_B(a u_n b)\|_2 \recht 0$ for all $a,b \in \cM$.
\end{enumerate}
\end{proposition}

Of course, if one wants to deal as well with the non-separable case, one should take a net instead of a sequence in statement \eqref{stat.four}.

\begin{proof}
\eqref{stat.one} $\Rightarrow$ \eqref{stat.two}. Taking a non-zero component of $v$, this is trivial.

\eqref{stat.two} $\Rightarrow$ \eqref{stat.three}. Since $P$ and $B$ are in the centralizer algebra $\cM^\vphi$ and $\vphi$ is almost periodic, we
can assume that $w,w_1,\ldots,w_n$ are all $\gamma$-eigenvectors for $\vphi$. Note that, whenever $w \in \cM$ is a $\gamma$-eigenvector, the
projection of $L^2(\cM)$ onto the closure of $wB$ yields a projection $f \in \la \cM,e_B \ra$ and $f$ is the range projection of $w e_B w^*$. It
follows that $\vphih(f) \leq \gamma$. In the same way, the projection onto the closure of $\sum_{k=1}^n w_k B$  has finite $\vphih$-weight. Hence,
the projection $f$ onto the closure of $P w B$ in $L^2(\cM)$ satisfies the requirements in \eqref{stat.three}.

\eqref{stat.three} $\Rightarrow$ \eqref{stat.one}. If $p_\gamma$ denotes the orthogonal projection of $L^2(\cM)$ onto the $\gamma$-spectral subspace
of $\vphi$, we know that $\vphih(a) = \sum_\gamma \vphih(p_\gamma a p_\gamma)$ and we can replace $a$ by $p_\gamma a p_\gamma \neq 0$. Taking a
spectral projection of the form $\chi_{[\delta,+\infty[}(a)$, we obtain an orthogonal projection $f \in \la \cM,e_B \ra^+ \cap P'$ with $\vphih(f) <
\infty$ and the range of $f$ contained in the $\gamma$-spectral subspace of $\vphi$. Hence, the range of $f$ is a non-zero $P$-$B$-sub-bimodule of
$L^2(\cM)_\gamma$ with finite trace over $B$. Cutting down by a central projection of $B$ (see Lemma \ref{lem.cut-down-central}), we get a
$P$-$B$-sub-bimodule $H \subset L^2(\cM)_\gamma$ which is finitely generated over $B$. Hence, we can take $n \geq 1$, a projection $p \in \M_n(\C)
\ot B$ and a $B$-module isomorphism
$$\psi : p L^2(B)^{\oplus n} \recht H \; .$$
Since $H$ is a $P$-module, we get a homomorphism $\theta : P \recht p (\M_n(\C) \ot B)p$ satisfying $x \psi(\xi) = \psi(\theta(x)\xi)$ for all $x \in
P$ and $\xi \in H$. Define $e_i \in L^2(B)^{\oplus n}$ as $e_i = (0,\ldots,1,\ldots,0)$ and $\xi \in \M_{1,n}(\C) \ot H$ as $\xi_i = \psi(pe_i)$. The
polar decomposition of the vector $\xi$ yields an isometry $v \in \M_{1,n}(\C) \ot \cM$ belonging to the $\gamma$-spectral subspace for $\vphi$. A
direct computation shows that $x v = v \theta(x)$ for all $x \in P$, as well as $v^*v \leq p$.

\eqref{stat.one} $\Rightarrow$ \eqref{stat.four}. Suppose that we have all the data of \eqref{stat.one}. If $(u_n)$ is a sequence of unitaries in $P$
such that $\|E_B(a u_n b)\|_2 \recht 0$ for all $a,b \in \cM$, it follows that $\|(\id \ot E_B)(v^* u_n v) \|_2 \recht 0$ when $n \recht \infty$.
But, $\|(\id \ot E_B)(v^* u_n v) \|_2 = \|(\id \ot E_B)(v^*v) \theta(u_n)\|_2 = \|(\id \ot E_B)(v^*v)\|_2$. We conclude that $v=0$, a contradiction.

\eqref{stat.four} $\Rightarrow$ \eqref{stat.three}. By \eqref{stat.four}, we can take $\eps > 0$ and $K \subset \cM$ finite such that for all
unitaries $u \in P$, $\max_{a,b \in K} \|E_B(a u b)\|_2 \geq
\eps$. Define the element $c = \sum_{b \in K} b e_B b^*$ in $\la
\cM,e_B \ra^+$. Note that $\vphih(c) < \infty$. Let $d \in \la
\cM,e_B \ra^+$ be the element of minimal $L^2$-norm (with respect to
$\vphih$) in the $L^2$-closed convex hull of $\{u c u^* \mid u \in
\cU(P) \}$. By uniqueness of the element of minimal $L^2$-norm, it
follows that $d \in \la \cM,e_B \ra^+ \cap P'$ and by construction
$\vphih(d) < \infty$. It remains to show that $d \neq 0$. But, for all
$u \in \cU(P)$, we have
$$\sum_{a \in K} \vphih(e_B a \; u c u^* \; a^* e_B) = \sum_{a,b \in K} \|E_B(a u
b)\|_2^2 \geq \eps^2 \; .$$
It follows that $\sum_{a \in K} \vphih(e_B a d a^* e_B) \geq \eps^2$
and $d \neq 0$.
\end{proof}

\begin{lemma} \label{lem.abelian}
Let $M$ be a finite von Neumann algebra and $B \subset M$ a maximal
abelian subalgebra.
\begin{itemize}
\item If $q \in M$ is an abelian projection, there exists $v \in M$
  satisfying $v^*v=q$ and $vMv^* \subset B$.
\item If $M$ is of finite type I and $P_0 \subset M$ an abelian von Neumann
  subalgebra, there exists a unitary $u \in M$ such that $uP_0 u^*
  \subset B$.
\end{itemize}
\end{lemma}
\begin{proof}
We do not provide a full proof of this classical lemma: see paragraph
6.4 in \cite{KR} for the necessary background. The following
indications shall allow the reader to fill in the proof.

For the first statement, it suffices to find a projection in $B$ which
is equivalent with $q$, i.e.\ $v \in M$ with $v^*v=q$ and $vv^* \in
B$. Since $B$ is maximal abelian, we have $vMv^* \subset B$.

For the second statement: since $M$ is of finite type I and $L^\infty(X)=B
\subset M$ is maximal abelian, the partial isometries in $M$ normalizing
$B$ induce an equivalence relation with finite orbits on $X$. Taking a
fundamental domain for this equivalence relation, we can easily
conclude. Of course, a proper proof can be given in operator algebraic
terms: if $M$ is of type I$_n$ and $B \subset M$ maximal abelian, we
can write $1$ as the sum of $n$ equivalent abelian projections
contained in $B$. Embedding $P_0 \subset P \subset M$ with $P$ maximal
abelian, we can do the same with $P$ and then, $P$ and $B$ are unitary conjugate.
\end{proof}

\begin{theorem}[Popa, \cite{P5}] \label{thm.intertwine-abelian}
Let $(M,\tau)$ be a finite von Neumann algebra and $P_0,B \subset M$
abelian subalgebras. Suppose that $B$ is maximal abelian and $P:= M
\cap P_0'$ abelian (hence, maximal abelian). The following statements
are equivalent.
\begin{enumerate}
\item\label{stat.a} There exists a non-zero $v \in M$ such that $P_0v
  \subset \sum_{k=1}^n v_k B$ for some finite set of elements $(v_k)$
  in $B$.
\item\label{stat.b} There exists a non-zero $a \in \la M,e_B \ra^+
  \cap P_0'$ satisfying $\tauh(a)< \infty$. Here $\la M,e_B \ra$
  denotes the basic construction for the inclusion $B \subset M$ and
  $\tauh$ is the canonical semi-finite trace on it.
\item\label{stat.c} There exists a non-zero partial isometry $v \in M$
  such that $v^*v \in P$, $p:=vv^* \in B$ and $vPv^* = Bp$.
\end{enumerate}
If moreover $M$ is a factor and $P$ and $B$ are Cartan subalgebras, a
fourth statement is equivalent:
\begin{enumerate}\setcounter{enumi}{3}
\item\label{stat.d} There exists a unitary $u \in M$ such that $uPu^*=B$.
\end{enumerate}
\end{theorem}
\begin{proof}
Given Proposition \ref{prop.translate}, it suffices to prove that
\eqref{stat.b} implies \eqref{stat.c} as well as \eqref{stat.d} under
the additional assumption that $M$ is factorial and $P$ and $D$ are Cartan.

Using Proposition \ref{prop.translate}, we take $n \geq 1$, a
projection $p \in \M_n(\C) \ot B$, a non-zero partial isometry $w \in
\M_{1,n}(\C) \ot M$ and a homomorphism $\theta : P_0 \recht p(\M_n(\C)
\ot B)p$ such that $x w= w \theta(x)$ for all $x \in P_0$. We can
 replace $p$ by an equivalent projection in $\M_n(\C) \ot B$ and
take $p = \operatorname{diag}(p_1,\ldots,p_n)$. Then,
$\diag(p_1B,\ldots,p_nB)$ is a maximal abelian subalgebra of the
finite type I algebra $p(\M_n(\C) \ot B)p$. Since $P_0$ is abelian,
Lemma \ref{lem.abelian} allows to suppose that $\theta(P_0) \subset
\diag(p_1B,\ldots,p_nB)$. Hence, we can cut down $\theta$ and $w$ by
one of the projections $(0,\ldots,p_i,\ldots,0)$ and suppose from the
beginning that $n=1$.

Write $q:=w^*w$, $e:=ww^* \in P$ and $A:=pMp \cap \theta(P_0)'$. Then, $q
\in A$ and $qAq =w^*( eMe \cap (Pe)')w = w^*Pw$, which is abelian. Since
$A$ is finite and $pB \subset A$ maximal abelian, Lemma
\ref{lem.abelian} gives $u \in A$ satisfying $uu^* = q$ and $u^*Au
\subset pB$. Writing $v=u^*w^*$, we have $vPv^* \subset B$ and
$v^*v=e$. Write $f:=vv^* \in B$. Hence, $eP \subset v^*Bv \subset
eMe$. Since $v^*Bv$ is abelian, it follows that $eP = v^*Bv$ and so,
$vPv^* = fB$.

Assume now that $M$ is a factor and that $P,B \subset M$ are Cartan
subalgebras. Whenever $u_1$ is a unitary in $M$ normalizing $P$ and
$u_2$ is a unitary in $M$ normalizing $B$, $u_2 v u_1$ moves as well a
corner of $P$ into a corner of $B$. A maximality argument yields \eqref{stat.d}.
\end{proof}

\section{Some results on (weakly) mixing actions} \label{sec.mixing}

An action of a countable group $G$ on $(\cA,\vphi)$ is said to be \emph{ergodic} if the scalars are the only $G$-invariant elements of $\cA$.
Equivalently, the multiples of $1$ are the only $G$-invariant vectors
in $L^2(\cA,\vphi)$. Stronger notions of ergodicity are the mixing and weak
mixing properties.

\begin{definition} \label{def.mixing}
An action of a countable group $G$ on $(\cA,\vphi)$ is said to be
\begin{itemize}
\item \emph{mixing} if for every $a,b \in \cA$, $\vphi(a \si_g(b)) \recht
  \vphi(a)\vphi(b)$
  when $g \recht \infty$;
\item \emph{weakly mixing} if for every $a_1,\ldots,a_n \in \cA$
and $\eps > 0$, there exists $g \in G$ such that $|\vphi(a_i\si_g(a_j)) - \vphi(a_i)\vphi(a_j)| < \eps$ for all $i,j =1,\ldots,n$.
\end{itemize}
\end{definition}

For the convenience of the reader, we prove the following classical equivalent characterizations for weakly mixing actions.

\begin{proposition} \label{prop.mixing}
Let a countable group $G$ act on the finite von Neumann algebra $(A,\tau)$ by automorphisms $(\si_g)$. Then, the following statements are equivalent.
\begin{enumerate}
\item\label{mix.1} The action $(\si_g)$ is weakly mixing.
\item\label{mix.2} For every $a_1,\ldots,a_k \in \cA$ with $\tau(a_i) = 0$, there exists a sequence $g_n \recht \infty$ in $G$ such that
$\si_{g_n}(a_i) \recht 0$ weakly for all $i=1,\ldots,k$.
\item\label{mix.3} $\C 1$ is the only finite-dimensional invariant subspace of $L^2(A)$.
\item\label{mix.4} $\C 1$ is the only finite-dimensional invariant subspace of $A$.
\item\label{mix.5} For every action $(\al_g)$ of $G$ on a finite von Neumann algebra $(M,\tau)$, $(A \ot M)^{\si \ot \al} = 1 \ot M^\al$.
\item\label{mix.6} The diagonal action of $G$ on $A \ot A$ is ergodic.
\end{enumerate}
\end{proposition}
\begin{proof}
The implications \eqref{mix.1} $\Rightarrow$ \eqref{mix.2} $\Rightarrow$ \eqref{mix.3} $\Rightarrow$ \eqref{mix.4}, as well as \eqref{mix.5}
$\Rightarrow$ \eqref{mix.6}, being obvious, we prove two implications below.

\eqref{mix.4} $\Rightarrow$ \eqref{mix.5}. Suppose that $X \in (A \ot M)^{\si \ot \al}$. Denote by $\eta$ the canonical embeddings $M \recht L^2(M)$
and $A \recht L^2(A)$. Define the Hilbert-Schmidt operator $T : \overline{L^2(M)} \recht L^2(A) : T\overline{\xi} = \eta\bigl((\id \ot
\om_{\xi,\eta(1)})(X)\bigr)$. Note that the image of $T$ is contained in $\eta(A)$ and that $TT^*$ commutes with the unitary representation $(\pi_g)$
on $L^2(A)$ given by $\pi_g \eta(a) = \eta(\si_g(a))$. Moreover, $TT^*$ is trace-class. Taking a spectral projection, we find a $G$-invariant
finite-dimensional subspace of $A$. By \eqref{mix.4}, the image of $T$ is included in $\C \eta(1)$, which means that $X \in 1 \ot M^\al$.

\eqref{mix.6} $\Rightarrow$ \eqref{mix.1}. Suppose that $(\si_g)$ is not weakly mixing. We find $\eps > 0$ and $a_1,\ldots,a_n$ with $\tau(a_i)=0$
and $\sum_{i,j=1}^n |\tau(a_j^* \si_g(a_i))|^2 \geq \eps$ for every $g \in G$. Define the vector $\xi = \sum_{i=1}^n a_i \ot a_i^*$ in $L^2(A \ot
A)$. Let $\xi_1$ be the element of minimal norm in the closed convex hull of $\{(\pi_g \ot \pi_g)\xi \mid g \in G\}$. Since for any $g \in G$,
$$\la \xi,(\pi_g \ot \pi_g)(\xi) \ra = \sum_{i,j=1}^n |\tau(a_j^* \si_g(a_i))|^2 \geq \eps$$
we conclude that $\xi_1 \neq 0$. Moreover, by the uniqueness of $\xi_1$, we get that $\xi_1$ is $(\pi_g \ot \pi_g)$-invariant. By construction
$\xi_1$ is orthogonal to $1$ and we have obtained a contradiction with \eqref{mix.6}.
\end{proof}

\begin{lemma} \label{lem.intertwiners}
Let $(\cM,\vphi)$ be an almost periodic von Neumann algebra and $P \subset B
\subset \cM^\vphi$ von Neumann subalgebras of the centralizer algebra $\cM^\vphi$. Suppose that there exists a sequence
of unitaries $(u_n)$ in $P$ such that
$$\|E_B(a u_n b)\|_2 \recht 0 \quad\text{whenever}\quad a,b \in \Ker
E_B \; ,$$
where $E_B : \cM \recht B$ is the $\vphi$-preserving conditional
expectation. If $x \in \cM$ is such that $P x \subset \sum_{k=1}^n x_k
B$ for a finite family of elements $x_k \in \cM$, then $x \in B$.

More generally, any $P$-$B$-sub-bimodule of $L^2(\cM)$ that is of finite trace as a $B$-module, is contained in $L^2(B)$.
\end{lemma}
\begin{proof}
Let $H_0 \subset L^2(\cM)$ be a $P$-$B$-sub-bimodule that is of finite trace as a $B$-module. We may assume that $H_0$ is contained in the
$\gamma$-eigenspace for $\vphi$ and that $H_0$ is orthogonal to $L^2(B)$. We have to prove that $H_0 = 0$. Suppose the contrary. As in Proposition
\ref{prop.translate}, we find $n \geq 1$, a non-zero partial isometry $v \in \M_{1,n}(\C) \ot \cM$, a projection $p \in \M_n(\C) \ot B$ and a
homomorphism $\theta : P \recht p(\M_n(\C) \ot B)p$ such that $a v = v \theta(a)$ for all $a \in P$, $v^*v \leq p$ and $(\id \ot E_B)(v) = 0$.

Using the $L^2$-norm with respect to the functional $\Tr \ot \vphi$ on $\M_n(\C) \ot \cM$, we have \linebreak $\|(\id \ot E_B)(v^* u_n v)\|_2 \recht
0$ when $n \recht \infty$. On the other hand, $v^* u_n v = \theta(u_n) v^*v$, implying that $\|(\id \ot E_B)(v^* u_n v)\|_2= \|(\id \ot
E_B)(v^*v)\|_2$. We conclude that $v^*v = 0$, a contradiction.
\end{proof}

\begin{theorem}[Popa, \cite{P1}] \label{thm.mixing-one}
Suppose that $G$ acts mixingly on an almost periodic $(\cN,\vphi)$ and write $\cM = \cN \rtimes G$. Let $p \in \M_n(\C) \ot \cL(G)$ a projection with
(non-normalized) trace $t$ and write $\cL(G)^t = p(\M_n(\C) \ot \cL(G))p$, $\cM^t =p(\M_n(\C) \ot \cM)p$. If $P \subset \cL(G)^t$ is a diffuse von
Neumann subalgebra, any $P$-$\cL(G)^t$-sub-bimodule of $L^2(\cM^t)$ that is of finite trace as an $\cL(G)^t$-module, is contained in $L^2(\cL(G)^t)$.
\end{theorem}

So, under the conditions of Theorem \ref{thm.mixing-one}, if $x \in \cM^t$ such that
$$P x \subset \sum_{k=1}^n x_k \cL(G)^t$$ for a finite family $x_k \in
\cM^t$, then $x \in \cL(G)^t$.

\begin{proof}
We claim that whenever $(x_n)$ is a bounded sequence in $\cL(G)$ that weakly tends to $0$,
$$\|E_{\cL(G)}(a x_n b)\|_2 \recht 0$$
when $n \recht \infty$, for all $a,b \in \Ker(E_{\cL(G)})$. Here
$E_{\cL(G)} : \cM \recht \cL(G)$ is the $\vphi$-preserving conditional
expectation. It suffices to prove the claim when $a,b \in \cN$ with $\vphi(a)=\vphi(b)=0$. Writing $x_n =
\sum_{g \in G} x_n(g) u_g$, we have
$$\|E_{\cL(G)}(a x_n b)\|_2^2 = \sum_{g \in G} |x_n(g) \vphi(a \si_g(b))|^2 \; .$$
Take $C > 0$ such that $\|x_n\| \leq C$ for all $n$. Choose $\eps > 0$. Since $(\si_g)$ is a mixing action, take $K \subset G$ finite such that
$|\vphi(a \si_g(b))|^2 \leq \eps/(2C^2)$ for all $g \in G \setminus K$. Since $x_n$ tends weakly to $0$, $x_n(g) \recht 0$ for every $g$. Hence, take
$n_0$ such that for $n \geq n_0$, $\sum_{g \in K} |x_n(g) \vphi(a \si_g(b))|^2 < \eps / 2$. Since $\sum_g |x_n(g)|^2 \leq C^2$ for all $n$, we obtain
that $\|E_{\cL(G)}(a x_n b)\|_2^2 \leq \eps$ for all $n \geq n_0$, which proves the claim.

It is then clear that any sequence of unitaries $(u_n)$ in $P$ tending weakly to $0$ satisfies the conditions of Lemma \ref{lem.intertwiners} with
$B=\cL(G)^t$ and $M=M^t$.
\end{proof}

\begin{proposition}[Popa, \cite{P1}] \label{prop.mixing-two}
Suppose that $G$ acts mixingly on the almost periodic $(\cN,\vphi)$
and arbitrarily on the almost periodic $(\cA,\psi)$. Consider the diagonal action on $\cA \ot \cN$. Write $M = \cA^\psi
\rtimes G$ and $\Mtil = (\cA \ot \cN)^{\psi \ot \vphi} \rtimes G$. Let $P \subset M$ be a diffuse subalgebra such that there is no non-zero
homomorphism from $P$ to an amplification of $\cA^\psi$. If $x \in \Mtil$ and $Px \subset \sum_{k=1}^n x_k M$, we have $x \in M$.
\end{proposition}
\begin{proof}
Write $A = \cA^\psi$. It follows from Proposition \ref{prop.translate} that there exists a sequence of unitaries $(u_n)$ in $P$ such that $\|E_A(u_n
u_g)\|_2 \recht 0$ for all $g \in G$. Let $E : (\cA \ot \cN) \rtimes G \recht \cA \rtimes G$ be the unique state-preserving conditional expectation.
By Lemma \ref{lem.intertwiners}, it suffices to check that $\|E(a u_n b)\|_2 \recht 0$ for all $a,b \in \Ker E$. It moreover suffices to check this
last statement for $a,b \in \cN$ with $\vphi(a)=\vphi(b)=0$. Writing $u_n = \sum_g u_n(g) u_g$ with $u_n(g) \in A$, we have
$$\|E(a u_n b)\|_2^2 = \sum_{g \in G} |\vphi(a \si_g(b))|^2 \, \|u_n(g)\|_2^2 \; .$$
We conclude the proof in exactly the same way as the proof of Theorem \ref{thm.mixing-one}.
\end{proof}

Finally, the notion of a $2$-mixing action is introduced. Definition
\ref{def.mixing} of a mixing action comes down to the notion of a $1$-mixing
action.

\begin{definition} \label{def.two-mixing}
An action of a countable group $G$ on $(\cA,\vphi)$ is said to be
\emph{$2$-mixing} if
$$\vphi(a \si_{g}(b) \si_h(c)) \recht \vphi(a)\vphi(b)\vphi(c)
\quad\text{when $g,h,g^{-1}h \recht \infty$.}$$
\end{definition}
Note that any $2$-mixing action is mixing and satisfies
$$|\vphi(a \si_{g}(b) \si_h(c)) - \vphi(a) \vphi(\si_{g}(b) \si_h(c))
| \recht 0 \quad\text{when $g,h \recht \infty$.}$$
Conversely, this last statement characterizes $2$-mixing actions.

\begin{lemma} \label{lem.two-mixing}
Let $(\si_g)_{g \in G}$ be a free $2$-mixing action of a countable group $G$ on $(X,\mu)$. Write $A=L^\infty(X,\mu)$. For every $\eps > 0$, there
exists a finite partition of $1$ in $A$ given by $1=q_1 + \cdots + q_n$ with $q_i$ projections in $A$ and satisfying
\begin{equation}\label{eq.tedoen}
\limsup_{g \recht \infty} \Bigl\| \sum_{k=1}^n \si_g(q_k) x \si_g(q_k) \Bigr\|_2^2 \leq \eps \|x\|_2^2
\end{equation}
for all $x \in A \rtimes G$ with $E_A(x) = 0$.
\end{lemma}
\begin{proof}
Choose $\eps > 0$. Combining freeness and the mixing property, we take a finite partition of $1$ in $A$ given by $1=q_1 + \cdots + q_n$ with $q_i$ projections
in $A$ and satisfying
$$\sum_{k=1}^n \tau(q_k \si_g(q_k)) \leq \eps$$
for all $g \neq e$. We claim that \eqref{eq.tedoen} holds for all $x \in A \rtimes G$ with $E_A(x) = 0$. It is sufficient to check this for $x =
\sum_{h \in F} a_h u_h$ for some finite subset $F \subset G$ not containing $e$. Then,
$$\Bigl\| \sum_{k=1}^n \si_g(q_k) x \si_g(q_k) \Bigr\|_2^2 = \sum_{h \in F,k=1}^n \tau(a_h^* a_h \si_g(q_k) \si_{hg}(q_k)) \; .$$
When $g \recht \infty$, the right hand side is arbitrary close to
$$\sum_{h \in F, k=1}^n \tau(a_h^* a_h) \tau(\si_g(q_k) \si_{hg}(q_k)) = \sum_{h \in F, k=1}^n \tau(a_h^* a_h) \tau(q_k
\si_{g^{-1}hg}(q_k)) \leq \eps \|x\|^2 \; .$$ So, we are done.
\end{proof}


\end{document}